\newif\ifejs
\ejsfalse

\newif\ifcolortwonew
\colortwonewtrue
\colortwonewfalse

\ifejs
\documentclass[ejs]{imsart}
\else
\documentclass[generic]{imsart}
\fi

\RequirePackage[OT1]{fontenc}
\RequirePackage{amsthm,amsmath,amssymb,mathrsfs,dsfont,stmaryrd,array} 
\RequirePackage[authoryear]{natbib}   
\usepackage{caption} 
\usepackage{graphicx} 
\usepackage{a4wide}
\usepackage{subcaption}
\usepackage{bm} 
\usepackage{color} 
\RequirePackage[colorlinks,citecolor=blue,urlcolor=blue]{hyperref}
\usepackage[capitalise,nameinlink
    ,noabbrev]{cleveref}
\crefname{assumption}{}{}

\ifejs
\doi{10.1214/154957804100000000}
\pubyear{0000}
\volume{0}
\firstpage{0}
\lastpage{0}
\fi

\startlocaldefs

\numberwithin{equation}{section}
\theoremstyle{plain}
\newtheorem{theorem}{Theorem}[section]
\newtheorem{corollary}{Corollary}[section]
\newtheorem{lemma}{Lemma}[section]
\newtheorem{heuristic}{Heuristic}[section]
\theoremstyle{definition}
\newtheorem{assumption}{Assumption}[section]
\newtheorem{definition}{Definition}[section]
\theoremstyle{remark}
\newtheorem{remark}{Remark}[section]
\DeclareMathOperator*{\argmax}{arg\,max}
\DeclareMathOperator{\Pow}{Pow}
\DeclareMathOperator{\FDR}{FDR}
\DeclareMathOperator{\FDP}{FDP}
\DeclareMathOperator{\MWBH}{MWBH}
\DeclareMathOperator{\GMWBH}{GMWBH}
\DeclareMathOperator{\BH}{BH}
\DeclareMathOperator{\ABH}{ABH}
\DeclareMathOperator{\WBH}{WBH}
\DeclareMathOperator{\HZZ}{HZZ}
\DeclareMathOperator{\ZZPro}{Pro}
\DeclareMathOperator{\DiffPow}{DiffPow}
\DeclareMathOperator{\IHW}{IHW}
\DeclareMathOperator{\ADDOW}{ADDOW}
\DeclareMathOperator{\GADDOW}{GADDOW}
\DeclareMathOperator{\crossADDOW}{crADDOW}
\DeclareMathOperator{\cross}{cross}

\DeclareMathOperator{\LCM}{LCM}
\newcommand{\Pb}{\mathbb{P}}
\newcommand{\comp}[1]{{#1}^{\mathsf{c}}} 
\newcommand{\Pro}[1]{\mathbb{P}\left(#1\right)} 
\newcommand{\Esp}[1]{\mathbb{E}\left[ #1 \right]}
\newcommand{\ind}[1]{\mathds{1}_{\left\{#1 \right\}}}
\newcommand{\Km}{\hat K}
\ifcolortwonew
\newenvironment{version2}{\hypersetup{allcolors=red}\color{red}}{\hypersetup{allcolors=blue}}
\else
\newenvironment{version2}{}{}
\fi

\ifcolortwonew
\newenvironment{version3}{\hypersetup{allcolors=red}\color{red}}{\hypersetup{allcolors=blue}}
\else
\newenvironment{version3}{}{}
\fi

\endlocaldefs

\begin{document}

\begin{frontmatter}

\title{
Adaptive $p$-value weighting with power optimality
\ifejs
\thanksref{t1}
\fi
}
\ifejs
\thankstext{t1}{This is an original survey paper}
\fi
\runtitle{Adaptive $p$-value weighting with power optimality}
\author{\fnms{Guillermo} \snm{Durand}\corref{}\ead[label=e1]{guillermo.durand@upmc.fr}}
\address{Laboratoire de Probabilit\'{e}s, Statistique et Mod\'{e}lisation,\\
Sorbonne Universit\'{e},\\
4 place Jussieu, 75252 Paris Cedex 05\\
\printead{e1}}
\runauthor{G. Durand}

\begin{abstract}
Weighting the $p$-values is a well-established strategy that improves the power of multiple testing procedures while dealing with heterogeneous data. However, how to achieve this task in an optimal way is rarely considered in the literature. This paper contributes to fill the gap in the case of group-structured null hypotheses, by introducing a new class of procedures named ADDOW (for Adaptive Data Driven Optimal Weighting) that adapts both to the alternative distribution and to the proportion of true null hypotheses. We prove the asymptotical FDR control and power optimality among all weighted procedures of ADDOW, which shows that it dominates all existing procedures in that framework. Some numerical experiments show that the proposed method preserves its optimal properties in the finite sample setting when the number of tests is moderately large.
\end{abstract}



\tableofcontents

\end{frontmatter}
\section{Introduction}

Recent high-throughput technologies bring to the statistical community new type of data being increasingly large, heterogeneous and complex. Addressing significance in such context is particularly challenging because of the number of questions that could naturally come up. A popular statistical method is to adjust for multiplicity by controlling the False Discovery Rate (FDR), which is defined as the expected proportion of errors among the items declared as significant. Once the amount of possible false discoveries is controlled, the question of increasing the power, that is the amount of true discoveries, arises naturally. In the literature, it is well-known that the power can be increased by clustering the null hypotheses into homogeneous groups. The latter can be derived in several ways:

\begin{itemize}

\item sample size: a first example is the well-studied data set of the Adequate Yearly Progress (AYP) study \citep{rogosa2005accuracy}, which compares the results in mathematics tests between socioeconomically advantaged and disadvantaged students in Californian high school. As studied by \citet{cai2009simultaneous}, ignoring the sizes of the schools tends to favor large schools among the detections, simply because large schools have more students and not because the effect is stronger. By grouping the schools in small, medium, and large schools, more rejections are allowed among the small schools, which increases the overall detection capability. This phenomenon also appears in more large-scale studies, as in GWAS (Genome-Wide Association Studies) by grouping hypotheses according to allelic frequencies, \citep{sun2006stratified} 
 or in microarrays experiments by grouping the genes according to the DNA copy number status \citep{roquain2009optimal}. \begin{version3} Common practice is generally used to build the groups from this type of covariate. 
\end{version3}

\item spatial structure: some data sets naturally involve a spatial (or temporal) structure into groups. A typical example is neuroimaging: in \citet{schwartzman2005cross}, a study compares diffusion tensor imaging brain scans on 15443 voxels of 6 normal and 6 dyslexic children. By estimating the densities under the null of the voxels of the front and back halves of the brain, some authors highlight a noteworthy difference which suggests that analysing the data by making two groups of hypotheses seems more appropriate, see \citet{efron2008simultaneous} and \citet{cai2009simultaneous}.

\item hierarchical relation: groups can be derived from previous knowledge on hierarchical structure, like pathways for genetic studies, based for example on known ontologies (see e.g. \citet{ashburner2000gene}). Similarly, in clinical trials, the tests are usually grouped in primary and secondary endpoints, see \citet{dmitrienko2003gatekeeping}. 

\end{itemize}

In these examples, while ignoring the group structure can lead to overly conservative procedures, this knowledge can easily be incorporated by using weights. This method can be traced back to \citet{holm1979simple} who presented a sequentially rejective Bonferroni procedure that controls the Family-Wise Error Rate (FWER) and added weights to the $p$-values. Weights can also be added to the type-I error criterion instead of the $p$-values, as presented in \citet{benjamini1997multiple} with the so-called weighted FDR. \citet{blanchard2008two} generalized the two approaches by weighting the $p$-values and the criterion, with a finite positive measure to weigh the criterion (see also \citet{ramdas2017unified} for recent further generalizations). \citet{genovese2006false} introduced the $p$-value weighted BH procedure (WBH) which has been extensively used afterwards with different choices for the weights. \citet{roeder2006using,roeder2009genome} have built the weights upon genomic linkage, to favor regions of the genome with strong linkage. 
\citet{hu2010false} calibrated the weights by estimating the proportion of true nulls inside each group (procedure named HZZ here). \citet{zhao2014weighted} went one step further by improving HZZ and BH with weights that maximize the number of rejections at a threshold  computed from HZZ and BH. They proposed two procedures Pro1 and Pro2 shown to control the FDR asymptotically and to have a better power than BH and HZZ. 

However, the problem of finding optimal weights (in the sense of achieving maximal averaged number of rejected false nulls) has been only scarcely considered in the literature. For FWER control and Gaussian test statistics, \citet{wasserman2006weighted} designed oracle and data-driven optimal weights, while \citet{dobriban2015optimal} considered a Gaussian prior on the signal. For FDR control, \citet{roquain2009optimal} \begin{version2}and \citet{habiger2014weighted}\end{version2} designed oracle optimal weights by using the knowledge of the distribution under the alternative of the hypotheses. Unfortunately, this knowledge is not reachable in practice. This leads to the natural idea of estimating the oracle optimal weights by maximizing the number of rejections. This idea has been followed by \citet{ignatiadis2016data} with a procedure called IHW. While they proved that IHW controls \begin{version3}asymptotically\end{version3} the FDR, its power properties have not been considered. In particular, it is unclear whether maximizing the overall number of rejections is appropriate in order to maximize power. \begin{version3}Other recent works \citep{li2016multiple,ignatiadis2017covariate,lei2018adapt} suggest weighting methods (with additional steps or different threshold computing rules) but they don't address the power question theoretically either.\end{version3}

In this paper, we present a general solution to the problem of optimal data-driven weighting of BH procedure in the case of grouped null hypotheses. The new class of procedures is called ADDOW (for Adaptive Data-Driven Optimal Weighting). \begin{version3}It relies on the computation of weights that maximize the number of detections at any rejection threshold, combined with the application of a step-up procedure with those weights. This is similar to IHW, however, by taking a larger weight space thanks to the use of estimators of true null proportion in each group, we allow for larger weights, hence more detections.\end{version3} With mild assumptions, we show that ADDOW asymptotically controls the FDR and has optimal power among all weighted step-up procedures. Interestingly, our study shows that the heterogeneity with respect to the proportion of true nulls should be taken into account in order to attain optimality. This fact seems to have been ignored so far: for instance we show that IHW \begin{version3}has optimality properties\end{version3} when the true nulls are evenly distributed across groups but \begin{version3}we also show that \end{version3}its performance can quickly deteriorate otherwise \begin{version3}with a numerical counterexample\end{version3}.

In Section~\ref{section_framework}, we present the mathematical model and assumptions. In Section~\ref{section_weighting}, we define what is a weighting step-up procedure \begin{version3}and discuss some procedures of the literature\end{version3}. In Section~\ref{section_addow}, we introduce ADDOW. 
Section~\ref{section_results} provides our main theoretical results. Our numerical simulations are presented in Section~\ref{section_simulations}, while \begin{version3}the overfitting problem is discussed in Section~\ref{section_overfitting} with the introduction of a variant of ADDOW. We\end{version3} conclude in Section~\ref{section_conclusion} with a discussion. The proofs of the two main theorems are given in Section~\ref{section_proof_thm} and more technical results are deferred to appendix. Let us underline that an effort has been made to make the proofs as short and concise as possible, while keeping them as clear as possible.

In all the paper, the probabilistic space is denoted $\left(\Omega,\mathcal{A},\mathbb{P}\right)$. The notations $\overset{a.s.}{\longrightarrow}$ and $\overset{\mathbb{P}}{\longrightarrow}$ stand for the convergence almost surely and in probability. 

\section{Setting}
\label{section_framework}

\subsection{Model}
\label{subsection_model}

We consider the following stylized grouped $p$-value modeling: let $G\geq2$ be the number of groups. 
\begin{version3}
Let us emphasize that $G$ is kept fixed throughout the paper. Because our study will be asymptotic in the number of tests $m$, for each $m$ we assume that we test $m_g$ hypotheses in group $g\in\{1,\dotsc,G\}$, where the $m_g$ are non-decreasing integer sequences depending on $m$ (the dependence is not written for conciseness) and such that $\sum_{g=1}^G m_g=m$. In each group $g\in\{1,\dotsc,G\}$, let $\left(H_{g,1}, \dotsc,H_{g,m_g}\right)$ be some binary variables corresponding to the null hypotheses to be tested in this group, with $H_{g,i}=0$ if it is true and $H_{g,i}=1$ otherwise. Consider in addition $\left(p_{g,1},\dotsc,p_{g,m_g}\right)$ some random variables in $[0,1]$ where each $p_{g,i}$ corresponds to the $p$-value testing $H_{g,i}$. Note also $m_{g,1}=\sum_{i=1}^{m_g}H_{g,i} $ the number of false nulls and $m_{g,0}=m_g-m_{g,1}$ the number of true nulls in group $g$. 
\end{version3}

We make the following marginal distributional assumptions for $p_{g,i}$. 
\begin{version3}
\begin{assumption}
If $H_{g,i}=0$, $p_{g,i}$ follows a uniform distribution on $[0,1]$.\label{ass1}
\end{assumption}
\end{version3}
We denote by $U:x\mapsto \ind{x>0} \times \min(x,1)$ its cumulative distribution function (c.d.f.). 
\begin{version3}
\begin{assumption}
If $H_{g,i}=1$, $p_{g,i}$ follows a common distribution corresponding to c.d.f. $F_g$, which is strictly concave on $[0,1]$.\label{ass2}
\end{assumption}
\end{version3}
In particular, note that the $p$-values are assumed to have the same alternative distribution within each group. \begin{version3}Note that the concavity assumption is mild (and implies continuity on $\mathbb{R}$ as proven in Lemma~\ref{Fgcont} for completeness).\end{version3} Furthermore, by concavity, $x\mapsto \frac{F_g(x)-F_g(0)}{x-0}$ has a right limit in 0 that we denote by $f_g(0^+)\in[0,\infty]$, and $x\mapsto \frac{F_g(x)-F_g(1)}{x-1}$ has a left limit in 1 that we denote by $f_g(1^-)\in[0,\infty)$.

\begin{version3}
\begin{assumption}
There exists $\pi_g>0$ and $\pi_{g,0}>0$ such that for all $g$, $m_g/m\to\pi_g$ and $m_{g,0}/m_g\to\pi_{g,0}$ when $m\to\infty$. Additionally, for each $g$, $ \pi_{g,1}=1- \pi_{g,0}>0$. \label{ass3}
\end{assumption}
\end{version3}
\begin{version3}The above assumption means\end{version3} that, asymptotically, no group, and no proportion of signal or sparsity, is vanishing. We denote $\pi_0=\sum_g \pi_g \pi_{g,0}$ the mean of the $\pi_{g,0}$'s and denote the particular case where the nulls are evenly distributed in each group by~\eqref{ED}:
\begin{equation}
\pi_{g,0}=\pi_0 , \:\:1\leq g\leq G.
\label{ED}\tag{ED}
\end{equation}

Let us \begin{version3}finally\end{version3} specify assumptions on the joint distribution of the $p$-values. \begin{version3}
\begin{assumption}
The $p$-values are weakly dependent within each group:
\begin{equation}
\frac{1}{m_{g,0}} \sum_{i=1}^{m_g}\ind{p_{g,i}\leq t, H_{g,i}=0} \overset{\Pb}{\longrightarrow} U(t) ,\:\: t \geq0,
\label{wd0}
\end{equation}
and
\begin{equation}
\frac{1}{m_{g,1}} \sum_{i=1}^{m_g}\ind{p_{g,i}\leq t, H_{g,i}=1} \overset{\Pb}{\longrightarrow} F_g(t) ,\:\: t \geq0 .
\label{wd1}
\end{equation}\label{ass4}
\end{assumption}
\end{version3}
This assumption is mild and classical, see \citet{storey2004strong}. Note that weak dependence is trivially achieved if the $p$-values are independent\begin{version3}, and that no assumption on the $p$-value dependence accross groups is made.\end{version3} \begin{version3}Finally note that there is a hidden dependence in $m$ in the joint distribution of the $p$-values $(p_{g,i})_{\substack{1\leq g\leq G  \\ 1\leq i\leq m_g}}$ but that does not impact the remaining of the paper as long as~\eqref{wd0} and~\eqref{wd1} are satisfied.\end{version3}

\subsection{$\pi_{g,0}$ estimation}
\label{subsection_estimation}

\begin{version3}
\begin{assumption}
For each $g$, we have at hand an (over-)estimator $\hat \pi_{g,0}\in(0,1]$ of $m_{g,0}/m_g$ 
such that $\hat \pi_{g,0}\overset{\mathbb{P}}{\longrightarrow} \bar \pi_{g,0}$ for some $\bar \pi_{g,0}\geq \pi_{g,0}$.\label{ass5}
\end{assumption}
\end{version3}
Let also $\bar \pi_0=\sum_g \pi_g \bar\pi_{g,0}$. \begin{version3}In the model of Section~\ref{subsection_model}\end{version3}, this assumption can be fulfilled by using the estimators introduced in \citet{storey2004strong}:
\begin{equation}
\hat \pi_{g,0}(\lambda)=\frac{1-\frac{1}{m_g}\sum_{i=1}^{m_g}\ind{p_{g,i}\leq\lambda}+\frac{1}{m} }{1-\lambda},
\label{storey}
\end{equation}
for a given parameter $\lambda\in(0,1)$ let arbitrary (the $\frac{1}{m}$ is here just to ensure $\hat \pi_{g,0}(\lambda)>0$). It is easy to deduce from~\eqref{wd0} and~\eqref{wd1} that $\frac{1}{m_g}\sum_{i=1}^{m_g}\ind{p_{g,i}\leq\lambda}\overset{\begin{version3}\mathbb{P}\end{version3}}{\longrightarrow} \pi_{g,0}\lambda+\pi_{g,1}F_g(\lambda)$, which provides our condition:
$$ \hat \pi_{g,0}(\lambda) \overset{\begin{version3}\mathbb{P}\end{version3}}{\longrightarrow} \pi_{g,0}+\pi_{g,1}\frac{1-F_g(\lambda)}{1-\lambda} \geq \pi_{g,0} .$$

While $(\bar\pi_{g,0})_g$ is let arbitrary in our setting, some particular cases will be of interest in the sequel. First is the Evenly Estimation case~\eqref{EE} one where 
\begin{equation}
\bar\pi_{g,0}=\bar\pi_0 ,\:\: 1\leq g\leq G.
\label{EE}\tag{EE}
\end{equation}
In that case, our estimators all share the same limit, and doing so they do not take in account the heterogeneity with respect to the proportion of true nulls. \begin{version3}Case\end{version3}~\eqref{EE} is relevant when the proportion of true nulls is homogeneous across groups, that is, when~\eqref{ED} holds. A particular subcase of~\eqref{EE} is the Non Estimation case~\eqref{NE} where:
\begin{equation}
\hat\pi_{g,0}=1 
\text{ which implies } \bar\pi_{g,0}=1 ,\:\: 1\leq g\leq G.
\label{NE}\tag{NE}
\end{equation}
\begin{version3}
Case~\eqref{NE} is basically the case where no estimation is intended, and the estimators are simply taken equal to 1.
\end{version3}

\begin{version2}Let us also\end{version2} introduce the Consistent Estimation case~\eqref{CE} for which the estimators $\hat \pi_{g,0}$ are assumed to be all consistent: 
\begin{equation}
\bar\pi_{g,0}=\pi_{g,0} ,\:\:1\leq g\leq G.
\label{CE}\tag{CE}
\end{equation}
While this corresponds to a favorable situation, this assumption can be met in classical situations, where $f_g(1^-)=0$ and $\lambda=\lambda_m$ tends to 1 slowly enough in definition~\eqref{storey}, see Lemma~\ref{PAindepgauss} in Section~\ref{subsection_prooflemmas}. The condition $f_g(1^-)=0$ is called "purity" in the literature. It has been introduced in \citet{genovese2004stochastic} and then deeply studied, along with the convergence of Storey estimators, in \citet{neuvial2013asymptotic}.

\begin{version2}
Finally, the main case of interest is the Multiplicative Estimation case~\eqref{ME} defined as the following:
\begin{equation}
\exists C\geq1,\,\bar \pi_{g,0}=C \pi_{g,0} ,\:\:1\leq g\leq G.
\label{ME}\tag{ME}
\end{equation}
Note that the constant $C$ above cannot depend on $g$. Interestingly, the~\eqref{ME} case covers the~\eqref{CE} case (in this respect, $C=1$) and also the case where~\eqref{ED} and~\eqref{EE} both hold (in this respect, $C=\frac{\bar\pi_0}{\pi_0}$). So the~\eqref{ME} case can be viewed as a generalization of previous cases.
\end{version2}

\subsection{Criticality}
\label{subsection_crit}

\begin{version3}
Depending on the choice of $\alpha$, multiple testing procedures may make no rejection at all when $m$ tends to $\infty$. This case is not interesting and we should focus on the other case.
\end{version3}
To this end, \citet{chi2007performance} introduced the notion of criticality: they defined some critical alpha level, denoted $\alpha^*$, for which BH procedure has no asymptotic power if $\alpha<\alpha^*$. \begin{version3}\citet{neuvial2013asymptotic} generalized this notion for any multiple testing procedure (see Section 2.5 therein) and also established a link between criticality and purity.\end{version3}

\begin{version3}
In Section~\ref{subsection_prooflemmas}, Definition~\ref{def_alpha}, we define $\alpha^*$ 
\end{version3}
in our heterogeneous setting and will focus in our results on the supercritical case. 
\begin{version3}
\begin{assumption}
The target level $\alpha$ lies in $(\alpha^*,1)$.\label{ass6}
\end{assumption}
\end{version3}
Lemma~\ref{alphacrit} states that $\alpha^*<1$ so such an $\alpha$ always exists. While the formal definition of $\alpha^*$ is reported to the appendix for the sake of clarity, let us emphasize that it depends on the \begin{version3}parameters of the model, that are\end{version3} $(F_g)_g$, $(\pi_g)_g$ and $(\pi_{g,0})_g$, and on the \begin{version3}parameters of the chosen estimators, that are\end{version3} $(\bar\pi_{g,0})_g$.

\subsection{Leading example}
\label{subsection_gaussian}
While our framework allows a general choice for $F_g$, a canonical example that we have in mind is the Gaussian one-sided framework where \begin{version3}the $p$-values are derived from Gaussian test statistics.\end{version3}

\begin{version3}
Formally, we assume that $p_{g,i}=\bar\Phi(X_{g,i})$, where $\bar\Phi(z)=\mathbb{P}\left(Z\geq z \right)$ for $Z\sim\mathcal{N}(0,1)$, and 
$$\boldsymbol{\mathrm{X}}=(X_{1,1},\dotsc,X_{1,m_1},\dotsc,X_{g,1},\dotsc,X_{g,m_g})$$
is a Gaussian vector with distribution $\mathcal{N}(\boldsymbol{\mu},\Sigma)$. Here, 
$$\boldsymbol{\mathrm{\mu}}=(\mu_{1,1},\dotsc,\mu_{1,m_1},\dotsc,\mu_{G,1},\dotsc,\mu_{G,m_G}),$$
with $\mu_{g,i}=0$ if $H_{g,i}=0$ and $\mu_{g,i}=\mu_g$ if $H_{g,i}=1$, and we assume that $\Sigma_{j,j}=1$ for all $1\geq j\geq m$. Hence $X_{g,i}\sim\mathcal{N}(0,1)$ under the null, and $X_{g,i}\sim\mathcal{N}(\mu_g,1)$ under the alternative.
\end{version3}

\begin{version3}
In this case, Assumption~\ref{ass1} is fulfilled, and 
$$F_g(\cdot)=\bar\Phi\left(\bar\Phi^{-1}(\cdot)-\mu_g \right),$$ with derivative $$f_g(\cdot) = \exp\left(\mu_g\left(\bar\Phi^{-1}(\cdot)-\frac{\mu_g}{2}   \right)   \right)>0,$$
hence $F_g$ is strictly concave and Assumption~\ref{ass2} is also fulfilled. Furthermore we easily check that  $f_g(0^+)=\infty$, so $\alpha^*=0$ and $f_g(1^-)=0$ which means that this framework is supercritical ($\alpha^*=0$, see Definition~\ref{def_alpha}) with purity and then can achieve consistent estimation~\eqref{CE} with additional independence assumptions.
\end{version3}

\begin{version3}
Two particular subcases of interest arise when $\Sigma$ has a particular form and can be written as
\begin{equation*}
\begin{pmatrix}
\Sigma_{(1)}& \boldsymbol{0} & \hdotsfor{2} &  \boldsymbol{0}\\
 \boldsymbol{0} &\Sigma_{(2)}& \boldsymbol{0} & \dots& \boldsymbol{0}\\
 \vdots&\vdots&\ddots&\vdots&\vdots\\
  \boldsymbol{0}&\dots&\boldsymbol{0}& \Sigma_{(G-1)}&  \boldsymbol{0} \\
  \boldsymbol{0} &  \hdotsfor{2} &  \boldsymbol{0} & \Sigma_{(G)}
\end{pmatrix},
\end{equation*}
where $\Sigma_{(g)}$ is a square matrix of size $m_g$. The first subcase is when $\Sigma_{(g)}$ is the identity matrix. In this case, the $p$-values are all independent and Assumption~\ref{ass4} is fulfilled by the law of strong numbers. The second subcase is when $\Sigma_{(g)}$ is a Toeplitz matrix with $\left(\Sigma_{(g)}\right)_{j,k}=\frac{1}{|j-k|+1}$. In this case, Assumption~\ref{ass4} is also fulfilled (see e.g. \citealp[Proposition 2.1, Equation (LLN-dep) and Theorem 3.1]{delattre2016empirical}).
\end{version3}



\subsection{Criterion}
\label{subsection_criterions}
The set of indices corresponding to true nulls is denoted by $\mathcal{H}_0$, that is $(g,i)\in\mathcal{H}_0$ if and only if $H_{g,i}=0$, and we also denote $\mathcal{H}_1=\comp{\mathcal{H}_0}$. 

In this paper, we define a multiple testing procedure $R$ as a set of indices that are rejected: $p_{g,i}$ is rejected if and only if $(g,i)\in R$. The False Discovery Proportion (FDP) of $R$, denoted by $\FDP(R)$, is defined as the number of false discoveries divided by the number of rejections if there are any, and 0 otherwise:
$$ \FDP(R)=\frac{\left|R \cap\mathcal{H}_0  \right|}{\left|R\right| \vee1} .$$
We denote $\FDR(R)=\mathbb{E}\left[ \FDP(R) \right]$ the FDR of $R$. Its power, denoted $\Pow(R)$, is defined as the mean number of true positives divided by $m$:
$$\Pow(R)=m^{-1}\mathbb{E}\left[\left|R \cap\mathcal{H}_1  \right|\right] .$$
Note that our power definition is slightly different than the usual one for which  the number of true discoveries is divided by $m_1=\sum_g m_{g,1}$ instead of $m$. This simplifies our expressions (see Section~\ref{subsection_notation}) and does not have any repercussion because the two definitions differ only by a multiplicative factor converging to $1-\pi_0\in(0,1)$ when $m\to\infty$. 

\begin{version3}Finally, let us emphasize that the power is the (rescaled) number of good rejections, that is, the number of rejected hypotheses that are false. The power is different from the number of total rejections, this distinction is fundamental and will be discussed all along this paper (like, for example, when discussing Heuristic~\ref{mainheur}, or in the simulations of Section~\ref{subsection_sim_me}).\end{version3}

\section{Weighting}
\label{section_weighting}

\subsection{Weighting the BH procedure}
\label{subsection_MWBH}

Say we want to control the FDR at level $\alpha$. Assume that the $p$-values are arranged in increasing order $p_{(1)}\leq \dotsc \leq p_{(m)} $ with $p_{(0)}=0$, the classic BH procedure consists in rejecting all $p_{g,i}\leq \alpha \frac{\hat k}{m}$ where $\hat k= \max\left\{k\geq0 : p_{(k)}\leq \alpha \frac{ k}{m}  \right\}$.

Take a nondecreasing function $h$ defined on $[0,1]$ such that $h(0)=0$ and $h(1)\leq 1$, we denote $\mathcal{I}(h)=\sup\left\{ u\in[0,1] : h(u)\geq u     \right\}.$ Some properties of the functional $\mathcal{I}(\cdot)$ are gathered in Lemma~\ref{IG}, in particular $h\left(\mathcal{I}(h)\right)=\mathcal{I}(h)$. We now reformulate BH with the use of $\mathcal{I}(\cdot)$, because it is more convenient when dealing with asymptotics. Doing so, we follow the formalism notably used in \citet{roquain2009optimal} and \citet{neuvial2013asymptotic}. Define the empirical function
$$\widehat G : u \mapsto m^{-1}\sum_{g=1}^G\sum_{i=1}^{m_g}\mathds{1}_{\{ p_{g,i}\leq \alpha u \}},$$
then ${\hat k}=m\times\mathcal{I}( \widehat G )$. This is a particular case of Lemma~\ref{hatuegalIG}. \begin{version3}Note that $\widehat G(u)$ is simply the number of $p$-values that are less than or equal to $\alpha u$, divided by $m$.\end{version3}

\begin{figure}
\centering\includegraphics[width=\linewidth]{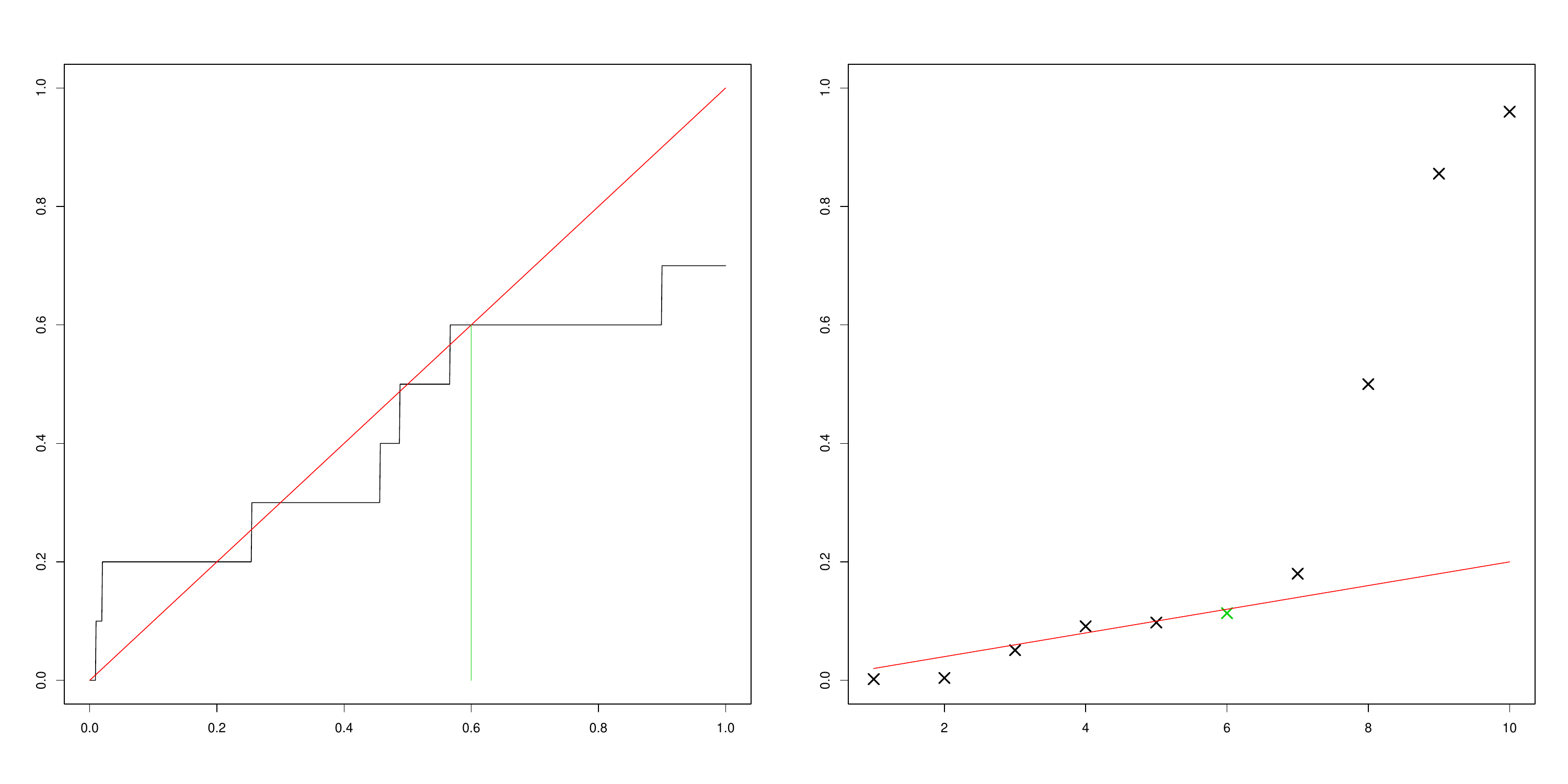}
\caption[Illustration of the BH procedure]{\footnotesize The BH procedure applied to a set of 10 $p$-values. Right plot: the $p$-values and the function $k\to\alpha k/m$. Left plot: identity function and $\widehat G$. Each plot shows that 6 $p$-values are rejected.}
\label{I_of_hat_G}
\end{figure}

The graphical representation of the two points of view for BH is depicted in Figure~\ref{I_of_hat_G} with $m=10$. The $p$-values are plotted on the right part of the figure along with the function $k\mapsto \alpha k/m$ and we see that the last $p$-value under the line is the sixth one. On the left, the function $\widehat G$ 
corresponding to these $p$-values is displayed alongside the identity function, with the last crossing point being located between the sixth and seventh jumps, thus $\mathcal{I}( \widehat G )=6/m$ and 6 $p$-values are rejected.

The weighted BH (WBH) with weight vector $w\in \mathbb{R}^G_+$ is defined by computing
$$\widehat G_w : u \mapsto m^{-1}\sum_{g=1}^G\sum_{i=1}^{m_g}\mathds{1}_{\{ p_{g,i}\leq \alpha u w_g \}}$$
and rejecting all $p_{g,i}\leq \alpha  \mathcal{I}\left(G_w\right) w_g$. We denote it $\WBH(w)$. Note that $w$ is authorized to be random, hence it can be computed from the $p$-values. In particular, $\BH=\WBH(\bm{1})$ where $\bm{1}=(1,\dots,1)\in\mathbb{R}^G_+$.

Following \citet{roquain2009optimal}, to deal with optimal weighting, we need to further generalize WBH into a multi-weighted BH (MWBH) procedure by introducing a weight function $W : [0,1]\to \mathbb{R}^G_+$, {which can be random}, such that the following function:
\begin{equation}
\widehat G_W : u \mapsto m^{-1}\sum_{g=1}^G\sum_{i=1}^{m_g}\mathds{1}_{\{ p_{g,i}\leq \alpha u W_g(u) \}},
\label{eq_def_gW}
\end{equation}
is nondecreasing. The resulting procedure rejects all the $p$-values such that $p_{g,i}\leq \alpha  \hat u_W W_g(\hat u_W)$ and is denoted $\MWBH(W)$ where, for the rest of the paper, we denote
\begin{equation}
\hat u_W=\mathcal{I}\left( \widehat G_W \right),
\label{def_uw}
\end{equation}
and name it the step-up threshold. One different weight vector $W(u)$ is associated to each $u$, hence the "multi"-weighting. Note that the class of MWBH procedures \begin{version3}is a straightforward generalization of the class of\end{version3} WBH procedures because \begin{version3}for\end{version3} any weight vector \begin{version3}$w$, $w$\end{version3} can be seen as a constant weight function \begin{version3}$u\mapsto w$ and $\widehat G_w$ is nondecreasing\end{version3}.

Note that, there is a simple way to compute $\hat u_W$. For each $r$ between 1 and $m$ denote the $W(r/m)$-weighted $p$-values $p_{g,i}^{[r]}=p_{g,i}/W_g(r/m)$ (with the convention $p_{g,i}/0=\infty$), order them $p_{(1)}^{[r]}\leq\dotsc\leq p_{(m)}^{[r]}$ and note $p_{(0)}^{[r]}=0$. Then $\hat  u_W= m^{-1}\max\left\{ r\geq0 :  p_{(r)}^{[r]} \leq \alpha \frac{r}{m} \right\}$ (this is Lemma~\ref{hatuegalIG}).

As in previous works \begin{version3}(see e.g. \citealp{genovese2006false} or \citealp{zhao2014weighted})\end{version3}, in order to achieve a valid FDR control, these procedures should be used with weights that satisfy some specific \begin{version3}constraints\end{version3}. 
\begin{version3}The following weight spaces will be used in the following of the paper:\end{version3}
\begin{equation}
\Km=\left\{w\in \mathbb{R}^G_+ : \sum_g \begin{version2}\frac{m_g}{m}\end{version2}  \hat \pi_{g,0} w_g\leq 1\right\},
\label{def_Km}
\end{equation}
\begin{equation}
\Km_{\text{NE}}=\left\{w\in \mathbb{R}^G_+ : \sum_g  \begin{version2}\frac{m_g}{m} \end{version2} w_g\leq 1\right\}.
\label{def_KmNE}
\end{equation}
Note that $\Km$ 
{may appear unusual} because it depends on the estimators $\hat\pi_{g,0}$, however it is completely known and usable in practice. \begin{version3}Some intuition about the choice of $\Km$ is given in next section.\end{version3} Note also that $\Km=\Km_{\text{NE}}$ in the~\eqref{NE} case.

Finally, for a weight function $W$ and a \begin{version3}rejection\end{version3} threshold $u\in[0,1]$, we denote by $R_{u,W}$ the double indexed procedure rejecting the $p$-values less than or equal to $\alpha u W_g(u)$, that is $R_{u,W}=\{(g,i) : p_{g,i}\leq \alpha u W_g(u) \}$. By~\eqref{eq_def_gW}, note that  $\widehat G_W(u)=m^{-1}\left|  R_{u,W}\right|$ \begin{version3}(which means that $\widehat G_W(u)$ is the number of rejections of $R_{u,W}$, divided by $m$)\end{version3} and \begin{version3}that $\MWBH(W)$ can also be written as $R_{\hat u_W,W}$\end{version3}. 

\subsection{Choosing the weights}
\label{subsection_previous}

Take $W$ and $u$, and let $P^{(m)}_W(u)=\Pow\left( R_{u,W}\right)$. We have 
\begin{align*}
P^{(m)}_W(u)&=m^{-1}\mathbb{E}\left[ \sum_{g=1}^G\sum_{i=1}^{m_g}\mathds{1}_{\{ p_{g,i}\leq \alpha u W_g(u) ,  H_{g,i}=1\} }  \right]\notag \\
&=\sum_{g=1}^G \frac{m_{g,1}}{m}F_g\left(\alpha u W_g(u)\right).
\end{align*}
Note that these relations are valid only if $W$ and $u$ are deterministic. In particular, they are not valid when used a posteriori with a data-driven weighting and $u=\hat u_W$.

In \citet{roquain2009optimal}, the authors define the oracle optimal weight function $W^{*}_{or}$ as:
\begin{equation}
W^{*}_{or}(u)=\argmax_{w\in \Km_{\text{NE}}} P^{(m)}_w(u) .
\label{def_or}
\end{equation}
Note that they defined $W^{*}_{or}$ only in case~\eqref{NE}, but their definition easily extends to the general case as above, \begin{version3}by replacing $\Km_{\text{NE}}$ by $\Km$\end{version3}. They proved the existence and uniqueness of $W^{*}_{or}$ \begin{version2}when both~\eqref{ED} and~\eqref{NE} hold\end{version2} and that, asymptotically, $\MWBH(W^{*}_{or})$ controls the FDR at level $\pi_0\alpha$ and has a better power than every $\MWBH(w^{(m)})$ for $w^{(m)}\in \Km_{\text{NE}}$ some deterministic weight vectors satisfying a convergence criterion.

However, computing $W^{*}_{or}$ requires the knowledge of the $F_g$, not available in practice, so the idea is to estimate $W^{*}_{or}$ with a data driven weight function $\widehat W^*$ and then apply MWBH with this random weight function. For this, consider the functional defined by, for any (deterministic) weight function $W$ and $u\in[0,1]$:
\begin{align}
G^{(m)}_W(u)=\mathbb{E}\left[\widehat G_W(u) \right]
 &=\sum_{g=1}^G \left(  \frac{m_{g,0}}{m}  U(\alpha u W_g(u))+ \frac{m_{g,1}}{m}  F_g(\alpha u W_g(u))   \right)\notag\\
&=P^{(m)}_W(u)+H^{(m)}_W(u),\label{eq_heur}
\end{align}
\begin{version3}
where
\begin{equation*}
H^{(m)}_W(u)=\sum_{g=1}^G  \frac{m_{g,0}}{m}  U(\alpha u W_g(u)).
\end{equation*}
\end{version3}
$G^{(m)}_W(u)$ is the mean ratio of rejections for the procedure rejecting each $p_{g,i}\leq \alpha u W_g(u)$. \begin{version3}$P^{(m)}_W(u)$ is the rescaled mean of the number of true positives (i.e. the power) of this procedure while $H^{(m)}_W(u)$ is the rescaled mean of the number of its false positives.\end{version3} 

\begin{version3}
\begin{heuristic}
Maximizing $G^{(m)}_W(u)$ should be close to maximizing $P^{(m)}_W(u)$.\label{mainheur}
\end{heuristic}
Indeed, consider weight functions $W$ such that $\sum_{g}  \frac{m_{g,0}}{m}  W_g(u)=1$ and then replace $U(x)$ by $x$ for all $x\in\mathbb{R}_+$ (whereas $U(x)=x$ only holds for $x\leq1$), then $H^{(m)}_W(u)$ becomes $\alpha u\sum_{g}  \frac{m_{g,0}}{m}  W_g(u)=\alpha u$ and it does not depend on the weights. So $P^{(m)}_W(u)$ is the only term depending on $W$ in~\eqref{eq_heur} and maximizing $P^{(m)}_W(u)$ or $G^{(m)}_W(u)$ is the same.
\end{version3} 

Now, we can evaluate the constraint \begin{version3}we just put\end{version3} on $W$ by estimating $\frac{m_{g,0}}{m}=\frac{m_{g}}{m}\frac{m_{g,0}}{m_g}$ by $\frac{m_{g}}{m}\hat\pi_{g,0}$ (which leads to the weight space $\Km$ defined in equation~\eqref{def_Km}), and $G^{(m)}_w(u)$ can be easily estimated by the (unbiased) estimator $\widehat G_w(u)$. As a result, maximizing the latter in $w$ should lead to good weights, not too far from $W^{*}_{or}(u)$.

\citet{zhao2014weighted} followed \begin{version3}Heuristic~\ref{mainheur}\end{version3} by applying a two-stage approach to derive two procedures, named Pro1 and Pro2. Precisely, in the first stage they use the weight vectors $\hat w^{(1)}=(\frac{1}{\hat \pi_0},\dotsc,\frac{1}{\hat \pi_0})$, where $\hat \pi_0=\sum_g \frac{m_g}{m}\hat \pi_{g,0}$, and $\hat w^{(2)}$ defined by $\hat w^{(2)}_g=\frac{\hat \pi_{g,1}}{\hat \pi_{g,0}(1-\hat \pi_0)}$, where $\hat \pi_{g,1}=1-\hat \pi_{g,0}$, and let $\hat u_M=\max(\hat u_{\hat w^{(1)}},\hat u_{\hat w^{(2)}})$. In the second stage, they maximize $\widehat G_w(\hat u_M)$ over $\Km$, which gives rise to the weight vector $\widehat W^*(\hat u_M)$ according to our notation. Then they define their procedures as the following: $$\ZZPro1=R_{\hat u_M,\widehat W^*(\hat u_M) },$$ 
and $$\ZZPro2=\WBH\left({\widehat W^*(\hat u_M)}   \right).$$ 
\begin{version3}$\ZZPro2$ comes from an additional step-up step compared to $\ZZPro1$, hence its rejection threshold, $\hat u_{\widehat W^*(\hat u_M)}$, is larger than $\hat u_M$ and allows for more detections.\end{version3} The caveat of this approach is that the initial thresholding, that is the definition of $\hat u_M$, seems somewhat arbitrary, which will result in sub-optimal procedures, see Corollary~\ref{cor_zz}. As a side remark, $\hat w^{(1)}$ and $\hat w^{(2)}$ are involved in other procedures of the literature. The HZZ procedure of \citet{hu2010false} is $\WBH(\hat w^{(2)})$, and $\WBH(\hat w^{(1)})$ is the classical Adaptive BH procedure (see e.g. Lemma~2 of~\citet{storey2004strong}) denoted here as ABH.

\citet{ignatiadis2016data} actually used  \begin{version3}Heuristic~\ref{mainheur}\end{version3} with multi-weighting (while their formulation differs from ours) which consists in maximizing $\widehat G_w(u)$ in $w$ for each $u$. However, their choice of the weight space is only suitable for the case~\eqref{NE} and can make \begin{version3}Heuristic~\ref{mainheur}\end{version3} break down, because in general \begin{version3}$H^{(m)}_W(u)$\end{version3} can still depend on $w$, see remark~\ref{rk_bh_best} \begin{version3}below\end{version3}. In the next section, we take the best of the two approaches to attain power optimality with data-driven weighting. Let us already mention that the crucial point is Lemma~\ref{ASTUCE}, that fully justifies  \begin{version3}Heuristic~\ref{mainheur}\end{version3}, \begin{version2}but only in case~\eqref{ME}\end{version2}. \begin{version3}When~\eqref{ME} does not hold, we must take care that Heuristic~\ref{mainheur} can fail for the same reason that it can fail with IHW. Thereby, in general, more detections do not necessarily imply more power.\end{version3}

\begin{remark}
\begin{version3}In particular, we\end{version3} can compute numerical counterexamples where BH has larger asymptotic power than IHW. For example, if we break \eqref{ED} by taking a small $\pi_{1,0}$ (almost pure signal) and a large $\pi_{2,0}$ (sparse signal), along with a small group and a large one ($\pi_1$ much smaller than $\pi_2$) and strong signal in both groups, \begin{version3}we can achieve a larger power with BH than with IHW. Our interpretation is that, in that case, IHW slightly favors group 2 because of its size, whereas the oracle optimal favors group 1 thanks to the knowledge of the true parameters.\end{version3}
BH, \begin{version3}by weighting uniformly,\end{version3} does not favor any group, \begin{version3}which allows its power to end up between the power of the oracle and the power of IHW.\end{version3} This example is \begin{version3}studied\end{version3} in Section\begin{version3}~\ref{subsection_sim_me} and illustrated in\end{version3} \begin{version3}Figures~\ref{fig_set2_pow} and~\ref{fig_set2_fdr}\end{version3}.
\label{rk_bh_best}
\end{remark}

\begin{version3}
\subsection{Recent weighting methods}
\label{subsection_sabha_etc}

Besides IHW, there are several recent methods putting weights on $p$-values. We briefly discuss three of them. The first is a variation of IHW by the same authors, IHWc \citep{ignatiadis2017covariate}, where the letter 'c' stands for 'censoring'. The method bring two innovations to IHW. First, the use of cross-weighting thanks to a subdivision of the hypotheses into folds: the weights of the $p$-values of a fold are computed by only using the $p$-values of the other folds. This approach reduces overfitting since, during the step-up procedure, the information brought by a given $p$-value is used only once instead of twice. The second innovation is the censoring, where a threshold $\tau$ is fixed and only $p$-values larger than $\tau$ are used to compute the weights, while only $p$-values lesser than $\tau$ can be rejected during the step-up. Together, these innovations allow IHWc to control the FDR in finite sample at level $\alpha$ if the $p$-values associated to true nulls are independent. However, using only large $p$-values to compute the weights seems somehow counterintuitive: large $p$-values are likely to be associated to true nulls and to be uniform, so they won't allow the weights to properly discriminate the groups and to increase the power compared to BH. We will verify this intuition in Section~\ref{subsection_sim_pow}. Finally, it is worth noting that IHWc allows for a kind of $\pi_{g,0}$ estimation \`{a} la Storey, with a variant called IHWc-Storey.

The censoring idea originates from the Structure Adaptive BH Algorithm (SABHA, \citealp{li2016multiple}), which has a group structured version with an FDR bounded by $\alpha C$ for a known constant $C>1$ when the $p$-values are independent. Hence, applying the group structured SABHA at level $\alpha/C$ gives FDR control at level $\alpha$, but using a target level $<\alpha$ can induce conservatism, especially since computing the weights only with the large $p$-values involve the same risks that we highlighted when discussing of IHWc.

Lastly, AdaPT \citep{lei2018adapt} introduces threshold surfaces $s_t(x)$ that can be considered as weights and adapted to group setting. AdaPT is not a WBH procedure, its whole philosophy is totally different and relies on symmetry properties of the true null distribution of the $p$-values by using an estimator of the FDP, different than the one implicitly used in BH-like methods, which also relies on symmetry and allow to mask $p$-values during the procedure (see also \citealp{barber2015controlling} and \citealp{arias2017distribution} for more details on this pioneering paradigm). We won't further consider AdaPT because of its fundamental differences with WBH procedures and because we are mainly interested by optimality among said WBH procedures.

For more discussion about IHW, IHWc, SABHA and AdaPT, see \citet[Section 6.2]{ignatiadis2017covariate} and \citet[Section 1.4]{lei2018adapt} .

\end{version3}

\section{New procedure: ADDOW}
\label{section_addow}


We exploit  \begin{version3}Heuristic~\ref{mainheur}\end{version3} and propose to estimate the oracle optimal weights $W^{*}_{or}$ by maximizing in $w\in \Km$ the empirical counterpart to $G^{(m)}_{w}(u)$, that is $\widehat G_w(u)$. 
\begin{definition}
We call an adaptive data-driven optimal weight function a random function $\widehat W^*:[0,1]\to \Km$ such that for all $u\in[0,1]$:
$$\widehat G_{\widehat W^*}(u)= \underset{w\in \Km}{\max}\:  \widehat G_{ w}(u) .$$
\end{definition}
\begin{version3}
Such maximum is guaranteed to exist because $\left\{ \widehat G_w(u),\, w \in \Km \right\}$ is a finite set. Indeed, it is a subset of $ \left\{ \frac{k}{m}, \,k\in \llbracket 0, m\rrbracket  \right\} $. However, for a given $u$, $\widehat W^*(u)$ may not be uniquely defined, hence there is no unique optimal weight function $\widehat W^*$ in general.
\end{version3}
So, in all the following, we \begin{version3}fix\end{version3} a certain $\widehat W^*$, and our results do not depend on the choice of $\widehat W^*$. An important fact is that $\widehat G_{\widehat W^*}$ is nondecreasing (see Lemma~\ref{croissance}) so $\hat u_{\widehat W^*}$ exists and the corresponding MWBH procedure is well-defined:

\begin{definition}
The ADDOW procedure is the MWBH procedure using $\widehat W^*$ as the weight function, that is, $\ADDOW=\MWBH\left( \widehat W^* \right)$.
\end{definition}
One shall note that ADDOW is in fact a class of procedures depending on the estimators $\hat\pi_{g,0}$ through $\Km$. \begin{version3}Its rationale is similar to IHW in that we intend to maximize the number of rejections, but incorporating the estimators $\hat\pi_{g,0}$ allows for larger weights and more detections.\end{version3} \begin{version3}Finally,\end{version3} note that, in the \eqref{NE} case, ADDOW reduces to IHW. 

\begin{remark}
It turns out that ADDOW \begin{version3}is equal\end{version3} to a certain WBH procedure. It comes from part 2 of the proof of Theorem~\ref{thm_pow} and Remark~\ref{MWBH<wbh}. Moreover, to every MWBH procedure, corresponds a WBH procedure with power higher or equal. \begin{version3}This fact does not limit the interest of the MWBH class, because computing the dominating WBH procedure of a given $\MWBH(\widehat W)$ procedure requires the knowledge of the step-up threshold $\hat u_{\widehat W}$ which is known by actually computing $\MWBH(\widehat W)$.\end{version3}
\end{remark}

\section{Results}
\label{section_results}
\subsection{Main results}
\label{subsection_addow_results}

Now we present the two main theorems of this paper. The two are asymptotical and justify the use of $\ADDOW$ when $m$ is large. The first is the control of the FDR at level at most $ \alpha$. The second shows that ADDOW has maximum power over all MWBH procedures in \begin{version2}the~\eqref{ME} case.\end{version2} The two are proven in Section~\ref{section_proof_thm}. 

\begin{theorem} Let us \begin{version3}assume that Assumptions \cref{ass1,ass2,ass3,ass4,ass5,ass6} are fulfilled\end{version3}. We have
\begin{equation}
\lim_{m\to\infty} \FDR\left( \ADDOW \right) \leq\alpha    .
\label{eq_thm_fdr}
\end{equation}
\begin{version2}Moreover, if $\alpha\leq\bar \pi_0$ and \eqref{ME} holds, 
\begin{equation}
\lim_{m\to\infty} \FDR\left( \ADDOW \right) =\frac{\alpha}{C} .    
\label{eq_thm_fdr_ref}
\end{equation}
\end{version2}
\label{thm_fdr}
\end{theorem}
\begin{version2}
\begin{remark}
Equation~\eqref{eq_thm_fdr_ref} means that in the~\eqref{CE} case (where $C=1$), exact asymptotic control is achieved.
\end{remark}
\end{version2}

\begin{theorem}  Let us \begin{version3}assume that Assumptions \cref{ass1,ass2,ass3,ass4,ass5,ass6} are fulfilled\end{version3}, with the additional assumption \begin{version2}that~\eqref{ME} holds.\end{version2} For any sequence of random weight functions $(\widehat W)_{m\geq1}$, such that $\widehat W:[0,1]\to \Km$ and $\widehat G_{\widehat W}$ is nondecreasing, we have
$$\lim_{m\to\infty} \Pow\left( \ADDOW \right)\geq\limsup_{m\to\infty} \Pow\left( \MWBH\left(\widehat W\right) \right).$$
\label{thm_pow}
\end{theorem}

\subsection{Relation to $\IHW$}
\label{subsection_ihw}

Recall that IHW reduces ADDOW in the~\eqref{NE} case, that~\eqref{NE} is a subcase of~\eqref{EE}\begin{version2}, and that when both~\eqref{EE} and~\eqref{ED} hold then~\eqref{ME} is achieved\end{version2}. Hence, as a byproduct, we deduce from Theorems~\ref{thm_fdr} and~\ref{thm_pow} the following result on IHW.

\begin{corollary}
Let us \begin{version3}assume that Assumptions \cref{ass1,ass2,ass3,ass4,ass5,ass6} are fulfilled\end{version3}, with the additional assumption \begin{version2}that~\eqref{ED} holds\end{version2}. Then
\begin{equation}
\lim_{m\to\infty} \FDR\left( \IHW \right) ={\pi_0}\alpha ,\label{eq:cor:fdr}
\end{equation}
and for any sequence of random weight functions $(\widehat W)_{m\geq1}$ such that  $\widehat W:[0,1]\to \Km_{\text{NE}}$ and $\widehat G_{\widehat W}$ is nondecreasing, we have
\begin{equation}
\lim_{m\to\infty} \Pow\left( \IHW \right) \geq  \limsup_{m\to\infty} \Pow\left( \MWBH\left(\widehat W\right) \right)  .\label{eq:cor:pow}
\end{equation}
\label{cor_ihw}
\end{corollary}

While equation~\eqref{eq_thm_fdr} of Theorem~\ref{thm_fdr} covers Theorem 4 of the supplementary material of \citet{ignatiadis2016data} (with slightly stronger assumption on the smoothness of the $F_g$’s), the FDR controlling result of Corollary~\ref{cor_ihw} gives a slightly sharper bound ($\pi_0\alpha$ instead of $\alpha)$ in~\eqref{ED} case. 

The power optimality stated in Corollary~\ref{cor_ihw} is new and was not shown in  \citet{ignatiadis2016data}. It thus supports the fact that IHW should be used under the assumption~\eqref{ED} and when $\pi_0$ is close to 1 or not estimated. 

\subsection{Comparison to other existing procedures}
\label{subsection_comparison}

For any estimators $\hat\pi_{g,0}\in[0,1]$, any weighting satisfying $\sum_g \frac{m_g}{m} w_g\leq1$ also satisfies $ \sum_g \frac{m_g}{m}\hat\pi_{g,0}w_g\leq1$, that is $\Km_{\text{NE}}\subset \Km$. Hence, any MWBH procedure estimating $\frac{m_{g,0}}{m_g}$ by 1 uses a weight function valued in $\Km$. This immediately yields the following corollary.

\begin{corollary}
Let us \begin{version3}assume that Assumptions \cref{ass1,ass2,ass3,ass4,ass5,ass6} are fulfilled\end{version3}, with the additional assumption \begin{version2}that~\eqref{ME} holds\end{version2}. Then
$$\lim_{m\to\infty} \Pow\left( \ADDOW \right)\geq\limsup_{m\to\infty} \Pow\left( R \right) ,  $$
for any $R\in\{\BH, \IHW \}$.
\end{corollary}

The next corollary simply states that ADDOW outperforms many procedures of the "weighting with $\pi_0$ adaptation" literature.

\begin{corollary}
Let us \begin{version3}assume that Assumptions \cref{ass1,ass2,ass3,ass4,ass5,ass6} are fulfilled\end{version3}, with the additional assumption \begin{version2}that~\eqref{ME} holds\end{version2}. Then
\begin{equation*}
\lim_{m\to\infty} \Pow\left( \ADDOW \right) \geq \limsup_{m\to\infty} \Pow\left(R \right),
\end{equation*}
for any $R\in\{\ZZPro1, \ZZPro2, \HZZ, \ABH \}$.
\label{cor_zz}
\end{corollary}
The results for Pro2, HZZ and ABH follow directly from Theorem~\ref{thm_pow} because these are MWBH procedures. The proof for Pro1 (which is not of the MWBH type) can be found in Section~\ref{subsection_proof_pro1}.

\section{Numerical experiments}
\label{section_simulations}

\subsection{Simulation setting}
\label{subsection_sim_setting}

\begin{version3}
FDR analysis and power analysis from Sections~\ref{subsection_sim_fdr} and~\ref{subsection_sim_pow} are conducted using simulations which setting we describe here. Section~\ref{subsection_sim_me} presents a counter-example using its own setting.
\end{version3}

We consider the one-sided Gaussian framework described in Section~\ref{subsection_gaussian} for $G=2$ groups. \begin{version3}We set $\alpha=0.05$, $m_1=m_2=4000$ (hence $m=8000$), $m_{1,0}=2800$ and $m_{2,0}=3200$, such that $\pi_{1,0}=0.7$ and $\pi_{2,0}=0.8$. The values of $\mu_1$ and $\mu_2$ are defined according to a varying parameter $\bar\mu$, which values are in $\{0.1,0.5,0.75,1,1.25,1.5,1.75,2,2.25,2.5,2.75,3\}$.\end{version3}

Our experiments have been performed by using the \begin{version3}four\end{version3} following scenarios Each simulation of each scenario is replicated 1000 times.
\begin{version3}

\begin{itemize}
\item Scenario 1: $\mu_1=\bar \mu$ and $\mu_2=2\bar \mu$ and the $p$-values are independent.
\item Scenario 2: $\mu_1=\bar \mu$ and $\mu_2=2\bar \mu$ and the dependence follows the Toeplitz pattern described in the end of Section~\ref{subsection_gaussian}.
\item Scenario 3: $\mu_1=\bar \mu$ and $\mu_2=0.01$ and the $p$-values are independent.
\item Scenario 4: $\mu_1=\bar \mu$ and $\mu_2=0.01$ and the dependence follows the Toeplitz pattern described in the end of Section~\ref{subsection_gaussian}.
\end{itemize}

\end{version3}

\begin{version3}In each scenario\end{version3}, \begin{version3}three\end{version3} groups of procedures are compared. \begin{version3}The difference between the three groups lies in the way $\pi_0$ is estimated. Group 1 corresponds to the~\eqref{NE} case: $\hat \pi_{g,0}=1$. Group 2 corresponds to the~\eqref{CE} case, with an oracle estimator: $\hat \pi_{g,0}=\pi_{g,0}$. Groups 3 use the Storey estimator  $\hat \pi_{g,0}(1/2)$ defined in Equation~\eqref{storey}. We choose $\lambda=1/2$ as it is a standard value (see e.g. \citealp{storey2002direct}).
The compared procedures are the following:
\begin{itemize}
\item ABH as defined in section~\ref{subsection_previous} (which is BH in Group 1),
\item HZZ as defined in section~\ref{subsection_previous} (except in Group 1 where it is not defined),
\item Pro2 as defined in section~\ref{subsection_previous} (for Group 1, we only use the BH threshold),
\item ADDOW (which is equal to IHW in Group 1),
\item An oracle ADDOW wich is the MWBH procedure using the oracle weights $W^{*}_{or}$ given by equation~\ref{def_or} (only in Groups 1 and 2),
\item IHWc (only in Groups 1 and 3). The version of IHWc used in Group 3 is IHWc-Storey.
\end{itemize}
For IHWc, the censoring level chosen is the default of the IHW R package, that is $\alpha$. 

In the following, only plots of scenarios 1 and 3 are shown, as the situation with Toeplitz dependence is found to be similar to the independent case, up to a slight increase of the FDR of most of the procedures.

\end{version3}


\subsection{FDR control}
\label{subsection_sim_fdr}

\begin{version3}

The FDR of all above procedures are compared in Figure~\ref{fig_mult_fdr} and Figure~\ref{fig_LARGEsmall_fdr}.

\begin{figure}
\centering
\includegraphics[width=0.99\linewidth]{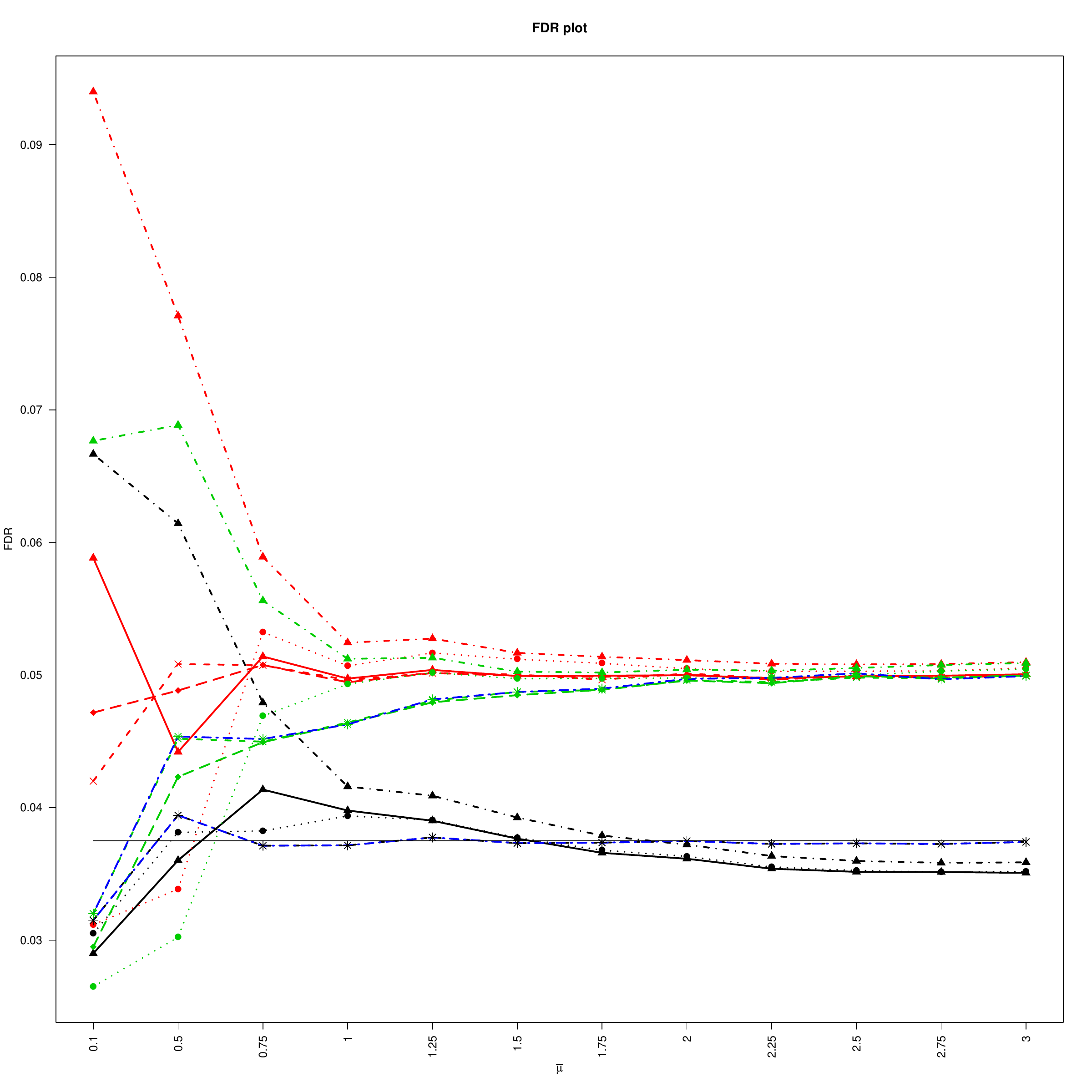}
\caption[FDR against \texorpdfstring{$\bar\mu$}{bar mu} in scenario 1]{\footnotesize \begin{version3} FDR against $\bar\mu$ in scenario 1. Group 1 in black; Group 2 in red; Group 3 in green. The type of procedure depends on the shape: Oracle ADDOW (triangles and solid line); ADDOW (triangles and dashed line); Pro2 (disks); HZZ (diamonds) and finally BH/ABH (crosses). IHWc and IHWc-Storey are in blue, respectively with black and green points. Horizontal lines: $\alpha$ and $\pi_0\alpha$ levels.  See Section~\ref{subsection_sim_setting}.\end{version3}}
\label{fig_mult_fdr}
\end{figure}

\begin{figure}
\centering
\includegraphics[width=0.99\linewidth]{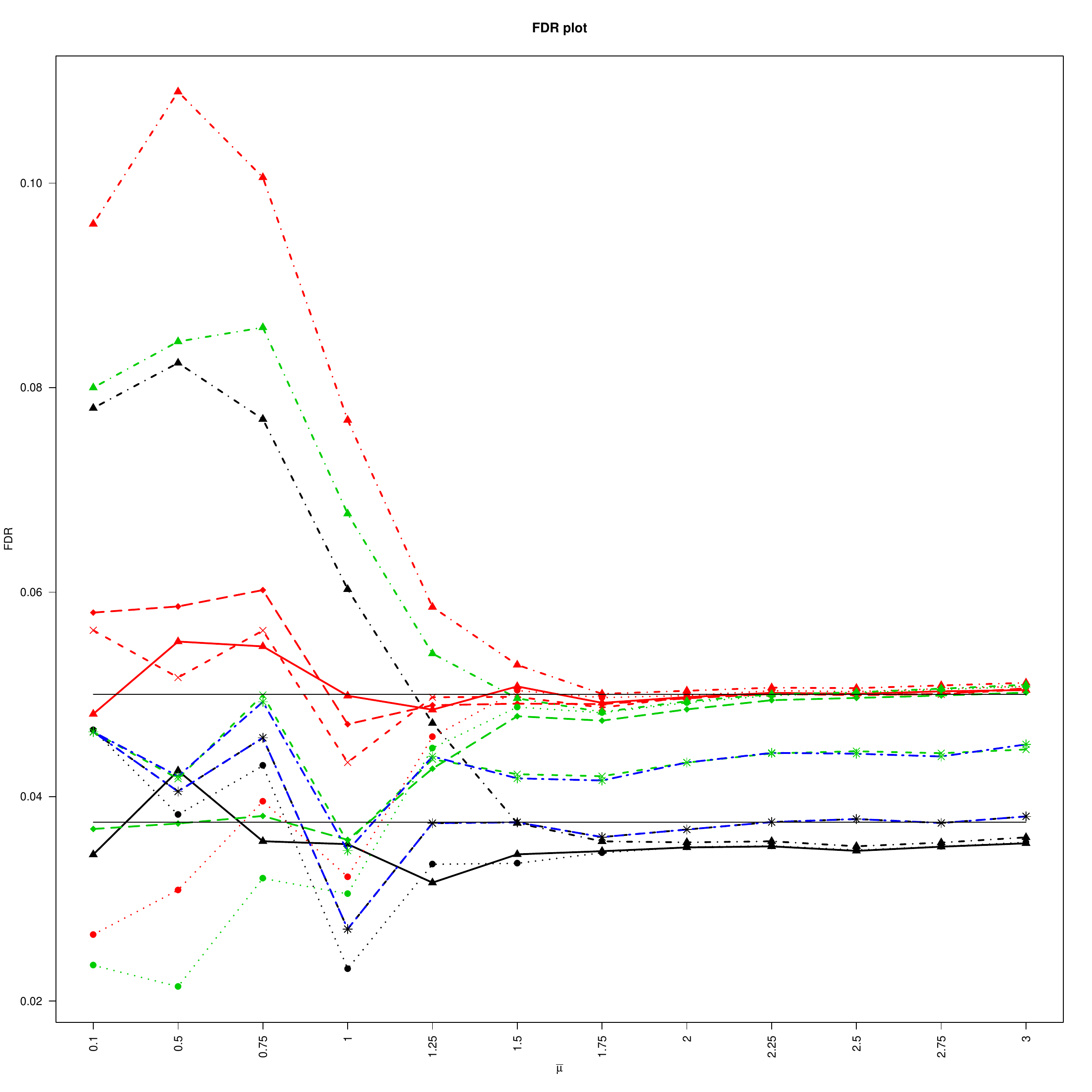}
\caption[FDR against \texorpdfstring{$\bar\mu$}{bar mu} in scenario 3]{\footnotesize \begin{version3} FDR against $\bar\mu$ in scenario 3. Same legend as in Figure~\ref{fig_mult_fdr}.\end{version3}}
\label{fig_LARGEsmall_fdr}
\end{figure}


In scenario 1, we can distinguish two different regimes depending on the signal strength. For $\mu \geq 1$ the signal strength is not weak in both groups (from $\mu_1=\bar\mu$ and $\mu_2=2\bar\mu$) and the FDR is controlled at level $\alpha$ for all procedures of Groups 2 \& 3 except ADDOW and Pro2, the two procedures using the data-driven weights, that is $\widehat W^*$. In particular, Oracle ADDOW in Group 2 controls the FDR at level $\alpha$. As the data driven weights converge to the oracle weights (see Lemma~\ref{chapeauchapeau}), we get an illustration of Theorem~\ref{thm_fdr} in the~\eqref{CE} case. The situation is similar for Group 1 and level $\pi_0\alpha$, except for Oracle ADDOW which controls the FDR only for $\mu \geq 2$.

The situation get more confused when the signal is weak ($\mu < 1$). The FDR of ADDOW (in each group) is largely inflated. The FDR control at level $\alpha$ also fails sometimes for Oracle ADDOW, Pro2, ABH and HZZ (only in Group 2).

In scenario 3, one group has always weak signal. The FDR inflation of ADDOW (in each group) and Group 2 is worse for small $\bar\mu$,  whereas, for large $\bar\mu$, the situation is similar to scenario 1, up to one exception: the FDR of ABH and IHWc in Group 3 does not reach $\alpha$ as it did in scenario 1, which suggests some sort of conservatism.

In both scenarios, procedures of Group 2 have a larger FDR than their equivalent in Group 3, which in turn have larger FDR than in Group 1.

As a side note, in both scenarios, and both Groups 1 and 3, the FDR plots of IHWc and ABH are nearly indistinguishable.

In both settings regarding $\bar\mu$ (large or small), procedures based on $\widehat W^*$ suffer from some sort of overfitting causing a loss of FDR control. This is discussed in Section~\ref{section_overfitting} with an attempt to stabilize the weights. Let us underline that this does not contradict Theorem~\ref{thm_fdr} because a small $\mu_g$ might imply a smaller convergence rate while $m$ stays $<10^4$ in our setting. 

\end{version3}

%

\subsection{Power analysis}
\label{subsection_sim_pow}

Now that the FDR control has been studied, let us compare the procedures in terms of power. First, to better emphasize the benefit of adaptation, the power is rescaled in the following way: we define the normalized difference of power with respect to BH, or DiffPow, by $$\DiffPow(R)=\frac{m}{m_1}\left(\Pow(R)-\Pow(\BH)\right),$$ for any procedure $R$. 

\begin{figure}
\centering
\includegraphics[width=0.99\linewidth]{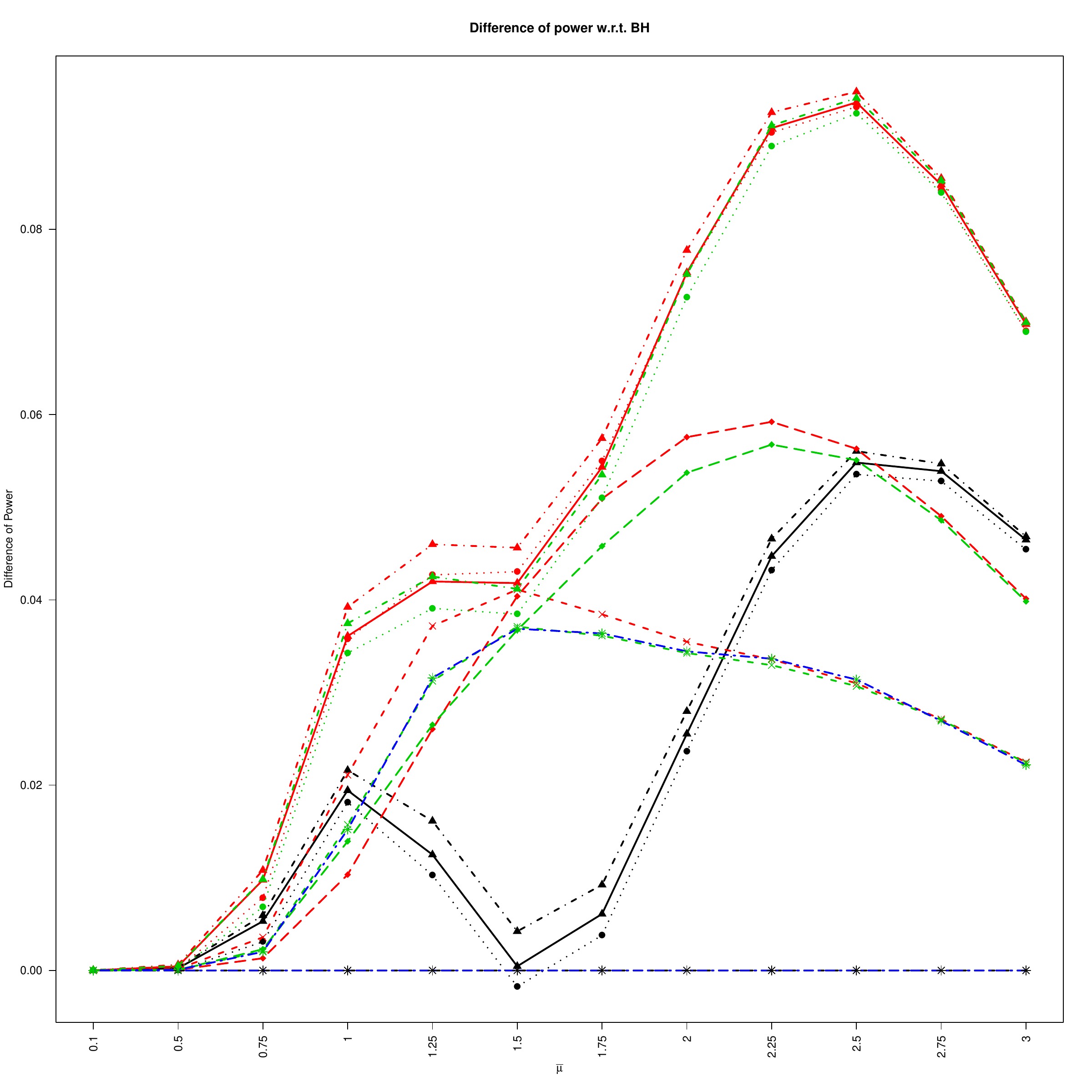}
\caption[DiffPow against \texorpdfstring{$\bar\mu$}{bar mu} in scenario 1]{\footnotesize \begin{version3}DiffPow against $\bar\mu$ in scenario 1. Same legend as Figure~\ref{fig_mult_fdr}.\end{version3}}
\label{fig_mult_pow}
\end{figure}

\begin{figure}
\centering
\includegraphics[width=0.99\linewidth]{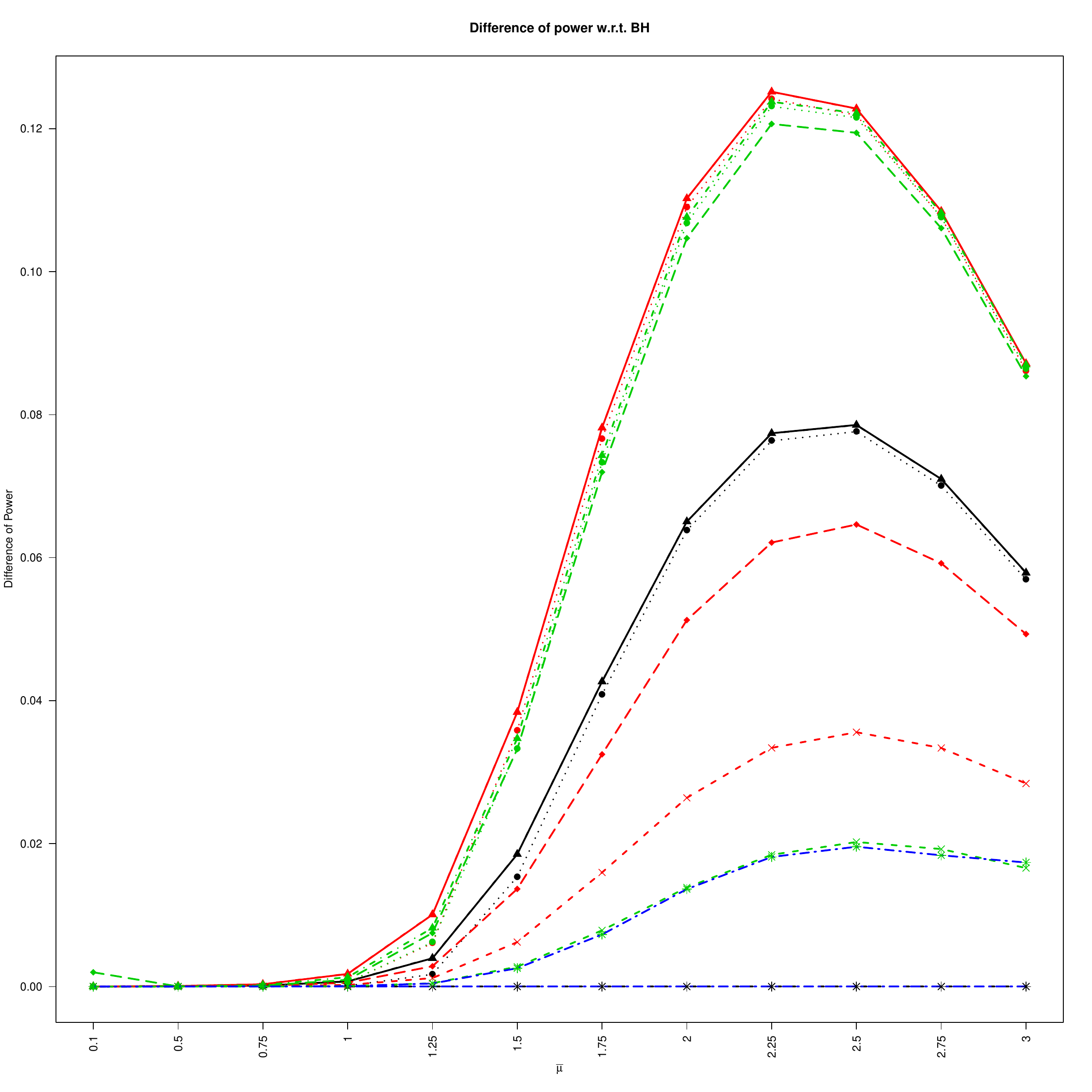}
\caption[DiffPow against \texorpdfstring{$\bar\mu$}{bar mu} in scenario 3]{\footnotesize \begin{version3}DiffPow against $\bar\mu$ in scenario 3. Same legend as Figure~\ref{fig_mult_fdr}.\end{version3}}
\label{fig_LARGEsmall_pow}
\end{figure}

\begin{figure}
\begin{minipage}[b]{.5\linewidth}
\centering
\includegraphics[scale=0.18]{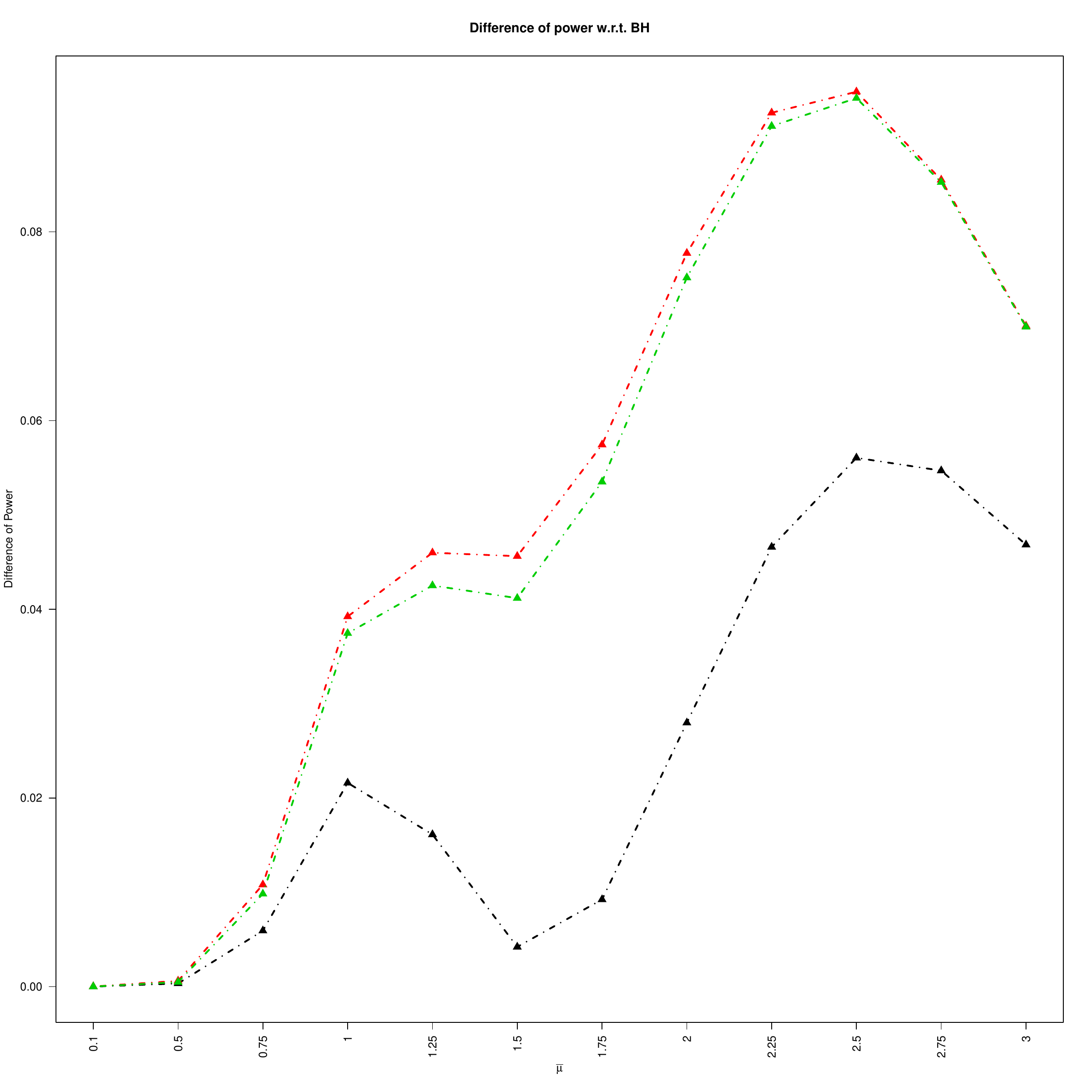}
\subcaption{\footnotesize \begin{version3}ADDOW in the four Groups\end{version3}}\label{fig_onlyADDOW}
\end{minipage}%
\begin{minipage}[b]{.5\linewidth}
\centering
\includegraphics[scale=0.18]{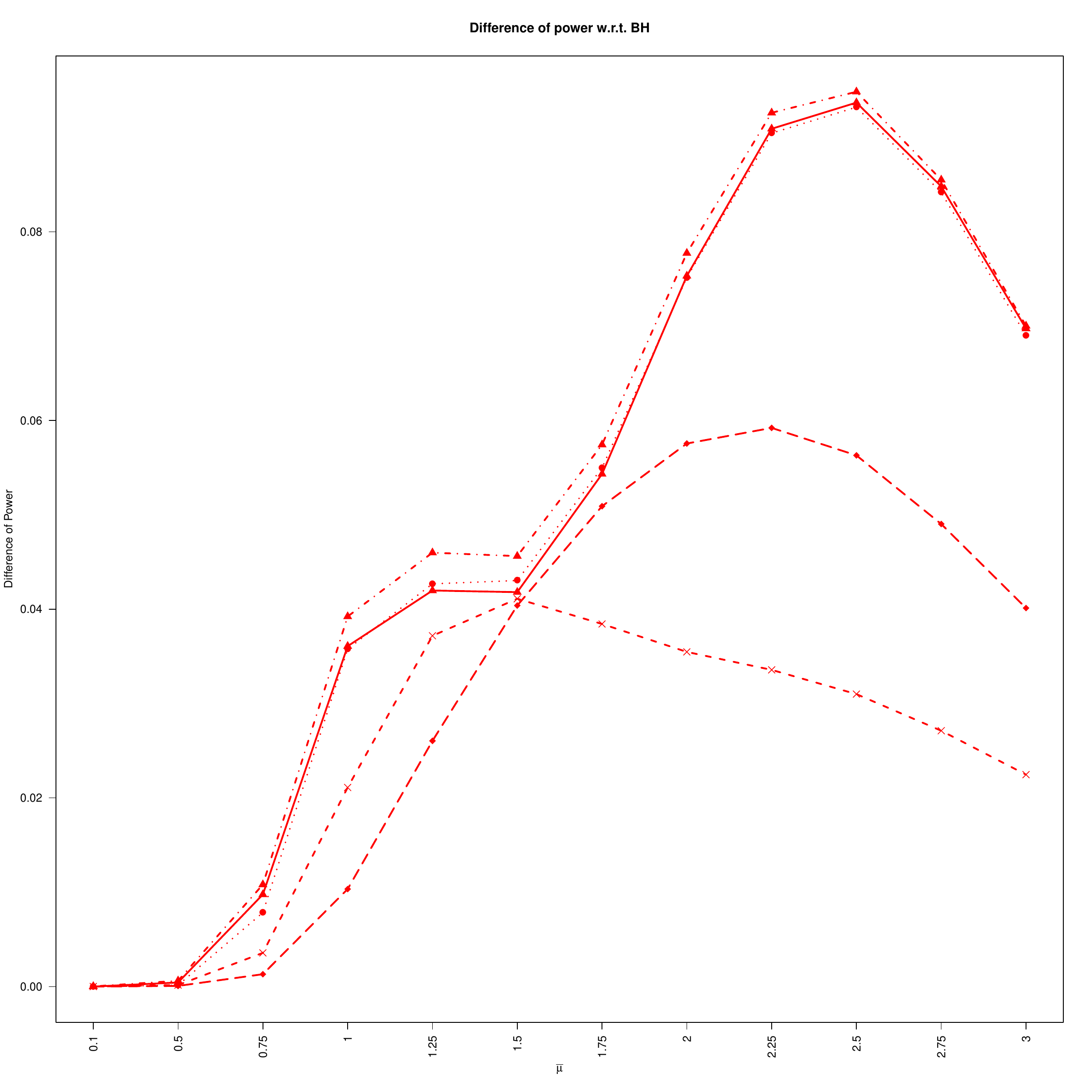}
\subcaption{\footnotesize \begin{version3}Procedures of Group 2\end{version3}}\label{fig_onlyCE}
\end{minipage}
\caption[Details of Figure~\ref{fig_mult_pow}]{\footnotesize \begin{version3}Details of Figure~\ref{fig_mult_pow} where only a subset of procedures is plotted.\end{version3}}\label{fig_stab}
\end{figure}

Figures~\ref{fig_mult_pow} and~\ref{fig_LARGEsmall_pow} display the power of all the procedures defined in Section~\ref{subsection_sim_setting}. Figures~\ref{fig_onlyADDOW} and~\ref{fig_onlyCE} display only a subset of them in Scenario 1, for clarity. We can make several observations:
\begin{itemize}
\item \begin{version3}In both scenarios, procedures of Group 2 are more powerful than their equivalent in Group 3, which are better than in Group 1 (up to one exception, see next point), see e.g. Figure~\ref{fig_onlyADDOW}. In particular, the difference between Group 2 and Group 1 is huge. This illustrates the importance of incorporating the knowledge of $\pi_0$ to improve power. \end{version3}
\item \begin{version3}In scenario 2, HZZ is largely better in Group 3 than in Group 2. Our interpretation is that the signal is so weak in the second group of $p$-values that the estimator $\hat \pi_{2,0}(1/2)$ is close to one, while $\hat \pi_{1,0}(1/2)$ stays close to $ \pi_{1,0}$. Hence $\hat w^{(2)}_1$ in Group 3 is larger than  $\hat w^{(2)}_1$ in Group 2 which allows for more good discoveries. The drawback of having $\hat w^{(2)}_2$  in Group 3 smaller than  $\hat w^{(2)}_2$ in Group 2 is not a real one since the signal is so small that it is impossible to detect no matter the weight. Recall that $\hat w^{(2)}$ is defined in Section~\ref{subsection_previous}.\end{version3}
\item \begin{version3}In every Group (that is for any choice of $\hat \pi_{g,0}$), and for both scenarios,\end{version3} ADDOW achieves the best power \begin{version3}(see e.g. Figure~\ref{fig_onlyCE})\end{version3}, which supports Theorem~\ref{thm_pow}. Additionnaly, maybe surprisingly, Pro2 behaves quite well, with a power close to the one of ADDOW \begin{version3}(sometimes larger than Oracle ADDOW)\end{version3} and despite its theoretical sub-optimality. 
\item \begin{version3}Inside Group 2 or Group 3, and for both scenarios, comparing ABH and HZZ 
to ADDOW and Pro2 shows the benefit of adding the $F_g$ adaptation to the $\pi_0$ adaptation:  the ADDOW and Pro2 have better power than ABH and HZZ for all signals (see e.g. Figure~\ref{fig_onlyCE}).\end{version3}
\begin{version3}In scenario 1, for Groups 2 and 3,\end{version3} we can see a zone of moderate signal (around $\bar\mu=1.5$) where the two categories of procedures are close. That is the same zone where HZZ becomes better than ABH. We deduce that in that zone the optimal weighting is the same as the uniform $\hat w^{(1)}$ weighting of ABH.
\item The comparison of the DiffPow between\begin{version3}, on the one hand, IHW and, on the other hand, ABH or HZZ from Group 2,\end{version3} in Figure~\ref{fig_mult_pow}, shows the difference between adapting only to the $F_g$'s versus adapting only to $\pi_0$. No method is generally better than the other: as we see in the plot, it depends on the signal strength. We also see that neither ABH nor HZZ is better than the other. 
\begin{version3}
\item In scenario 1, for all signals, methods of Group 3 are close to their equivalent of Group 2, which indicates that using $\lambda=1/2$ gives a good estimate of $\pi_{g,0}$ in practice (see e.g. Figure~\ref{fig_onlyADDOW}). Furthermore, the larger the signal is, the more methods of Group 3 get closer to Group 2.
\item In both scenarios, once again IHWc and ABH are nearly indistinguishable, which confirms the intuition given in Section~\ref{subsection_sabha_etc} that IHWc performs badly in terms of power due to using only large $p$-values to compute the weights. See in particular how the power of IHW is larger than the power of IHWc (and even than the power of IHWc-Storey) in Figure~\ref{fig_mult_pow}.
\end{version3}
\end{itemize}

\begin{figure}
\centering
\includegraphics[width=0.65\linewidth]
{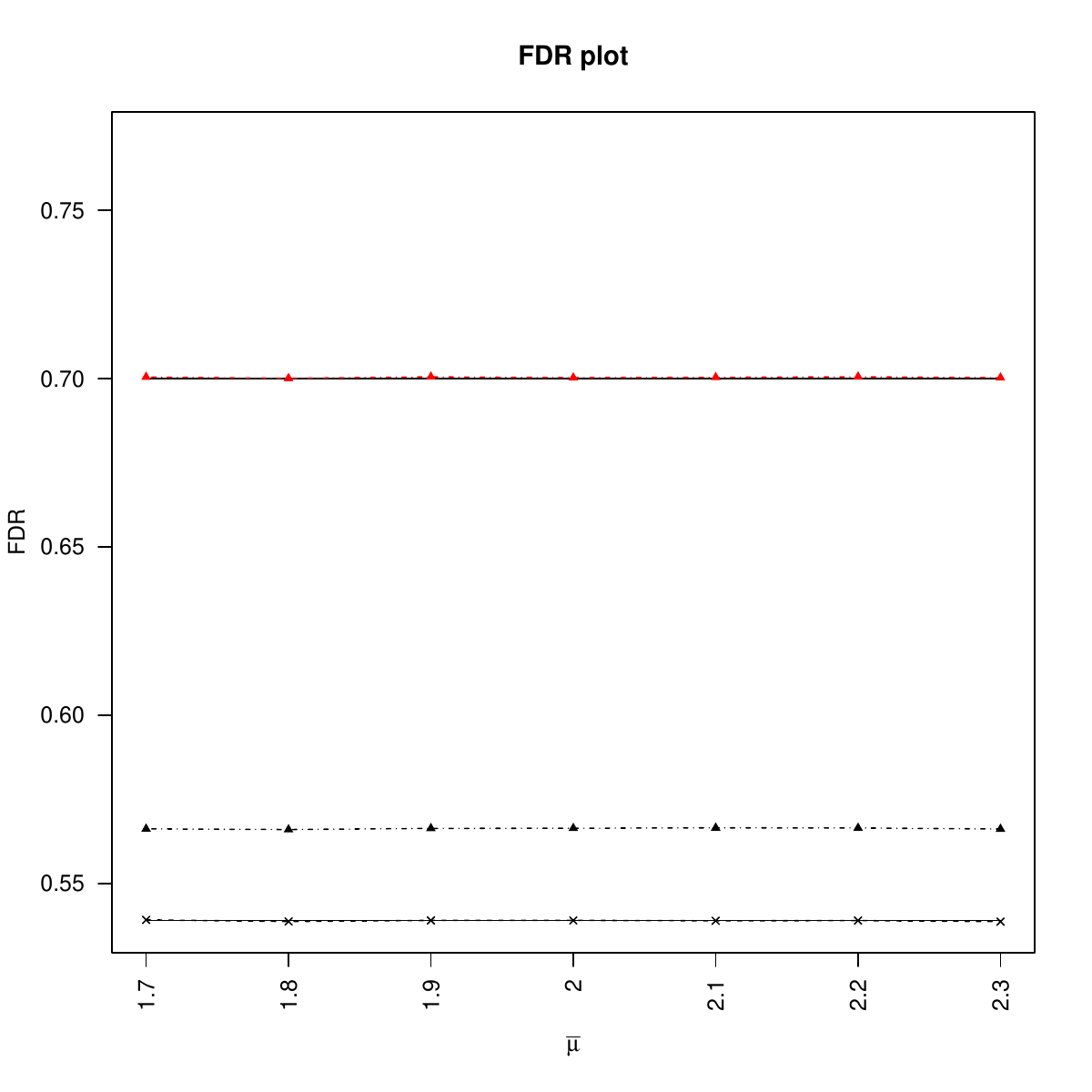}
\caption[FDR of ADDOW and BH against \texorpdfstring{$\bar\mu$}{bar mu} in the simulation of Section~\ref{subsection_sim_me}]{\footnotesize \begin{version3}FDR of ADDOW and BH against $\bar\mu$ in \begin{version3}the simulation of Section~\ref{subsection_sim_me}\end{version3}. \begin{version3}The two solid lines are the $\alpha$ and $\pi_0\alpha$ levels, the FDR of BH is confounded with the $\pi_0\alpha$ level. ADDOW in the~\eqref{NE} case is given by the black triangles and ADDOW in the~\eqref{CE} case is given by the red triangles.\end{version3}\end{version3}}
\label{fig_set2_fdr}
\end{figure}

\begin{figure}
\centering
\includegraphics[width=0.65\linewidth]{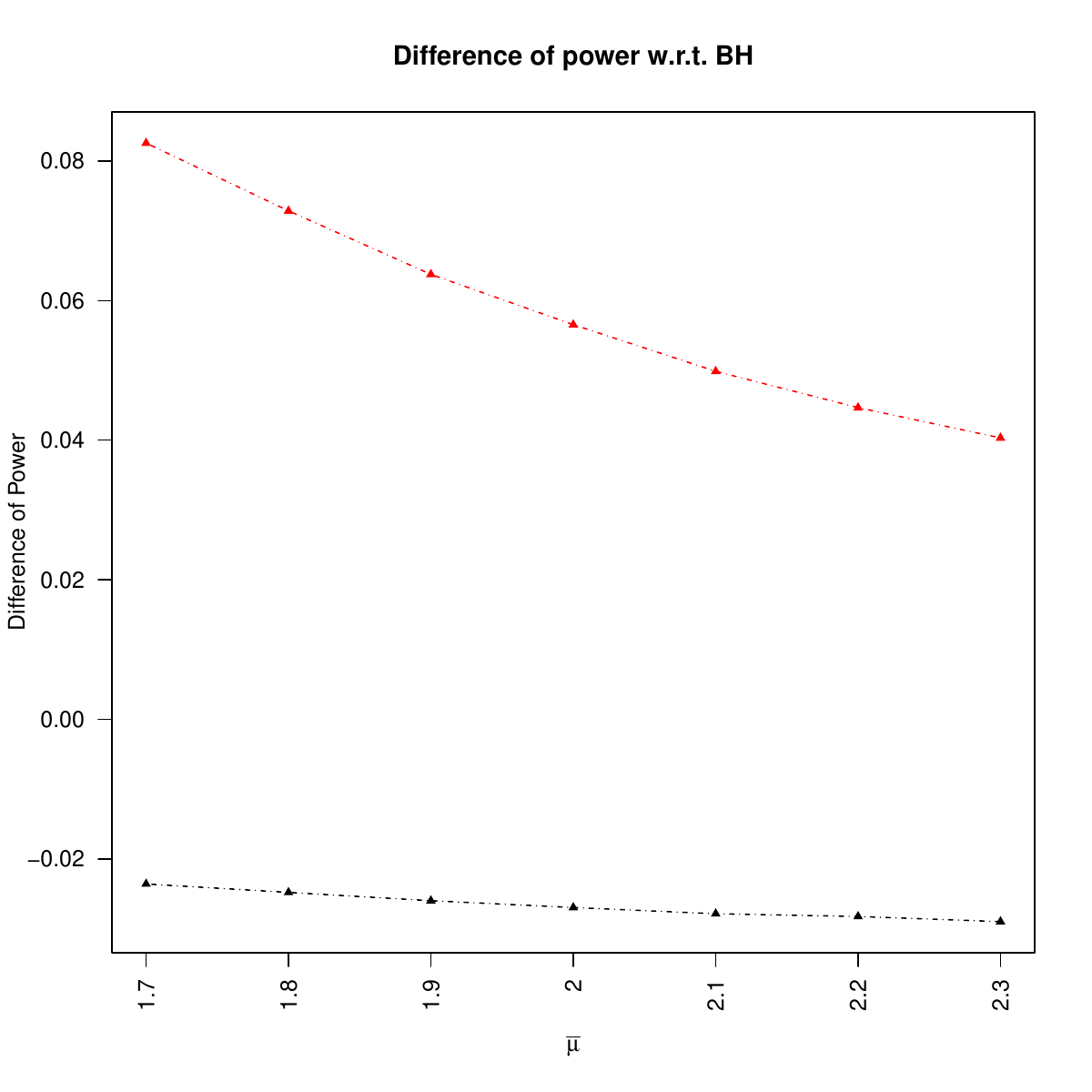}
\caption[DiffPow of ADDOW against \texorpdfstring{$\bar\mu$}{bar mu} in the simulation of Section~\ref{subsection_sim_me}]{\footnotesize DiffPow of ADDOW against $\bar\mu$ in \begin{version3}the simulation of Section~\ref{subsection_sim_me}\end{version3}. \begin{version3}ADDOW in the~\eqref{NE} case is given by the black triangles and ADDOW in the~\eqref{CE} case is given by the red triangles.\end{version3}}
\label{fig_set2_pow}
\end{figure}

\begin{version3}

\subsection{Importance of~\eqref{ME} for optimality results}\label{subsection_sim_me}

We provide here a setting and a simulation where Corollary~\ref{cor_ihw} fails because~\eqref{ED} does not hold, to illustrate the importance of~\eqref{ME} in Theorem~\ref{thm_pow} and in Theorem~\ref{thm_fdr} (to get~\eqref{eq_thm_fdr_ref}). The setting is chosen according to what we sketched in Remark~\ref{rk_bh_best} and is the following.

We consider again the one-sided Gaussian framework described in Section~\ref{subsection_gaussian} for $G=2$ groups and independent $p$-values. The parameters are the same as in Section~\ref{subsection_sim_setting} and each simulation of each scenario is replicated 1000 times. We choose a large value for $\alpha$ ($\alpha=0.7$) which is unlikely to appear in practice but allows us to get our counterexample. We set  $m_1=1000$ and $m_2=9000$, $m_{1,0}/m_1=0.05$ and $m_{2,0}/m_2=0.85$. So group 1 is small and has a lot of signal, while group 2 is large but has not much signal. The signal strength is given by $\mu_1=2$ and $\mu_2=\bar \mu$, and $\bar\mu\in\{1.7,1.8,1.9,2,2.1,2.2,2.3\}$, so the signal is strong and almost equal in both groups.

We compare only BH, ADDOW in the~\eqref{CE} case (with $\hat \pi_{g,0}=\pi_{g,0}$) and ADDOW in the~\eqref{NE} case (that is, IHW, with $\hat \pi_{g,0}=1$).
The simulation is illustrated with an FDR plot in Figure~\ref{fig_set2_fdr} and a DiffPow plot in Figure~\ref{fig_set2_pow}. 

In Figure~\ref{fig_set2_fdr}, the FDR of BH is $\pi_0\alpha$ as expected, and we see that the FDR of IHW is above that level, hence Equation~\eqref{eq:cor:fdr} is violated. On a side note, we see that, thanks to a large $m$ ($10^4$) and a rather strong signal, ADDOW in~\eqref{CE} does not overfit  and we get an illustration of Equation~\eqref{eq_thm_fdr_ref} with $C=1$.

Figure~\ref{fig_set2_pow} is rather unequivocal and shows that our parameter choice implies that IHW has a power smaller than BH (ADDOW in~\eqref{CE} case stays better as expected), hence Equation~\eqref{eq:cor:pow} is violated. Let us recall our interpretation proposed  in Remark~\ref{rk_bh_best}: IHW favors the large and sparse second group of hypotheses whereas the optimal power is achieved by favoring the small first group of hypotheses which contains almost only signal.\end{version3}  As a WBH procedure with weights (1,1), BH does not favor any group.
%
\begin{version3}Figure~\ref{fig_set2_pow}\end{version3} demonstrates the limitation of \begin{version3}Heuristic~\ref{mainheur} by providing a direct counterexample\end{version3}, and underlines the necessity of estimating the $\pi_{g,0}$ when nothing lets us think that~\eqref{ED} may be met. 

\begin{version3}

\section{Stabilization for overfitting}
\label{section_overfitting}

\subsection{Overfitting phenomena}\label{subsection_overfitting}

Since ADDOW uses the data both through the $p$-values and the weights, it suffers from an overfitting phenomena where the FDR in finite samples is above the target level $\alpha$, as we saw in Section~\ref{subsection_sim_fdr}.  In our setting, if the signal is strong enough, this drawback is proved to vanish when $m$ is large enough, see the simulations and Theorem \ref{thm_fdr}. However, the latter is not true for weak signal: if the data are close to be random noise, making the weight optimization leads ADDOW to \begin{version3}assign its weighting budget at random, and giving large weights to the wrong groups increases the FDP.\end{version3} 

As said before, our intuition is that the overfitting is at least partly due to using each $p$-value twice in the step-up procedure of ADDOW: in the expression $\ind{p_{g,i}\leq \alpha u \widehat W^*_g(u)}$, $p_{g,i}$ appears in both sides of the inequality because it is used to compute $\widehat W^*_g(u)$. Following this, we propose a variation of ADDOW that uses the same cross-weighting trick as IHWc.

\subsection{The $\crossADDOW$ variant}\label{subsection_crossADDOW}

The main idea is to split the $p$-values into $F$ folds, where $F$ is some fixed integer $\geq2$, and to use only $p$-values of the remaining $F-1$ folds to compute the weights assigned to the $p$-values of a given fold. The resulting procedure can be seen as a WBH procedure using $F\times G$ groups.

Formally, for each $m$ we have a random function $\mathbb{F}_m: (g,i)\mapsto \mathbb{F}_m(g,i)\in\{1,\dotsc,F\}$ such that, for each $f\in \{1,\dotsc,F\}$ and each $g\in \{1,\dotsc,G\}$, $| \{1\leq i\leq m_g : \mathbb{F}_m(g,i) =f \} |\geq \lfloor\frac{m_g}{F}\rfloor$, which simply means that the $p$-values of each group $g$ are evenly distributed between the $F$ folds. Some dependence assumptions are required:

\begin{assumption}\label{dep_folds}
The $\sigma$-algebra generated by $(\mathbb{F}_m)_m$ and the $\sigma$-algebra generated by $\left((p_{g,i})_{(g,i)} \right)_m$ are independent.
\end{assumption}

\begin{assumption}\label{weak_dep_folds}
Conditionally to $(\mathbb{F}_m)_m$, we have weak dependence (as in Assumption~\ref{ass4}) inside each fold.
\end{assumption}

For each fold $f\in\{1,\dotsc,F\}$, we compute $\ADDOW_{-f}$, that is ADDOW but using only $p$-values for the folds in $\{1,\dotsc,F\}\setminus\{f\}$. This is done by constructing the empirical function
$$   \widehat G^{-f}_w : u \mapsto | \{(g,i) : \mathbb{F}_m(g,i) \neq f \} |^{-1}\sum_{\substack{(g,i):\\\mathbb{F}_m(g,i)\neq f}}\mathds{1}_{\{ p_{g,i}\leq \alpha u w_g \}},$$
and then maximizing it in $w\in \Km^{-f}$ for each $u\in[0,1]$, where:
$$\Km^{-f}=\left\{w\in \mathbb{R}^G_+ : \sum_g \begin{version2}\frac{| \{1\leq i\leq m_g : \mathbb{F}_m(g,i) \neq f \} |}{| \{(g,i) : \mathbb{F}_m(g,i) \neq f \} |}\end{version2}  \hat \pi_{g,0} w_g\leq 1\right\}.$$
While this expression seems complicated, note that if $F$ divides each $m_g$, then $| \{1\leq i\leq m_g : \mathbb{F}_m(g,i) =f \} |=\frac{m_g}{F}$ and $\Km^{-f}=\Km$. The maximization provides a weight function $\widehat W^{*-f}$ and the MWBH procedure provides a step-up threshold $\hat u_{\widehat W^{*-f}}=\mathcal{I}\left( \widehat G^{-f}_{\widehat W^{*-f}}  \right)$. To lighten notation, let $w^*_{g,f}=\widehat W^{*-f}(\hat u_{\widehat W^{*-f}})$.

Our ADDOW variant, named $\crossADDOW$ for cross-ADDOW, is the WBH procedure which assigns the weight $w^*_{g,f}$ to all $p$-values $p_{g,i}$ such that $\mathbb{F}_m(g,i)=f$. Now, in $\ind{p_{g,i}\leq \alpha u w^*_{g,f}}$, $p_{g,i}$ is only used once. While we don't have a finite-sample result about $\crossADDOW$, we expect it to have a lesser FDR than ADDOW, especially for weak signal. We expect $\crossADDOW$ to act like a stabilization of ADDOW and to not lose the good performances of ADDOW when the signal is not weak. Those intuitions are verified in the simulations of Section~\ref{subsection_simu_cross}. Still, $\crossADDOW$ has the nice property of being asymptotically equivalent to ADDOW.

\begin{theorem}
Let us assume that Assumptions \cref{ass1,ass2,ass3,ass4,ass5,ass6,dep_folds,weak_dep_folds} are fulfilled. Assume also that $\alpha\leq\bar\pi_0$. We have
\begin{equation}
\lim_{m\to\infty} \FDR\left( \crossADDOW \right) =  \lim_{m\to\infty} \FDR\left( \ADDOW \right),  
\label{eq_thm_cross_fdr}
\end{equation}
and
\begin{equation}
\lim_{m\to\infty} \Pow\left( \crossADDOW \right) =  \lim_{m\to\infty} \Pow\left( \ADDOW \right). 
\label{eq_thm_cross_pow}
\end{equation}
\label{thm_cross}
\end{theorem}

This Theorem is proved in Section~\ref{section_proof_cross}.

\subsection{Simulations with $\crossADDOW$}
\label{subsection_simu_cross}

The simulations presented here are the same as the simulations depicted in Section~\ref{subsection_sim_setting}, with the addition of $\crossADDOW$ in each Group.

\begin{figure}
\centering
\includegraphics[width=0.65\linewidth]{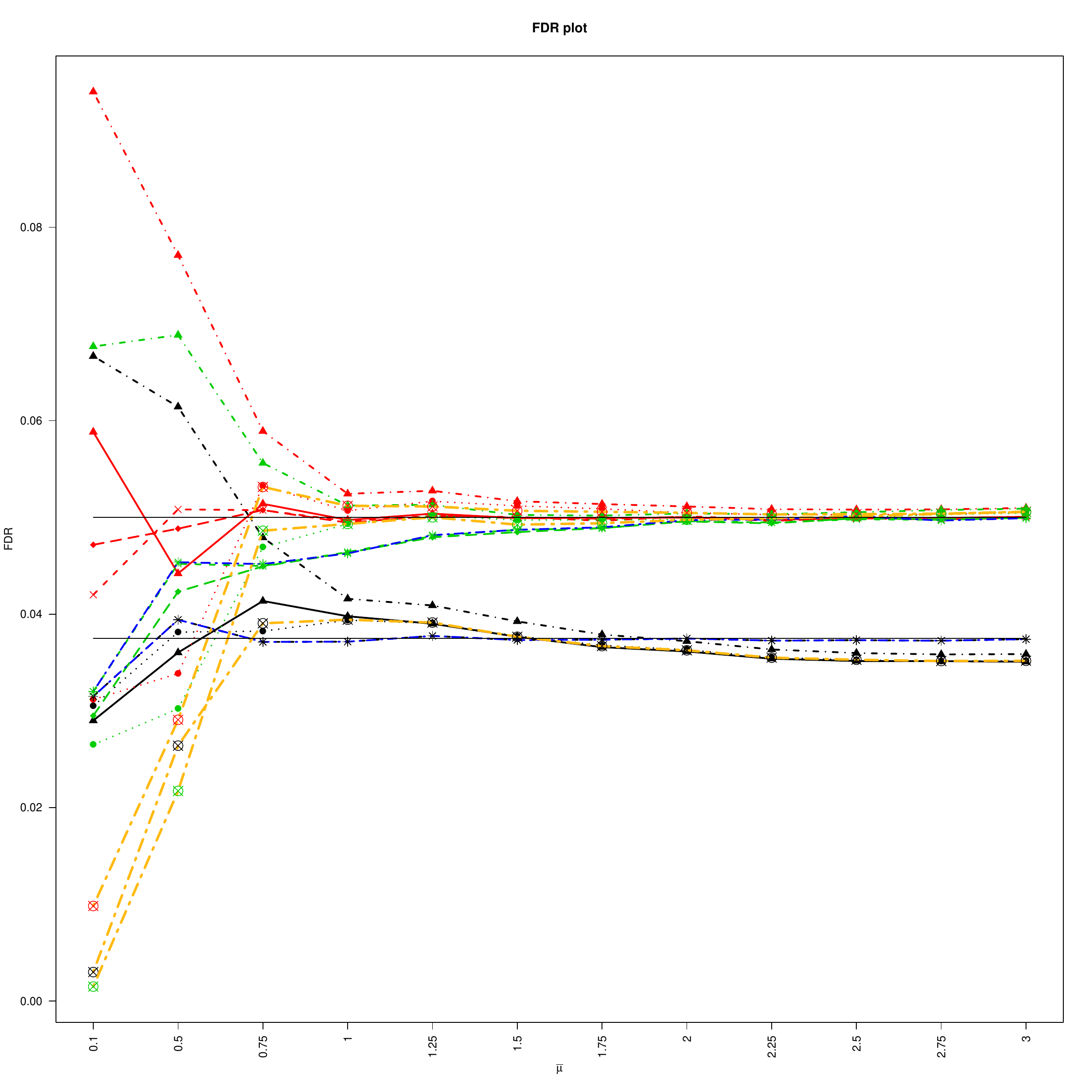}
\caption[FDR against \texorpdfstring{$\bar\mu$}{bar mu} in scenario 1]{\footnotesize \begin{version3}FDR against $\bar\mu$ in scenario 1. Same legend as in Figure~\ref{fig_mult_fdr_cr}, with the addition of $\crossADDOW$ (yellow lines). The color of the points (black, red, green) indicates the Group (respectively, 1, 2 and 3).\end{version3}}
\label{fig_mult_fdr_cr}
\end{figure}

\begin{figure}
\centering
\includegraphics[width=0.65\linewidth]{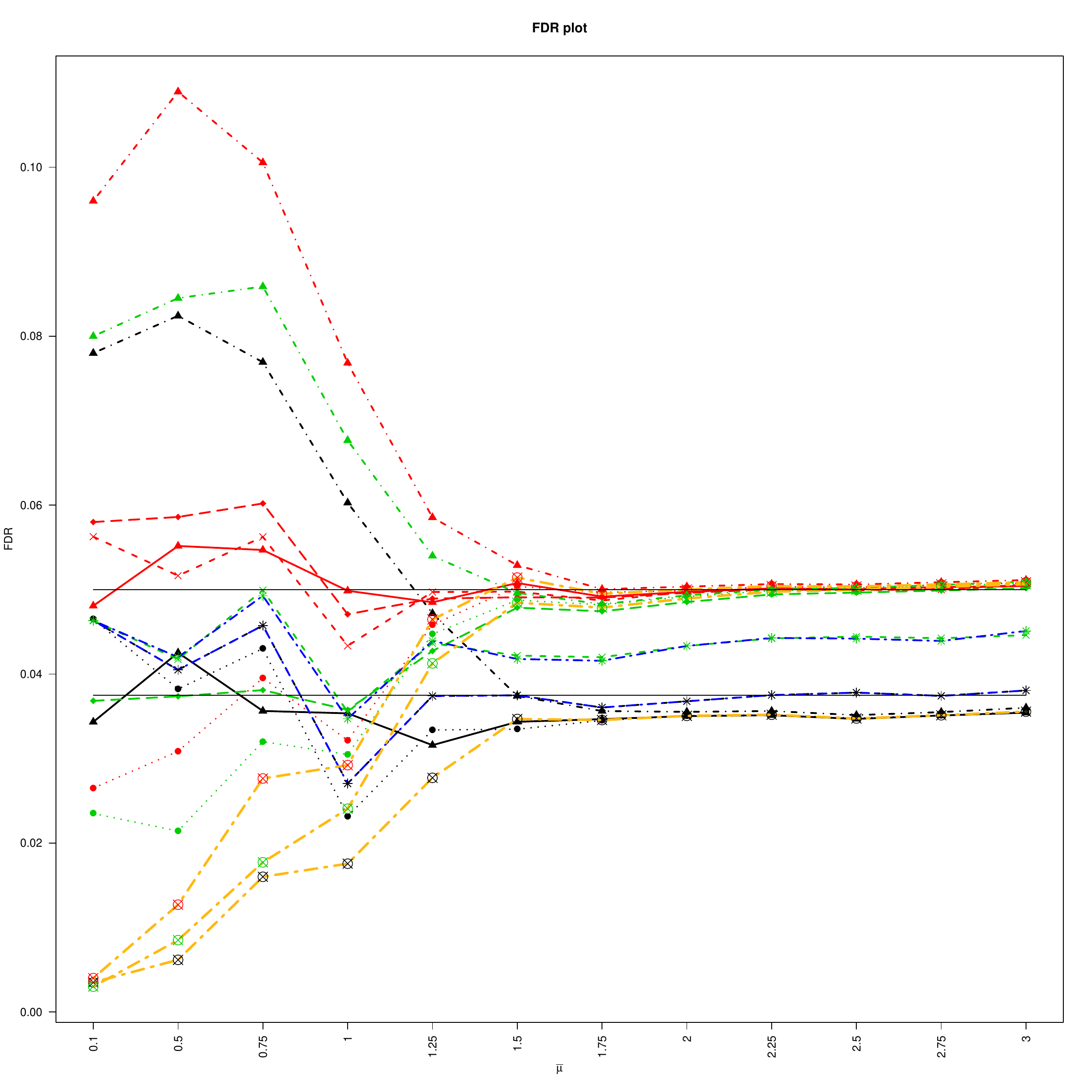}
\caption[FDR against \texorpdfstring{$\bar\mu$}{bar mu} in scenario 3]{\footnotesize \begin{version3} FDR against $\bar\mu$ in scenario 3. Same legend as in Figure~\ref{fig_mult_fdr_cr}.\end{version3}}
\label{fig_LARGEsmall_fdr_cr}
\end{figure}

\begin{figure}
\centering
\includegraphics[width=0.65\linewidth]{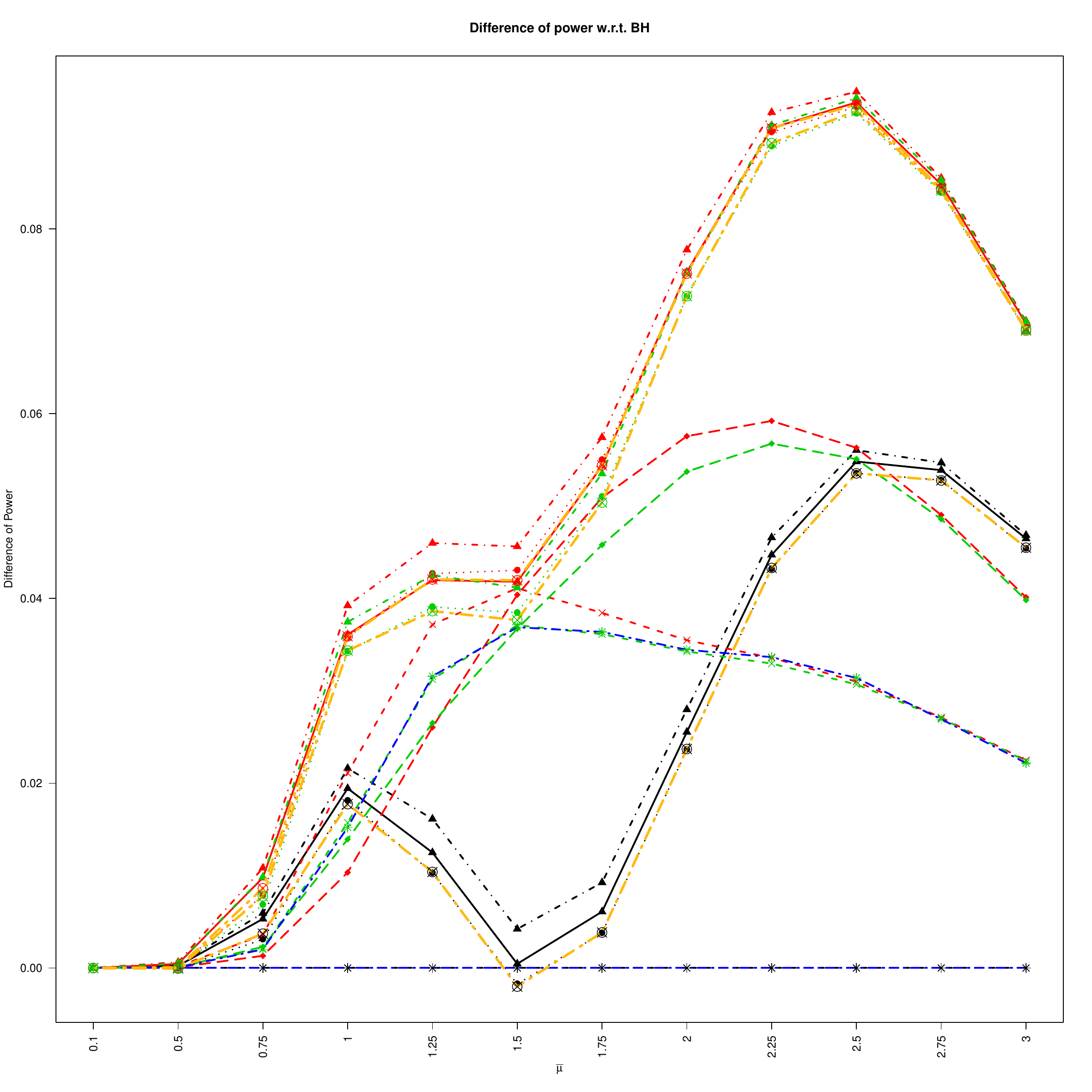}
\caption[DiffPow against \texorpdfstring{$\bar\mu$}{bar mu} in scenario 1]{\footnotesize \begin{version3}DiffPow against $\bar\mu$ in scenario 1. Same legend as Figure~\ref{fig_mult_fdr_cr}.\end{version3}}
\label{fig_mult_pow_cr}
\end{figure}

\begin{figure}
\centering
\includegraphics[width=0.65\linewidth]{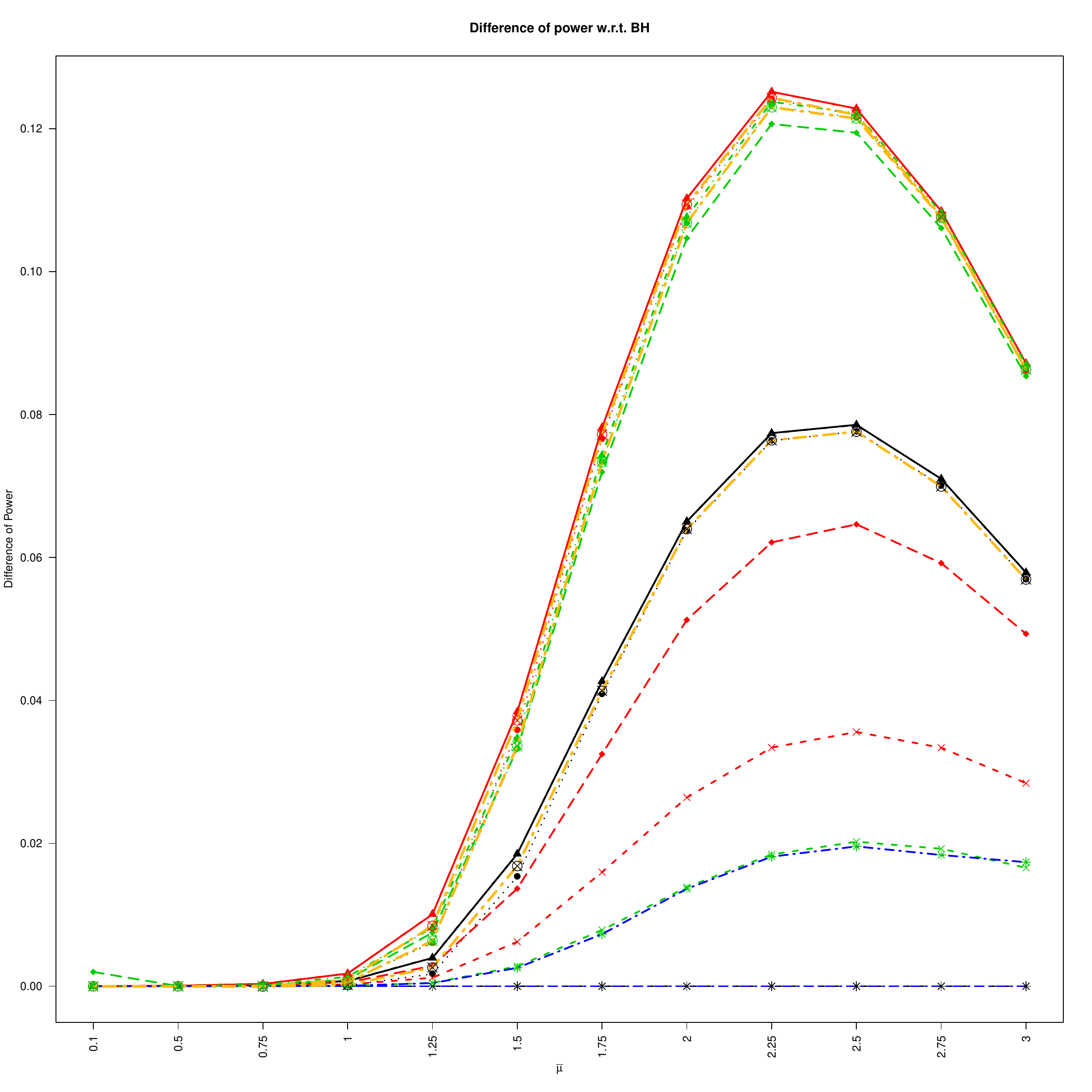}
\caption[DiffPow against \texorpdfstring{$\bar\mu$}{bar mu} in scenario 3]{\footnotesize \begin{version3}DiffPow against $\bar\mu$ in scenario 3. Same legend as Figure~\ref{fig_mult_fdr_cr}.\end{version3}}
\label{fig_LARGEsmall_pow_cr}
\end{figure}

From the FDR plots, we see that the FDR is hugely deflated and is now controlled at level $\alpha$ for weak $\bar\mu$ in each scenario, while for large $\bar\mu$ we are still slightly above the target level but with a small improvement over ADDOW. In scenario 1 there is a small window between large and small $\bar\mu$, around $\bar\mu=0.75$, where $\crossADDOW$ in Group 2 overfits more than for really large $\bar\mu$, but even there we see a large improvement over ADDOW.

As for the power, we see that $\crossADDOW$ is less powerful than ADDOW, as expected since we reject less hypotheses, but we see that in most Groups and scenarios the loss of power is almost negligible and $\crossADDOW$ remains even as powerful as Oracle ADDOW (with the exception of Group 1 in scenario 1). The difference of power between $\crossADDOW$ and Pro2 is even smaller and $\crossADDOW$ is better in most configurations, with the exception of Groups 2 and 3 around $\bar\mu=1.5$, which is the zone that we identified in Section~\ref{subsection_sim_pow} as the zone where the optimal weights are given by the uniform $\hat w^{(1)}$ weighting of ABH.

The simulations hence confirm our intuitions about the stabilization properties of $\crossADDOW$ especially for weak signal where ADDOW was totally unreliable. Studying the finite sample properties of $\crossADDOW$, especially its FDR, is an interesting direction for future works.

\end{version3}

\section{Concluding remarks}
\label{section_conclusion}

In this paper we presented a new class of data-driven step-up procedures, ADDOW, that generalizes IHW by incorporating $\pi_{g,0}$ estimators in each group. We showed that while this procedure asymptotically controls the FDR at the targeted level, it has the best power among all MWBH procedures when the $\pi_0$ estimation can be made consistently. In particular it dominates all the existing procedures of the weighting literature and solves the $p$-values weighting issue in a group-structured multiple testing problem. As a by-product, our work established the optimality of IHW in the case of homogeneous $\pi_0$ structure. Finally we proposed a stabilization variant designed to deal with the case where only few discoveries can be made (very small signal strength or sparsity). Some numerical simulations illustrated that our properties are also valid in a finite sample framework, provided that the number of tests \begin{version3}and the signal strength are\end{version3} large enough. \begin{version3}We also introduced $\crossADDOW$, a variant of ADDOW that uses cross-weighting to reduce the overfitting while having the exact same asymptotic properties.\end{version3}
 
\paragraph{Assumptions}
Our assumptions are rather mild: basically we only added the concavity of the $F_g$ to the assumptions of \citet{ignatiadis2016data}. Notably we dropped the other {regularity} assumptions on $F_g$ that were made in \citet{roquain2009optimal} while keeping all the useful properties on $W^*$ in the \eqref{NE} case. Note that the criticality assumption is often made in the literature, see \citet{ignatiadis2016data} (assumption 5 of the supplementary material), \citet{zhao2014weighted} (assumption A.1), or the assumption of Theorem 4 in \citet{hu2010false}. Finally, the weak dependence assumption is extensively used in our paper.  An interesting direction could be to extend our result to some strong dependent cases, for instance by assuming the PRDS (positive regression dependence), as some previous work already studied properties of MWBH procedures under that assumption, see \citet{roquain2008multi}.


{
\paragraph{Computational aspects}
The actual maximization problem of ADDOW is difficult, it involves a mixed integer linear programming that may take a long time to resolve. Some regularization variant may be needed for applications.}{ To this end, we can think to use the least concave majorant (LCM) instead of the empirical c.d.f. in equation~\eqref{eq_def_gW} {(as proposed in modification (E1) of IHW in \citealp{ignatiadis2016data})}. As we show in Section~\ref{section_proof_thm}, ADDOW can be extended to that case (see especially Section~\ref{subsection_notation}) and our results are still valid for this new regularized version of ADDOW.
}

\paragraph{Toward nonasymptotic results}
Interesting direction for future research can be to investigate the convergence rate in our asymptotic results. One possible direction can be to use the work of \citet{neuvial2008asymptotic}. However, it would require to compute the Hadamard derivative of the functional involved in our analysis, which might be very challenging. Finally, another interesting future work could be to develop other versions of ADDOW that ensure finite sample FDR control property: this certainly requires to use a different optimization process, which will make the power optimality difficult to maintain. \begin{version3}A possible such variation is $\crossADDOW$, whose FDR in finite sample has yet to be investigated.\end{version3}

\section{Proofs of Theorems~\ref{thm_fdr} and~\ref{thm_pow}}
\label{section_proof_thm}

\subsection{{Further generalization}}
\label{subsection_notation}

Define, for any $u$ and $W$,
$$\widehat H_{W}(u)=m^{-1}\left|R_{u,W} \cap\mathcal{H}_0  \right| =m^{-1}\sum_{g=1}^G   \sum_{i=1}^{m_g}\mathds{1}_{\{p_{g,i}\leq \alpha u W_g(u) ,H_{g,i}=0\}}  , $$
and
$$\widehat P_W(u)=m^{-1}\left|R_{u,W} \cap\mathcal{H}_1  \right| =\widehat G_{W}(u)-\widehat H_{W}(u) ,$$
so that $\FDP\left(  R_{u,W}\right)=\frac{\widehat H_{W}(u)}{\widehat G_{W}(u) \vee m^{-1}}$ and $\Pow\left(R_{u,W}\right)=\mathbb{E}\left[  \widehat P_W(u) \right]$ (recall that $\MWBH\left(W\right)$ \begin{version2}is\end{version2} $R_{\hat u_W,W}$). 
Also define $\widehat D_g(t)=m_g^{-1}\sum_{i=1}^{m_g}\ind{p_{g,i}\leq t}$ so that $\widehat G_{W}(u)=\sum_g \frac{m_g}{m}\widehat D_g(\alpha u W_g(u))$. 

For the sake of generality $\widehat D_g$ is not the only estimator of $D_g$ {(defined in equation~\eqref{def_Dg})} that we will use to prove our results (for example, we can use the LCM of $\widehat D_g$, denoted $\LCM(\widehat D_g)$, see Section~\ref{section_conclusion}). So let us increase slightly the scope of the {MWBH} class by defining $ \widetilde G_W(u)=\sum_g\frac{m_g}{m}\widetilde D_g(\alpha u W_g(u))$ for any estimator $\widetilde D_g$ such that $\widetilde D_g$ is nondecreasing, $\widetilde D_g(0)=0$, $\widetilde D_g(1)=1$ and $\left\|\widetilde D_g-D_g  \right\|\overset{\mathbb{P}}{\longrightarrow}0$, where $\|\cdot\|$ is the sup norm for the bounded functions on their definition domain. Note that at least $(D_g)_g$, $(\widehat D_g)_g$ (by Lemma~\ref{supas}), and $\left(\LCM(\widehat D_g)\right)_g$ {(by Lemma~\ref{lm_lcm})}  are eligible. 

If $W$ is such that $ \widetilde G_W$ is nondecreasing, we then define the generalized MWBH as
\begin{equation}
\GMWBH\left( (\widetilde D_g)_g ,W\right) =R_{\tilde u_{ W},  W} \text{ where }
 \tilde u_{W} = \mathcal{I}\left( \widetilde G_{W}  \right).\notag
\end{equation}

If $(\widetilde D_g)_g$ is such that we can define, for all $u\in[0,1]$, 
\begin{equation}
\widetilde W^*(u)\in\argmax_{w\in \Km}\widetilde G_w(u),
\label{def_wtilde}
\end{equation}
we define the generalized ADDOW by
$$\GADDOW\left( (\widetilde D_g)_g \right) = \GMWBH\left( (\widetilde D_g)_g,  \widetilde W^* \right),$$
the latter being well defined because $ \widetilde G_{\widetilde W^*}$ is nondecreasing (by a proof similar to the one of  Lemma~\ref{croissance}). Note that for any continuous $\widetilde D_g$, such as $\LCM(\widehat D_g)$ or $D_g$ itself, the arg max in~\eqref{def_wtilde} is non empty and GADDOW can then be defined.

What we show below are more general theorems, valid for any $\GADDOW\left( (\widetilde D_g)_g \right) $. Our proofs combined several technical lemmas deferred to Sections~\ref{subsection_asymptotical} and~\ref{subsection_technic}, which are based on the previous work of \citet{roquain2009optimal,hu2010false,zhao2014weighted}.

\begin{remark}
$\GADDOW\left( (\widetilde D_g)_g \right) $ when $\widetilde D_g=\LCM(\widehat D_g)$ and $\hat\pi_{g,0}=1$ \begin{version2}is exactly the same as IHW with modification (E1) defined in the supplementary material of \citet{ignatiadis2016data}.
In our notation, the latter is $\WBH\left(  \widetilde W^*\left(   \tilde u_{\widetilde W^*}  \right)  \right)$, which is the same as $\GADDOW\left( (\widetilde D_g)_g \right) $ because $\tilde u_{\widetilde W^*}     =   \tilde u_{ \widetilde W^*\left(   \tilde u_{\widetilde W^*}\right)}      $ (same proof as in Remark~\ref{MWBH<wbh}).
\end{version2}
\end{remark}

\subsection{Proof of Theorem~\ref{thm_fdr}}
\label{subsection_proof_thm_fdr}

We have
$$\FDP\left(  \GMWBH\left(\left(\widetilde D_g\right)_g, \begin{version2}\widetilde W^*\end{version2} \right) \right)=\frac{\widehat H_{\widetilde W^*}(\tilde u)}{\widehat G_{\widetilde W^*}(\tilde u)\vee m^{-1}}\in[0,1],$$
where {$\tilde u$ is defined as in~\eqref{def_uwtilde}} so by Lemma~\ref{BOUMH} we deduce that
$$\FDP\left(\GADDOW\left((\widetilde D_g)_g \right)\right)\underset{m\to\infty}{\overset{\mathbb{P}}{\longrightarrow}} \frac{ H^\infty_{W^*}(u^*)}{G^\infty_{W^*}( u^*)}=\frac{ H^\infty_{W^*}(u^*)}{u^*} , $$
and then
$$\lim_{m\to\infty}\FDR\left(\GADDOW\left((\widetilde D_g)_g \right)\right)={u^*}^{-1}   H^\infty_{W^*}(u^*),$$
where $G^\infty_{W^*}$, $H^\infty_{W^*}$ and $u^*$ are defined in Section~\ref{subsection_asymptotical}.

If $\alpha \geq\bar \pi_0$, $u^*=1$ by Lemma~\ref{prop_max} and $\alpha u^* W^*_g(u^*)\geq1$ by Lemma~\ref{prop_argmax} so ${u^*}^{-1} H^\infty_{W^*}(u^*)=\pi_0\leq\bar \pi_0\leq\alpha$.

If $\alpha \leq\bar \pi_0$, $\alpha u^* W^*_g(u^*)\leq1$ by Lemma~\ref{prop_argmax} so $U(\alpha u^* W^*_g(u^*))=\alpha u^* W^*_g(u^*)$ for all $g$ and then
\begin{align}
{u^*}^{-1} H^\infty_{W^*}(u^*)&=\alpha\sum_g\pi_g\pi_{g,0}W^*_g(u^*)\notag\\
&\leq\alpha\sum_g\pi_g\bar \pi_{g,0}W^*_g(u^*)=\alpha.\label{egg1}
\end{align}

Moreover \begin{version2}if~\eqref{ME} holds (that is, there exists $C\geq1$ such that $\bar \pi_{g,0}=C\pi_{g,0}$ for all $g$), we write
\begin{align}
{u^*}^{-1} H^\infty_{W^*}(u^*)&=\alpha\sum_g\pi_g\pi_{g,0}W^*_g(u^*)\notag\\
&=\frac{\alpha}{C}\sum_g\pi_g\bar \pi_{g,0}W^*_g(u^*)=\frac{\alpha}{C}.\label{egg2}
\end{align}
\end{version2}
The equalities in~\eqref{egg1} and~\eqref{egg2} are due to $\sum_g\pi_g\bar \pi_{g,0} W^*_g(u^*)=1$ (by Lemma~\ref{prop_argmax}).

\subsection{Proof of Theorem~\ref{thm_pow}}
\label{subsection_proof_thm_pow}

First, in any case,
$$\widehat P_{\widetilde W^*}(\tilde u)= \widehat G_{\widetilde W^*}(\tilde u)-\widehat H_{\widetilde W^*}(\tilde u)\overset{a.s.}{\longrightarrow} G^\infty_{W^*}(u^*)-H^\infty_{W^*}(u^*)=P^\infty_{W^*}(u^*) $$
by Lemma~\ref{BOUMH}, where $P^\infty_{W^*}$ is defined in Section~\ref{subsection_asymptotical}. Hence:
\begin{version2}
\begin{equation*}
\lim_{m\to\infty} \Pow\left( \GADDOW\left((\widetilde D_g)_g \right) \right)= P^\infty_{W^*}(u^*).
\end{equation*}
\end{version2}

For the rest of the proof, we assume we are in \begin{version2}case~\eqref{ME}\end{version2}, which implies by Lemma~\ref{ASTUCE} that $W^*(u)\in\argmax_{w\in K^\infty}P^\infty_w(u)$ for all $u$, and that $P^\infty_{W^*}$ is nondecreasing. We also split the proof in two parts. For the first part we assume that for all $m$, $\widehat W$ is a weight vector $\hat w \in \Km$ therefore not depending on $u$. In the second part we will conclude with a general sequence of weight functions.

\paragraph{Part 1}$\widehat W=\hat w \in \Km$ for all $m$. Let $\ell=\limsup \Pow\left(\MWBH\left(\hat w  \right) \right)$. Up to extracting a subsequence, we can assume that $\ell=\lim\Esp{\widehat P_{\hat w}(\hat u_{\hat w})}$, $\hat\pi_{g,0}\overset{a.s.}{\longrightarrow}\bar \pi_{g,0}$ for all $g$\begin{version3}, and that the convergences of Lemma~\ref{supas} are almost surely.\end{version3}.
Define the event
\begin{equation*}
\widetilde \Omega=\left\{
\begin{array}{rcc}
\forall g,\,\hat\pi_{g,0} &\longrightarrow   &\bar \pi_{g,0}  \\
\sup_{w\in  \mathbb{R}_+^G } \left\| \widehat P_w-P^\infty_w    \right\|  &\longrightarrow   &0   \\
\sup_{w\in  \mathbb{R}_+^G } \left\| \widehat G_w-G^\infty_w    \right\|  &\longrightarrow   &  0 
\end{array}
\right\}
\end{equation*}
then $\Pro{\widetilde\Omega}=1$, $\ell=\lim
\Esp{\widehat P_{\hat w}(\hat u_{\hat w})\mathds{1}_{\widetilde\Omega}}$ and by reverse Fatou Lemma $\ell\leq\Esp{\limsup
\widehat P_{\hat w}(\hat u_{\hat w})\mathds{1}_{\widetilde\Omega}}$. 

Now consider that $\widetilde\Omega$ occurs and fix a realization of it, the following of this part 1 is deterministic. Let $\ell'=\limsup\widehat P_{\hat w}(\hat u_{\hat w})$. The sequences $\left(\frac{m}{m_g\hat\pi_{g,0}}\right)$ are converging and then bounded, hence the sequence $(\hat w)$ is also bounded. By compacity, once again up to extracting a subsequence, we can assume that $\ell'=\lim\widehat P_{\hat w}(\hat u_{\hat w})$ and that $\hat w$ converges to a given $w^{cv}$. By taking $m\to\infty$ in the relation $\sum\frac{m_g}{m}\hat\pi_{g,0}\hat w_g\leq1$, it appears that $w^{cv}$ belongs to $K^\infty$. $\| \widehat G_{\hat w}-G^\infty_{w^{cv}}    \| \leq \sup_w \| \widehat G_w-G^\infty_w\| +\|  G^\infty_{\hat w}  -G^\infty_{w^{cv}}   \|  \to0$ so by Remark~\ref{weight_vector} $\hat u_{\hat w}\to u^\infty_{w^{cv}}$ 
and finally
\begin{align*}
\left| \widehat P_{\widehat w}(\hat u_{\widehat w})-P^\infty_{w^{cv}}(u^\infty_{w^{cv}})    \right|  &\leq \sup_{w\in  \mathbb{R}_+^G } \left\| \widehat P_w-P^\infty_w    \right\|  +\left|P^\infty_{\widehat w}(\hat u_{\widehat w}) -P^\infty_{w^{cv}}(u^\infty_{w^{cv}})     \right|\\
   &\overset{}{\longrightarrow}0, 
\end{align*}
by continuity of $F_g$ and because \begin{version2}$\widetilde\Omega$ is realized.\end{version2} So $\ell'=P^\infty_{w^{cv}}(u^\infty_{w^{cv}})\leq P^\infty_{W^*}(u^\infty_{w^{cv}})$ by maximality. Note also that $G^\infty_{w^{cv}}(\cdot)\leq G^\infty_{W^*}(\cdot)$ which implies that $u^\infty_{w^{cv}}\leq u^\infty_{W^*}=u^*$ so $\ell'\leq P^\infty_{W^*}(u^*)$ because $P^\infty_{W^*}$ is nondecreasing. Finally $\limsup
\widehat P_{\hat w}(\hat u_{\hat w})\mathds{1}_{\widetilde\Omega}\leq P^\infty_{W^*}(u^*)$ for any realization of \begin{version3}$\Omega$\end{version3}, by integrating we get that $\ell\leq P^\infty_{W^*}(u^*)$ which concludes that part 1.

\paragraph{Part 2}Now consider the case where $\widehat W$ is a weight function $u\mapsto \widehat W(u)$. Observe that
$$\hat u_{\widehat W}= \widehat G_{\widehat W}(\hat u_{\widehat W})= \widehat G_{\widehat W(\hat u_{\widehat W})}(\hat u_{\widehat W}), $$
so by definition of $\mathcal{I}(\cdot)$, $\hat u_{\widehat W}\leq \hat u_{\widehat W(\hat u_{\widehat W})}$, and then
$$\widehat P_{\widehat W}(\hat u_{\widehat W})=\widehat P_{\widehat W(\hat u_{\widehat W})}(\hat u_{\widehat W})\leq  \widehat P_{\widehat W(\hat u_{\widehat W})}\left(\hat u_{\widehat W(\hat u_{\widehat W})}\right)  .$$
As a consequence, $\Pow\left(\MWBH\left(\widehat W\right)\right)\leq \Pow\left(\WBH\left(\widehat W(\hat u_{\widehat W})\right)\right)$. Finally, apply part 1 to the weight vector sequence $\left(\widehat W(\hat u_{\widehat W})\right)$ to conclude.

\begin{remark}
We just showed that for every MWBH procedure, there is a corresponding WBH procedure with better power. In particular, by defining $\hat u=u_{\widehat W^*}$ the ADDOW threshold, we showed that $\hat u\leq \hat u_{\widehat W^*(\hat u)}$. But $\widehat G_{\widehat W^*}\geq \widehat G_{\hat w}$ and then $\hat u\geq u_{\hat w}$ for any $\hat w$. Hence $\hat u = \hat u_{\widehat W^*(\hat u)}$ and ADDOW is \begin{version3}equal to\end{version3} the WBH procedure associated to the weight vector $\widehat W^*(\hat u)$.
\label{MWBH<wbh}
\end{remark}

\begin{version3}
\begin{remark}
We actually proved a stronger result, as we can replace the statement $\widehat W:[0,1]\to \Km$ by $\widehat W:[0,1]\to \Km^{\mathrm{alt}}$ where $\Km^{\mathrm{alt}}=\left\{w\in \mathbb{R}^G_+ : \sum_g \frac{m_g}{m}  \hat \pi_{g,0}^{\mathrm{alt}} w_g\leq 1\right\}$ and the $\hat \pi_{g,0}^{\mathrm{alt}}$ are such that $\hat \pi_{g,0}^{\mathrm{alt}}\overset{\mathbb{P}}{\longrightarrow} \bar \pi_{g,0}^{\mathrm{alt}}$ for some $\bar \pi_{g,0}^{\mathrm{alt}} \geq \bar \pi_{g,0}$. That is, the weight space $\widehat W$ belongs to does not have to be the same weight space where we apply ADDOW, as long as it uses over-estimators of the limits of the over-estimators used in $\Km$.
\end{remark}
\end{version3}

\section*{Acknowledgments}

\begin{version3}
Thanks to my former PhD advisor Etienne Roquain for his help in many areas of the paper. Notably, he simplified and strengthened the proof of Lemma~\ref{supG}, found the proof of the concavity in Lemma~\ref{prop_max} by introducting $\tilde w$, and improved my algorithm to compute $\widehat W^*$. Thanks to my other former PhD advisor Pierre Neuvial for his advices and for introducing me to Rcpp.
\end{version3}
Also thanks to Christophe Giraud and Patricia Reynaud-Bouret for useful discussions \begin{version3}and advices\end{version3}. \begin{version2}Finally, thanks to two referees and an associate editor for improving comments. In particular, the idea of $\crossADDOW$ comes from both one referee and from Etienne.\end{version2} 

This work has been supported by CNRS (PEPS FaSciDo) and ANR-16-CE40-0019 \begin{version3}(SansSouci)\end{version3}.

\bibliographystyle{imsart-nameyear}
\bibliography{hal_these.bib}
\appendix
\section{Lemmas and proofs of Section~\ref{section_framework}} 
\label{subsection_prooflemmas}

\begin{lemma}
For all $g$, $F_g$ is continuous.
\label{Fgcont}
\end{lemma}
\begin{proof}
\begin{version3}
$F_g$ is concave so it is continuous over $\mathbb{R}\setminus\{0,1\}$. $F_g$ is continuous in 0 because it is c\`adl\`ag. $F_g$ is continuous in 1 by concavity and monotonicity.
\end{version3}
\end{proof}

\begin{lemma}
Take a real valued sequence $(\lambda_m)$ with $\lambda_m\in(0,1)$, converging to 1, such that $\frac{1}{\sqrt{m}}=o(1-\lambda_m)$ and $\frac{m_{g,0}}{m_g}=\pi_{g,0}+o(1-\lambda_m)$ for all $g$. If $f_g(1^-)=0$ for all $g$ and the $p$-values inside each group are mutually independent, then $$\forall g\in\{1,\dotsc,G\},\,\hat \pi_{g,0}(\lambda_m)\overset{\mathbb{P}}{\longrightarrow} \pi_{g,0}.$$
\label{PAindepgauss}
\end{lemma}

\begin{proof}
First note that $\frac{m_{g,1}}{m_g}-\pi_{g,1}=\pi_{g,0}-\frac{m_{g,0}}{m_g}=o(1-\lambda_m)$.

Thus we have
\begin{align*}
|\hat \pi_{g,0}(\lambda_m)-\pi_{g,0}|&=\left|\frac{1-\frac{1}{m_g}\sum_i\ind{p_{g,i}\leq\lambda_m}+\frac{1}{m} }{1-\lambda_m}-\pi_{g,0}\right|\\
&\leq \frac{\lambda_m\left|\pi_{g,0}- \frac{m_{g,0}}{m_g} \right| + \frac{m_{g,0}}{m_g}\left| \lambda_m-\frac{1}{m_{g,0}}\sum_i \ind{p_{g,i}\leq\lambda_m, H_{g,i}=0}   \right|}{1-\lambda_m} \\
&\quad+\frac{\left|\pi_{g,1}- \frac{m_{g,1}}{m_g} \right|    + \frac{m_{g,1}}{m_g}\left| F_g(\lambda_m)-\frac{1}{m_{g,1}}\sum_i \ind{p_{g,i}\leq\lambda_m, H_{g,i}=1}   \right|}{1-\lambda_m} \\
&\quad+\frac{m_{g,1}}{m_g}\frac{1-F_g(\lambda_m)}{1-\lambda_m}+\frac{1}{m(1-\lambda_m)}\\
&\leq  \frac{m_{g,0}}{m_g}\frac{\sup_{x\in[0,1]}\left| x-\frac{1}{m_{g,0}}\sum_i \ind{p_{g,i}\leq x, H_{g,i}=0}   \right|}{1-\lambda_m} \\
&\quad+\frac{m_{g,1}}{m_g}\frac{\sup_{x\in[0,1]}\left| F_g(x)-\frac{1}{m_{g,1}}\sum_i \ind{p_{g,i}\leq x, H_{g,i}=1}   \right|}{1-\lambda_m} +o(1).\\
\end{align*}
The two suprema of the last display, when multiplied by $\sqrt{m}$, converge in distribution (by Kolmogorov-Smirnov's theorem). So when divided by $1-\lambda_m$ they converge to 0 in distribution and then in probability (because $\frac{1}{1-\lambda_m}=o(\sqrt{m})$).
\end{proof}

\begin{definition}
The critical alpha value is
$$\alpha^*=\inf_{w\in K^\infty}\frac{1}{\sum_g\pi_g w_g\left(\pi_{g,0}+\pi_{g,1}f_g(0^+)\right) },$$
where $K^\infty=\{w\in \mathbb{R}^G_+ : \sum_g \pi_g \bar \pi_{g,0} w_g\leq 1\}$.
\label{def_alpha}
\end{definition}

\begin{lemma} $\alpha^*$ \begin{version2}always satisfies\end{version2} $\alpha^*<1.$
\label{alphacrit}
\end{lemma}
\begin{proof}
We only need to show that for one $w\in K^\infty$, we have 
$$\sum_g\pi_g w_g\left(\pi_{g,0}+\pi_{g,1}f_g(0^+)\right)>1.$$ Let us show that this is true for every $w\in K^\infty$ such that $ \sum_g \pi_g \bar \pi_{g,0} w_g= 1$, e.g. the $w$ defined by $w_g=\frac{1}{\bar \pi_{g,0}}$ for all $g$. We use the fact that $f_g(0^+)>\frac{F_g(1)-F_g(0)}{1-0}=1$ by the strict concavity of $F_g$. Then $\pi_{g,0}+\pi_{g,1}f_g(0^+)>1$ and
\begin{equation*}
\sum_g\pi_g w_g\left(\pi_{g,0}+\pi_{g,1}f_g(0^+)\right)>\sum_g\pi_g w_g\geq\sum_g\pi_g\bar \pi_{g,0} w_g=1.\qedhere
\end{equation*}
\end{proof}

Recall that $\mathcal{I}(\cdot)$ is defined as $\mathcal{I}(h)=\sup\left\{ u\in[0,1] : h(u)\geq u     \right\}$ on the function space:
\begin{equation}
\mathcal{F}=\left\{h:[0,1]\to[0,1] : h(0)=0, \,h\text{ is nondecreasing}   \right\}
\label{def_F}
\end{equation}
which has the natural order $h_1\leq h_2 \iff h_1(u)\leq h_2(u)\, \forall u\in[0,1]$. $\mathcal{F}$ is also normed with the sup norm $\|\cdot\|$.
\begin{lemma}
For all $h\in\mathcal{F}$, $\mathcal{I}(h)$ is a maximum and $h\left(\mathcal{I}(h)\right)=\mathcal{I}(h)$. Moreover, $\mathcal{I}(\cdot)$, seen as a map on $\mathcal{F}$, is nondecreasing and continous on each continuous $h_0\in\mathcal{F}$ such that either $u\mapsto h_0(u)/u$ is decreasing over $(0,1]$, or $\mathcal{I}(h_0)=0$.
\label{IG}
\label{seuil}
\end{lemma}
\begin{proof}
$\mathcal{I}(h)$ is a maximum because there exists $\epsilon_n\to0$ such that
$$h\left(\mathcal{I}(h)\right)\geq h\left(\mathcal{I}(h)-\epsilon_n\right)\geq \mathcal{I}(h)-\epsilon_n\to\mathcal{I}(h).$$

So $h\left(\mathcal{I}(h)\right)\geq\mathcal{I}(h)$. Then
$ h\left(h\left(\mathcal{I}(h)\right)\right) \geq h\left(\mathcal{I}(h)\right)  $
thus $h\left(\mathcal{I}(h)\right)\leq\mathcal{I}(h)$ by the definition of $\mathcal{I}(h)$ as a supremum.

Next, if $h_1\leq h_2$, $\mathcal{I}(h_1)=h_1\left( \mathcal{I}(h_1) \right)\leq h_2\left( \mathcal{I}(h_1) \right)$ so $\mathcal{I}(h_1)\leq \mathcal{I}(h_2)$ by defintion of $\mathcal{I}(h_2)$.

Now take a continuous $h_0\in\mathcal{F}$ such that either $u\mapsto h_0(u)/u$ is decreasing or $\mathcal{I}(h_0)=0$, and $h$ any element of $\mathcal{F}$. Let $\gamma>0$, let $u_-=\mathcal{I}(h_0)-\gamma$ and $u_+=\mathcal{I}(h_0)+\gamma$. We want to prove that there exists an $\eta_\gamma$ such that $\|h-h_0\|\leq\eta_\gamma$ implies $u_-\leq\mathcal{I}(h)\leq u_+$.

If $u_+>1$ then obviously $\mathcal{I}(h)\leq u_+$. If not, let $s_\gamma=\underset{u'\in[u_+,1]}{\max}\left(h_0(u')-u'  \right)$. It is a maximum by continuity over a compact and is such that $s_\gamma<0$, because $s_\gamma\geq0$ would contradict the maximality of $\mathcal{I}(h_0)$.

Then, for all $u'\in[u_+,1]$,
\begin{equation*}
h(u')-u'\leq h_0(u')-u' + \left\|h -h_0   \right\|,
\end{equation*}
and then
\begin{equation*}
\sup_{u'\in[u_+,1]}\left( h(u')-u' \right)\leq s_\gamma+  \left\|h -h_0   \right\|   .\end{equation*}
Hence, as soon as $\left\|h-h_0\right\|\leq\frac{1}{2}|s_\gamma|$, $\sup_{u'\in[u_+,1]}\left(   h(u')-u' \right)<0$ and $\mathcal{I}(h)< u_+$.

If $u_-\leq0$, which is always the case if $\mathcal{I}(h_0)=0$, then $\mathcal{I}(h)\geq u_-$. If $u_->0$, $u\mapsto h_0(u)/u$ is decreasing and
$$\frac{h_0(u_-)}{u_-}> \frac{h_0\left(\mathcal{I}(h_0)\right)}{\mathcal{I}(h_0)}=1,$$
so $h_0(u_-)>u_-$. We can then write the following:
\begin{align*}
h(u_-)-u_-&\geq h_0(u_-)-u_- -  \left\|h-h_0   \right\|    >0,
\end{align*}
as soon as $\left\|h-h_0  \right\|  \leq \frac{1}{2}\left(h_0(u_-)-u_-\right)$. 
This implies $\mathcal{I}(h)>u_-$. Taking
$$\eta_\gamma=\frac{1}{2}\min\left(|s_\gamma|\ind{u_+\leq1}+\ind{u_+>1}, \left( h_0(u_-)-u_-\right)\ind{u_->0}+\ind{u_-\leq0}  \right)  $$
completes the proof.
\end{proof}

\begin{lemma}
Let a weight function $W : [0,1]\to \mathbb{R}^G_+$. For each $r$ between 1 and $m$ denote the $W(r/m)$-weighted $p$-values $p_{g,i}^{[r]}=p_{g,i}/W_g(r/m)$ (with the convention $p_{g,i}/0=\infty$), order them $p_{(1)}^{[r]}\leq\dotsc\leq p_{(m)}^{[r]}$ and note $p_{(0)}^{[r]}=0$.

Then $\hat  u_W= m^{-1}\max\left\{ r\geq0 :  p_{(r)}^{[r]} \leq \alpha \frac{r}{m} \right\}. $
\label{hatuegalIG}
\end{lemma}
\begin{proof}
Let us denote $\hat r = \max\left\{ r\geq0 :  p_{(r)}^{[r]} \leq \alpha \frac{r}{m} \right\}$ and show $\hat u_W = \hat r /m$ by double inequality. First, we have
\begin{align*}
\widehat G_{W}\left(\frac{\hat r}{m}\right)&=m^{-1}\sum_{g=1}^G\sum_{i=1}^{m_g}\mathds{1}_{\left\{p_{g,i} \leq \alpha \frac{\hat r}{m} W_g\left(\frac{\hat r}{m}\right)  \right\}}\\
&=m^{-1}\sum_{g=1}^G\sum_{i=1}^{m_g}\mathds{1}_{\left\{p^{[\hat r ]}_{g,i} \leq \alpha \frac{\hat r}{m} \right\}}\\
&=m^{-1}\sum_{r=1}^m\mathds{1}_{\left\{p_{(r)}^{[\hat r]} \leq \alpha \frac{\hat r}{m}  \right\}}\geq \hat r/m,
\end{align*}
because $p_{(1)}^{[\hat r]} ,\dotsc,p_{(\hat r)}^{[\hat r]} \leq \alpha \frac{\hat r}{m} $. Then $\hat r/m\leq \hat u_W$ by definition of $\hat u_W$. Second, we know that $\hat u_W$ can be written as $\hat \kappa/m$ because $\hat u_W=\widehat G_W(\hat u_W)$, so we want to show that $\hat \kappa\leq \hat r$ which is implied by $\hat r$, $p_{(\hat \kappa)}^{[\hat \kappa]} \leq \alpha \frac{\hat \kappa}{m}$. The latter is true because
$$ \sum_{r=1}^m \mathds{1}_{\left\{ p^{[\hat \kappa]}_{(r)}\leq\alpha \frac{\hat \kappa}{m}   \right\} }= m\widehat G_W\left(\frac{\hat \kappa}{m}\right)=m\widehat G_W\left(\hat u_W\right)\geq \hat \kappa.$$

\end{proof}

\begin{lemma}
$\widehat G_{\widehat W^*}$ is nondecreasing.
\label{croissance}
\end{lemma}
\begin{proof}
Let $u\leq u^\prime$. $\widehat G_{\widehat W^*}(u^\prime )= \underset{w\in \Km}{\max}  \widehat G_{ w}(u^\prime) $ so by denoting $w= \widehat W^*(u)$ we have $\widehat G_{\widehat W^*}(u^\prime )\geq \widehat G_{w}(u^\prime )$. Furthermore, 
\begin{equation*}
\widehat G_{w}(u^\prime )=\frac{1}{m}\sum_{g=1}^G \sum_{i=1}^{m_g} \mathds{1}_{\left\{p_{g,i}\leq \alpha u^\prime w_g\right\}} \geq \frac{1}{m}\sum_{g=1}^G \sum_{i=1}^{m_g} \mathds{1}_{\left\{p_{g,i}\leq \alpha uw_g\right\}} = \widehat G_{\widehat W^*}(u),
\end{equation*}
which entails $\widehat G_{\widehat W^*}(u^\prime )\geq \widehat G_{\widehat W^*}(u)$.
\end{proof}

\section{Asymptotical weighting}
\label{subsection_asymptotical}

Define, for a weight function $W:[0,1]\to \mathbb{R}^G_+$, possibly random,
$$P^\infty_W:u\mapsto\sum_{g=1}^G\pi_g \pi_{g,1}F_g\left( \alpha u W_g(u) \right);$$
$$G^\infty_W:u\mapsto\sum_{g=1}^G\pi_g D_g\left( \alpha u W_g(u) \right);$$
and
$$H^\infty_W(u)=G^\infty_W(u)-P^\infty_W(u),$$
where \begin{equation}
D_g:t\mapsto \pi_{g,0}U(t)+\pi_{g,1}F_g(t)
\label{def_Dg}
\end{equation} is strictly concave on $[0,1]$ because $F_g$ is and $\pi_{g,1}>0$. Note that, if $W$ is a fixed deterministic weight function, $P^\infty_W$ and $G^\infty_W$ are the uniform limits of $P^{(m)}_W$ and $G^{(m)}_W$ when $m\to\infty$. If $W$ is such that $G^\infty_W$ is nondecreasing, we also define
\begin{equation}
u^\infty_W=\mathcal{I}\left(G^\infty_W\right).
\label{def_uwinf}
\end{equation}

Recall that $K^\infty=\{w\in \mathbb{R}^G_+ : \sum_g \pi_g \bar \pi_{g,0} w_g\leq 1\}$. It is the asymptotic version of $\Km$. We now define oracle optimal weights over $K^\infty$ for $G^\infty_\cdot(u)$ and $P^\infty_\cdot(u)$, for all $u>0$. 
\begin{lemma}
Fix an $u\in[0,1]$. Then $\argmax_{w\in K^\infty} G^\infty_w(u)$ is non empty.

If $0<\alpha u \leq \bar \pi_0$, it is a singleton. \begin{version2}In this case, its\end{version2} only element $w^*$ belongs to $[0,\frac{1}{\alpha u}]^G$ and \begin{version2}satisfies\end{version2} $\sum_g\pi_g\bar \pi_{g,0}w^*_g=1$. If $\alpha u \geq \bar \pi_0$ it is included in $[\frac{1}{\alpha u},\infty)^G$.

Finally, \begin{version2}$\max_{w\in K^\infty} G^\infty_w(u)\leq1$ with equality\end{version2} if and only if $\alpha u \geq \bar \pi_0$.

The same statements are true for $P^\infty_\cdot$, except that \begin{version2}the upper bound of $\max_{w\in K^\infty} P^\infty_w(u)$, which is achieved  if and only if $\alpha u \geq \bar \pi_0$, is not $1$ but $1-\pi_0$\end{version2}.
\label{prop_argmax}
\end{lemma}

\begin{proof}
The function $w\mapsto G^\infty_w(u)$ is continuous over the compact $K^\infty$ so it has a maximum. Note that $\max_{w\in K^\infty} G^\infty_w(0)=0$ and $\argmax_{w\in K^\infty} G^\infty_w(0)=K^\infty$. For the rest of the proof $u$ is greater than 0.

First we show that any $w^*\in\argmax_{w\in K^\infty} G^\infty_w(u)$ belongs to $[0,\frac{1}{\alpha u}]^G$ or $[\frac{1}{\alpha u},\infty)^G$. If not, there is $w^*\in\argmax_{w\in K^\infty} G^\infty_w(u)$ such that $\alpha u w^*_{g_1}>1$ and $\alpha u w^*_{g_2}< 1$ for some $g_1,g_2\leq G$. Now then we define $\tilde w$ such that $\tilde w_g=w^*_g$ for all $g\not\in \{ g_1, g_2\}$, $\tilde w_{g_1}=\frac{1}{\alpha u}$ and 
$$\tilde w_{g_2}=  w^*_{g_2}+\left(w^*_{g_1}-\frac{1}{\alpha u}\right)\frac{\pi_{g_1}\bar \pi_{g_1,0}}{\pi_{g_2}\bar \pi_{g_2,0}}>w^*_{g_2}  .$$
So $\tilde w$ belongs to $K^\infty$ and satisfies
\begin{align*}
G^\infty_{\tilde w}(u)&=\sum_{g\neq g_1,g_2}\pi_g D_g(\alpha u w^*_g)+\pi_{g_1}+\pi_{g_2}D_{g_2}(\alpha u \tilde w_{g_2})\\
&>\sum_{g\neq g_1,g_2}\pi_g D_g(\alpha u w^*_g)+\pi_{g_1}+\pi_{g_2}D_{g_2}(\alpha u w^*_{g_2})=G^\infty_{w^*}(u),
\end{align*}
because $D_{g}$ is increasing over $[0,1]$ and then constant equal to 1. This contradicts the definition of $w^*$ so is impossible.

Next we distinct three cases.

\emph{(i)} $\alpha u = \bar \pi_0$. Then $w_0=(\frac{1}{\alpha u},\dotsc,\frac{1}{\alpha u})=(\frac{1}{\bar \pi_0},\dotsc,\frac{1}{\bar \pi_0})$ is obviously an element of $\argmax_{w\in K^\infty} G^\infty_w(u)$ because
\begin{equation*}
G^\infty_{w_0}(u)=\sum_{g=1}^G\pi_g D_g\left( 1 \right)=1,
\end{equation*}
and we easily check that $\sum_g\pi_g\bar \pi_{g,0}(w_{0})_g=1$.
Thus for every $w\in K^\infty$ distinct from $w_0$, there must exist a $g_1\in\{1,\dotsc,G\}$ such that $\alpha u w_{g_1}< 1$, so $D_{g_1}\!\left( \alpha u w_{g_1} \right)<1$ and
$
G^\infty_{w}(u)<\sum_{g}\pi_g=1
$ : $w_0$ is the only element of $\argmax_{w\in K^\infty} G^\infty_w(u)$.

\emph{(ii)} $\alpha u < \bar \pi_0$. If a $w^*\in\argmax_{w\in K^\infty}G^\infty_w(u)$ exists in $[\frac{1}{\alpha u},\infty)^G$, then $w^*_g\geq\frac{1}{\alpha u}>\frac{1}{\bar \pi_{0}}$ and $\sum_g\pi_g\bar \pi_{g,0}w^*_g>1$ which is impossible. So 
$$\underset{w\in K^\infty}{\argmax}\,G^\infty_w(u)=\underset{w\in K^\infty\cap[0,\frac{1}{\alpha u}]^G }{\argmax}G^\infty_w(u).$$
The function $w\mapsto G^\infty_w(u)$ is strictly concave over the convex set $K^\infty\cap[0,\frac{1}{\alpha u}]^G$ because $\pi_{g,1}>0$ and $D_g$ is strictly concave over $[0,1]$ for all $g$, hence the maximum is unique.

We showed that the only $w^*\in \underset{w\in K^\infty}{\argmax}\,G^\infty_w(u)$ is not in $[\frac{1}{\alpha u},\infty)^G$ so there exists $g_1\leq G$ such that $\alpha u w^*_{g_1}<1$ thus $G^\infty_{w^*}(u)<1$. Furthermore $\sum_g\pi_g\bar \pi_{g,0}w^*_g=1$ : if not there \begin{version2}would exist\end{version2} a $\tilde w$ with $\tilde w_{g_1}> w^*_{g_1}$ \begin{version2}(for the same $g_1$ as in previous sentence)\end{version2} and $\tilde w_{g}=w^*_{g}$ for all $g\neq g_1$ such that $\tilde w\in K^\infty$ and $G^\infty_{\tilde w}(u)>G^\infty_{w^*}(u)$ which is impossible.

\emph{(iii)} $\alpha u > \bar \pi_0$. So $u>\frac{\bar \pi_0}{\alpha}$ and obviously 
$$\max_{w\in K^\infty} G^\infty_w(u)\geq \max_{w\in K^\infty} G^\infty_w\left(\frac{\bar \pi_0}{\alpha}\right)=G^\infty_{w_0}\left(\frac{\bar \pi_0}{\alpha}\right) =1,$$
as stated in case \emph{(i)}. So $\max_{w\in K^\infty} G^\infty_w(u)=1$ and the vectors $w^*$ of $\argmax_{w\in K^\infty} G^\infty_w(u)$ are the ones fulfilling $D_g(\alpha u w^*_g)=1$ for all $g$ that is $w^*\in [\frac{1}{\alpha u},\infty)^G$.

The proof is similar for $P^\infty_\cdot$, by replacing $D_g$ by $\pi_{g,1}F_g$.
\end{proof}

From now on, $W^*(u)$ denotes an element of $\argmax_{w\in K^\infty} G^\infty_w(u)$ (just like we write $\widehat W^*(u)$ as an element of $\argmax_{w\in \Km}\widehat G_w(u)$), our results will not depend on the chosen element of the argmax. Next Lemma gives some properties on the function $G^\infty_{W^*}$, among them $G^\infty_{W^*}$ is nondecreasing which allow us to define
\begin{equation}
u^*=u^\infty_{W^*}=\mathcal{I}\left(G^\infty_{W^*}\right).
\label{def_uwinfoptim}
\end{equation}

\begin{lemma}
$G^\infty_{W^*}$ is nondecreasing and $u^*>0$. $G^\infty_{W^*}$ is strictly concave over $[0,\frac{\bar \pi_0}{\alpha}\wedge1]$ and, if $\alpha\geq\bar \pi_0$, constant equal to 1 over $[\frac{\bar \pi_0}{\alpha},1]$. 

In particular, (i) $u^*=1$ if and only if $\alpha\geq\bar \pi_0$, (ii) the function $u\mapsto G^\infty_{W^*}(u)/u$ is decreasing over $(0,1]$, (iii) $G^\infty_{W^*}$ is continuous over $[0,1]$.
\label{prop_max}
\end{lemma}
\begin{proof}
$G^\infty_{W^*}$ is nondecreasing by exactly the same argument as in the proof of Lemma~\ref{croissance}. The result can be strengthened thanks to Lemma~\ref{prop_argmax}, by writing, for $u<u'\leq\frac{\bar \pi_0}{\alpha}\wedge1$, that $G^\infty_{W^*(u)}(u')>G^\infty_{W^*(u)}(u)$ because $1>G^\infty_{W^*}(u)$. So $G^\infty_{W^*}$ is increasing on $[0,\frac{\bar \pi_0}{\alpha}\wedge1]$.

To prove that $u^*>0$, take some $w\in K^\infty$ such that
$$\alpha>\frac{1}{\sum_g\pi_g w_g\left(\pi_{g,0}+\pi_{g,1}f_g(0^+)\right)}   \geq\alpha^* .$$
Because the expression above is continuous of the $w_g$, they can always be chosen nonzero. We have $u^*\geq u^\infty_w$ because $G^\infty_{W^*}\geq G^\infty_{w}$. Then we have, for $x>0$, $x\to0^+$,
\begin{align*}
 \frac{G^{\infty}_{w}(x)-G^{\infty}_{w}(0)}{x-0}=\frac{G^{\infty}_{w}(x)}{x}&=\sum_g\pi_g\pi_{g,0}\alpha w_g +\sum_g\pi_g\pi_{g,1}\alpha w_g\frac{F_g(\alpha x w_g)}{\alpha x w_g}
 \\
 &\to\alpha\sum_g\pi_g w_g\left(\pi_{g,0}+\pi_{g,1}f_g(0^+)\right)>1,
 \end{align*}
so $G^\infty_{w}(u)>u$ in the neighborhood of $0^+$, which entails $u^\infty_{w}>0$.

Now take $a,b\in[0,\frac{\bar \pi_0}{\alpha}\wedge1]$ with $a<b$ and $\lambda\in(0,1)$, by Lemma~\ref{prop_argmax}, \begin{version2}we have that\end{version2} $\alpha aW^*_g(a),\alpha bW^*_g(b)\leq1$ and then, for all $g$:
$$D_g\left(\lambda \alpha aW^*_g(a) +(1-\lambda) \alpha bW^*_g(b)\right)\geq\lambda D_g\left( \alpha aW^*_g(a) \right)+(1-\lambda)D_g\left( \alpha bW^*_g(b) \right).$$
Moreover, because $G^\infty_{W^*}(a)<G^\infty_{W^*}(b)$, for at least one $g_1$ we have $aW^*_{g_1}(a)\neq bW^*_{g_1}(b)$ and by strict concavity of $D_{g_1}$ the inequality above is strict for $g_1$. Then define $\tilde w_g=\frac{\lambda  aW^*_g(a) +(1-\lambda) bW^*_g(b)}{\lambda a +(1-\lambda) b}$. We have $\tilde w\in K^\infty$ and then for all $g$:
$$\pi_gD_g\left(\alpha(\lambda a +(1-\lambda)  b )\tilde w_g\right)\geq\lambda \pi_g D_g\left( \alpha aW^*_g(a) \right)+(1-\lambda)\pi_g D_g\left( \alpha bW^*_g(b) \right),$$
the inequality being strict for $g_1$. Finally by summing:
$$G^\infty_{W^*}(\lambda a +(1-\lambda)  b)\geq G^\infty_{\tilde w}(\lambda a +(1-\lambda)  b) > \lambda G^\infty_{W^*}(a)+(1-\lambda)G^\infty_{W^*}(b) .$$

Additionally, $G^\infty_{W^*}(u)=1$ for $\alpha u\geq\bar \pi_0$ comes from Lemma~\ref{prop_argmax}. The fact that $u^*=1 \iff \alpha\geq\bar \pi_0$ follows directly from the previous statements and Lemma~\ref{prop_argmax}. The decreasingness of $u\mapsto G^\infty_{W^*}(u)/u$ is straightforward from strict concavity properties because it is the slope of the line between the origin and the graph of $G^\infty_{W^*}$ at abscissa $u>0$. Previous statements imply that $G^\infty_{W^*}$ is continuous at least over $(0,\frac{\bar \pi_0}{\alpha}\wedge 1)$ and, if $\alpha\geq\bar \pi_0$, over $[\frac{\bar \pi_0}{\alpha},1]$. $K^\infty$ is bounded, let $B$ such that $|w_g|\leq B$ for all $w\in K^\infty$, then $G^\infty_{W^*}(u)\leq\sum_g\frac{m_g}{m}D_g(\alpha u B)\to0$ when $u\to0$ which gives the continuity in 0. As in the proof of Lemma~\ref{Fgcont}, the continuity in $\frac{\bar \pi_0}{\alpha}\wedge1$ is given by the combination of concavity and nondecreasingness.
\end{proof}

\begin{remark}
The case $\alpha\geq \bar \pi_0$ is rarely met in practice because $\alpha$ is chosen small and the signal is assumed to be sparse (so $\bar \pi_0$ is large) but it is kept to cover all situations. It confirms the intuitive idea that in this situation the best strategy is to reject all hypotheses because then the FDP is equal to $\pi_0\leq\bar \pi_0\leq\alpha$.
\label{rk_grand_alpha}
\end{remark}

\begin{remark}
For a weight vector $w\in \mathbb{R}^G_+$, $G^\infty_w$ is obviously continuous. Moreover if $w\neq0$, let $M=\max_{0\leq u\leq1} G^\infty_w(u)\leq 1$ and $u^\diamond=\min\{u :G^\infty_w(u)=M  \}>0$, then $G^\infty_w$ is strictly concave over $[0,u^\diamond]$ and constant equal to $M$ on $[u^\diamond,1]$, hence $u\mapsto G^\infty_{w}(u)/u$ is decreasing. So whether $w=0$ or not, $\mathcal{I}(\cdot)$ is continuous in $G^\infty_w$ by Lemma~\ref{seuil}. 
\label{weight_vector}
\end{remark}

{
\begin{remark}
The proof of the strict concavity of $G^\infty_{W^*}$ can easily be adapted to show the (non necessary strict) concavity of $\widetilde G_{\widetilde W^*}$ when $\widetilde D_g=\LCM\left(\widehat D_g\right)$.
\label{rk_concav_lcm}
\end{remark}
}

\begin{figure}
\centering
\includegraphics[width=1\linewidth]{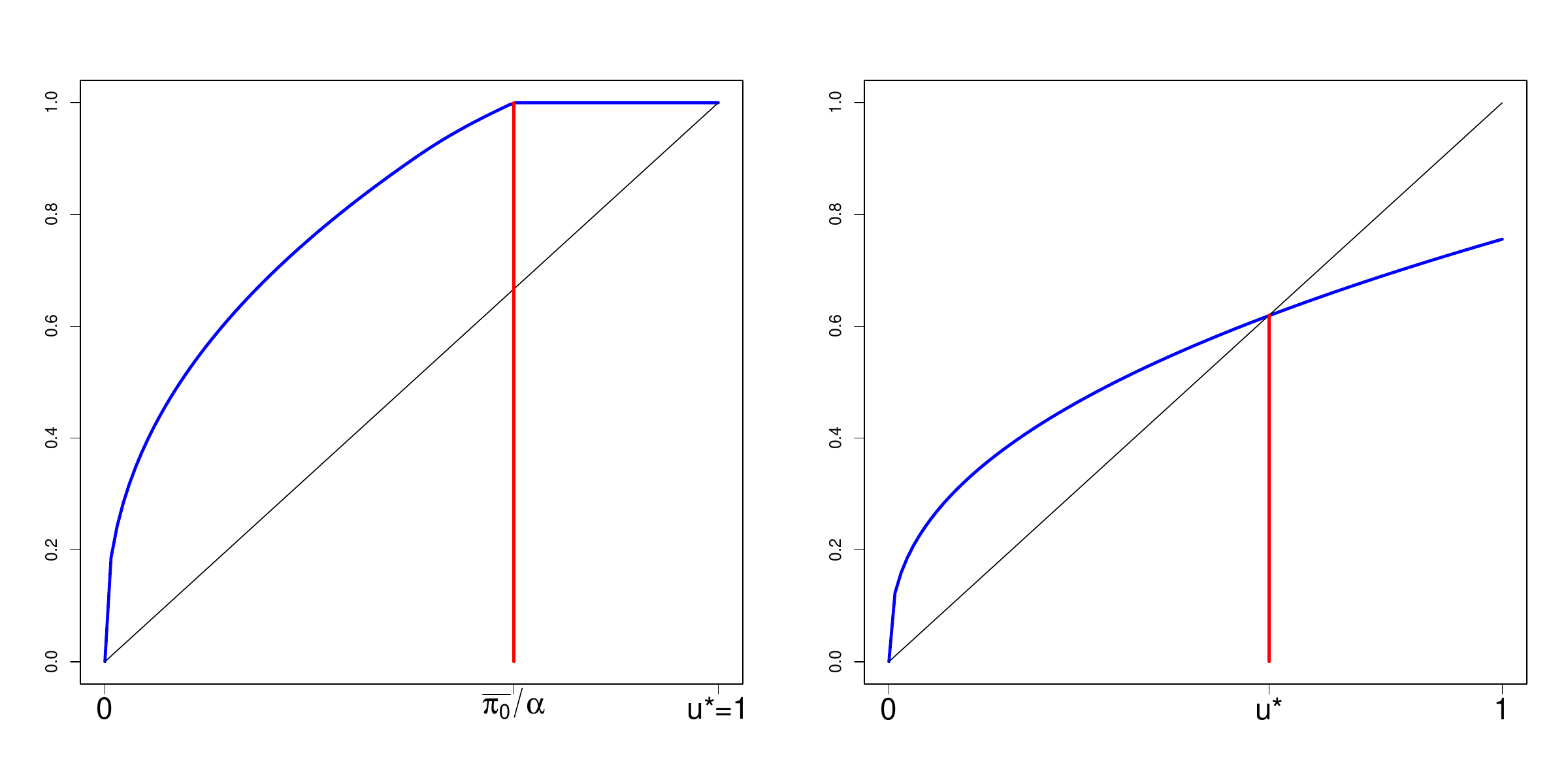}
\caption[Plot of \texorpdfstring{$u\mapsto G^\infty_{W^*}(u)$}{functionG} in both cases]{\footnotesize Plot of $u\mapsto G^\infty_{W^*}(u)$ when $\alpha\geq \bar \pi_0$ (left panel) and $\alpha< \bar \pi_0$ (right panel).}
\label{fig_lmb2}
\end{figure}

Figure~\ref{fig_lmb2} illustrates all the properties stated in Lemma~\ref{prop_max}, with the two cases $\alpha\geq \bar \pi_0$ and $\alpha< \bar \pi_0$.

The next Lemma justifies the intuitive idea that maximizing the rejections and the power is the same thing (as exposed in Section~\ref{subsection_previous}), but only \begin{version2}under~\eqref{ME}.\end{version2}

\begin{lemma}
\begin{version2}If~\eqref{ME} holds\end{version2}, for all $u\in[0,1]$, $$\argmax_{w\in K^\infty} G^\infty_w (u)=\argmax_{w\in K^\infty} P^\infty_w (u).$$
In particular, $P^\infty_{W^*}$ is continuous nondecreasing.
\label{ASTUCE}
\end{lemma}
\begin{proof}
First, $\argmax_{w\in K^\infty} G^\infty_w (0)=\argmax_{w\in K^\infty} P^\infty_w (0)=K^\infty$, so assume $u>0$. If $\alpha u\geq\bar \pi_0$, $\max_{w\in K^\infty} G^\infty_w (u)=1$ and $\max_{w\in K^\infty} P^\infty_w (u)=1-\pi_0$ by Lemma~\ref{prop_argmax}, thus $\argmax_{w\in K^\infty} G^\infty_w (u)$ and $\argmax_{w\in K^\infty} P^\infty_w (u)$ are both equal to the set of weights $w\in K^\infty$ such that $\alpha u w_g\geq 1$ for all $g$.

Now if $\alpha u\leq\bar \pi_0$, both arg max are singletons. Take $w^*$ the only element of $\argmax_{w\in K^\infty} P^\infty_w (u)$. \begin{version2}Recall that there exists $C\geq1$ such that, for all $1\leq g\leq G$, $\bar\pi_{g,0}=C\pi_{g,0}$,\end{version2} and write, for all $w\in K^\infty$,

\begin{version2}
\begin{align*}
G^\infty_w (u)&=\sum_g\pi_g\pi_{g,0} U(\alpha u w_g)+P^\infty_w (u)\\
&\leq \alpha u\sum_g\pi_g\pi_{g,0} w_g+P^\infty_{w^*} (u)\\
&\quad= \frac{\alpha u}{C}\sum_g\pi_g\bar\pi_{g,0} w_g+P^\infty_{w^*} (u)\\
&\quad\leq \frac{\alpha u}{C} \times 1+P^\infty_{w^*} (u)\\
&\quad\quad= \sum_g\pi_g\pi_{g,0}U(\alpha u w^*_g)+P^\infty_{w^*} (u)=G^\infty_{w^*} (u),
\end{align*}
\end{version2}
\begin{version2}because  $\sum_g\pi_g\bar\pi_{g,0}w^*_g=1$ and $\alpha u w^*_g\leq 1$ for all $g$, by Lemma~\ref{prop_argmax}.\end{version2} This means that $w^*$ is also the unique element of $\argmax_{w\in K^\infty} G^\infty_w (u)$. Finally the properties on $P^\infty_{W^*}$ are obtained by the same proof as Lemma~\ref{prop_max}.
\end{proof}

The next lemma is only a deterministic tool used in the proof of Lemma~\ref{chapeauchapeau}. {Define the distance $d$ of a weight vector $w$ to a subset $S$ of $\mathbb{R}^G_+$ by $d( w,S)=\inf_{\bar w\in S}\max_g|w_g-\bar w_g|$. Let $M_u=\argmax_{w\in K^\infty} G^\infty_{ w}(u) $ to lighten notations.}

\begin{lemma}
Take some $u\in (0,1]$. Then we have:
{
\begin{equation*}
\forall \epsilon>0, \exists  \xi >0, \forall  w\in K^\infty , \left| G^\infty_{ w}(u) -G^\infty_{W^*}(u)\right|\leq\xi \Rightarrow d\left( w,M_u  \right)   {<} \epsilon   .
\end{equation*}
}
{In particular, if} $\alpha u\leq\bar \pi_0$,
\begin{equation}
\forall \epsilon>0, \exists  \xi >0, \forall w\in K^\infty , \left| G^\infty_{w}(u) -G^\infty_{W^*}(u)\right|\leq\xi \Rightarrow \max_g\left|w_g-W^*_g(u)\right| {<} \epsilon  {,}
\label{cas1}
\end{equation}
{and if} $\alpha u\geq\bar \pi_0$,
\begin{equation}
\forall \epsilon>0, \exists  \xi >0, \forall w\in K^\infty , \left| G^\infty_{w}(u) -G^\infty_{W^*}(u)\right|\leq\xi \Rightarrow \left(\forall g,\, \alpha uw_g {>}1-\epsilon \right)   .
\label{cas2}
\end{equation}
\label{lm_cv_w}
\end{lemma}
\begin{proof}
{
If the statement is false, there exists some $\epsilon>0$ and some sequence $(w_n)_{n\geq1}$ converging to a $w^\ell$ in $K^\infty$ (because $K^\infty$ is compact), such that $d\left(w_n,M_u \right)\geq \epsilon $ and 
$$\left| G^\infty_{w_n}(u)  -  G^\infty_{W^*}(u)     \right|\to0.$$
By continuity of $D_g$, $G^\infty_{w^\ell}(u)=G^\infty_{W^*}(u)$ so $w^\ell\in M_u  $ which contradicts $d\left(w^\ell,M_u \right)\geq \epsilon $. If $\alpha u\leq\bar \pi_0$, $M_u$ is a singleton by Lemma~\ref{prop_argmax}, hence~\eqref{cas1}. However, if $\alpha u\geq\bar \pi_0$, $M_u=\{w\in\begin{version2}K^\infty\end{version2} : \alpha u w_g\geq1 \,\forall g \}$ by Lemma~\ref{prop_argmax}, hence~\eqref{cas2}. 
}
\end{proof}

\section{Convergence lemmas}
\label{subsection_technic}

Recall that $\|\cdot\|$ is the sup norm for the bounded functions on their definition domain: $\|f\|=\sup_{u\in[0,1]}|f(u)|$ or $\|f\|=\sup_{t\in\mathbb{R}}|f(t)|$.

\begin{lemma}
The following quantities converge to 0 \begin{version3}in probability\end{version3}:\\
$\sup_{w\in \mathbb{R}_+^G } \left\| \widehat H_w-H^\infty_w    \right\|$, $\sup_{w\in  \mathbb{R}_+^G } \left\| \widehat P_w-P^\infty_w    \right\|$, $\sup_{w\in  \mathbb{R}_+^G } \left\| \widehat G_w-G^\infty_w    \right\|$, and $\left\|\widehat D_g-D_g  \right\|$, for all $g\in\{1,\dotsc,G\}$.

Furthermore, for any $(\widetilde D_g)_g$ such that $\left\|\widetilde D_g-D_g  \right\|\overset{\mathbb{P}}{\longrightarrow}0$, 
\begin{equation}
\sup_{w\in  \mathbb{R}_+^G } \left\| \widetilde G_w-G^\infty_w    \right\| \overset{\mathbb{P}}{\longrightarrow}0.
\label{eq_cv_p_g}
\end{equation}
\label{supas}
\end{lemma}
\begin{proof}
By using the same proof as the one of the Glivenko-Cantelli theorem \begin{version3}(which adapts trivially to convergence in probability instead of almost surely)\end{version3}, we get from~\eqref{wd0} and~\eqref{wd1} that, for all $g$, 
$$\left\|\frac{1}{m_{g,0}} \sum_{i=1}^{m_g}\ind{p_{g,i}\leq \cdot, H_{g,i}=0}-U  \right\|\overset{\begin{version3}\Pb\end{version3}}{\longrightarrow}0,$$
and
$$ \left\|\frac{1}{m_{g,1}} \sum_{i=1}^{m_g}\ind{p_{g,i}\leq \cdot, H_{g,i}=1}-F_g \right\|  \overset{\begin{version3}\Pb\end{version3}}{\longrightarrow}0.$$
Next, we write that 
\begin{align*}
\left\|\frac{1}{m_{g}} \sum_{i=1}^{m_g}\ind{p_{g,i}\leq \cdot, H_{g,i}=0}-\pi_{g,0}U  \right\|&\leq \left|\frac{m_{g,0}}{m_g}-\pi_{g,0}\right| \\&\quad+ \pi_{g,0}\left\| \frac{1}{m_{g,0}} \sum_{i=1}^{m_g}\ind{p_{g,i}\leq \cdot, H_{g,i}=0}-U   \right\|   \\
&\overset{\begin{version3}\Pb\end{version3}}{\longrightarrow}0,
\end{align*}
and similarly $\left\|\frac{1}{m_{g}} \sum_{i=1}^{m_g}\ind{p_{g,i}\leq \cdot, H_{g,i}=1}-\pi_{g,1}F_g  \right\|\overset{\begin{version3}\Pb\end{version3}}{\longrightarrow}0$. So by summing, $\left\|\widehat D_g-D_g  \right\|\overset{\begin{version3}\Pb\end{version3}}{\longrightarrow}0$. Apply the triangular inequality once again to get $\left\|\frac{1}{m} \sum_{i=1}^{m_g}\ind{p_{g,i}\leq \cdot, H_{g,i}=0}-\pi_g\pi_{g,0}U  \right\|\overset{\begin{version3}\Pb\end{version3}}{\longrightarrow}0$, which implies
\begin{align*}
\sup_{w\in \mathbb{R}_+^G } \left\| \widehat H_w-H^\infty_w    \right\| &\leq \sum_{g=1}^G  \left\|\frac{1}{m} \sum_{i=1}^{m_g}\ind{p_{g,i}\leq \cdot, H_{g,i}=0}-\pi_g\pi_{g,0}U  \right\|\\
 &\overset{\begin{version3}\Pb\end{version3}}{\longrightarrow} 0.
\end{align*}
Similarly $\sup_{w\in  \mathbb{R}_+^G } \left\| \widehat P_w-P^\infty_w    \right\|  \overset{\begin{version3}\Pb\end{version3}}{\longrightarrow}0$ and $\sup_{w\in  \mathbb{R}_+^G } \left\| \widehat G_w-G^\infty_w    \right\|  \overset{\begin{version3}\Pb\end{version3}}{\longrightarrow}0$ by \begin{version3}summation\end{version3}.

Finally, 
\begin{equation*}
\sup_{w\in  \mathbb{R}_+^G } \left\| \widetilde G_w-G^\infty_w    \right\|\leq \sum_g \left( \left|\frac{m_g}{m}-\pi_g \right|  +\pi_g \left\|\widetilde D_g-D_g  \right\|   \right)\overset{\mathbb{P}}{\longrightarrow}0.\qedhere
\end{equation*}
\end{proof}

From now on $\widetilde D_g$ is assumed to converge uniformly to $D_g$ in probability and that $\widetilde W^*(u)\in\argmax_{w\in \Km}\widetilde G_w(u)$ exists for all $u$.

Next Lemma is the main technical one (with the longest proof).
\begin{lemma}
We have the following convergence in probability:
\begin{equation}
\left\|\widetilde G_{\widetilde W^*}-G^\infty_{W^*}\right\|\overset{\mathbb{P}}{\longrightarrow}0.\notag
\label{eq_supG_hat}
\end{equation}
\label{supG}
\end{lemma}
\begin{proof}

First,
\begin{align*}
\left\|  \widetilde G_{\widetilde W^*} -G^\infty_{W^*}   \right\|&\leq\sup_{w\in\mathbb{R}_+^G}\left\|  \widetilde G_{ w} -G^\infty_{w}    \right\| +\left\|  G^\infty_{ \widetilde W^*} -G^\infty_{W^*}  \right\|,
\end{align*}
where the first term tends to 0 by~\eqref{eq_cv_p_g}, so we work on the second term.

The main idea is to use the maximality of $\widetilde G_w(u)$ in $\widetilde W^*(u)$ and the maximality of $G^\infty_w(u)$ in $W^*(u)$. The problem is that one is a maximum over $\Km$ and the other is over $K^\infty$. The solution consists in defining small variations of $\widetilde W^*(u)$ and $W^*(u)$ to place them respectively in $K^\infty$ and $\Km$.

Let $\widetilde W^\dagger_g(u)=\frac{m_g\hat\pi_{g,0}}{m\pi_g\bar \pi_{g,0}} \widetilde W^*_g(u)$. Then $\widetilde W^\dagger(u)\in K^\infty$ and
\begin{align*}
\left\|\widetilde W^\dagger_g-\widetilde W^*_g\right\|&=\left|\frac{m_g\hat\pi_{g,0}}{m\pi_g\bar \pi_{g,0}}-1\right|\left\|\widetilde W^*_g \right\| \\
&\leq \left|\frac{m_g\hat\pi_{g,0}}{m\pi_g\bar \pi_{g,0}}-1\right| \frac{m}{m_g\hat\pi_{g,0}} \overset{\mathbb{P}}{\longrightarrow} 0 \text{ because }  \frac{m_g}{m}\hat\pi_{g,0}\overset{\mathbb{P}}{\longrightarrow}  \pi_g\bar \pi_{g,0},
\end{align*}
which in turn implies that
\begin{align}
\left\|  G^\infty_{\widetilde W^\dagger}  -G^\infty_{ \widetilde W^*}     \right\|&\leq\sum_g \pi_g \sup_u\left|D_g\left(\alpha u \widetilde W^\dagger_g(u)\right)  -D_g\left(\alpha u  \widetilde W^*_g(u)\right)      \right|\notag\\
&\overset{\mathbb{P}}{\longrightarrow} 0,
\label{eqvague}
\end{align}
because $D_g$ is uniformly continuous over $\mathbb{R}_+$. Likewise, we define $W^\dagger_g(u)=\frac{m\pi_g\bar \pi_{g,0}}{m_g\hat\pi_{g,0}} W^*_g(u)$. Therefore $W^\dagger(u)\in \Km$, 
$$\left\| W^\dagger_g-W^*_g\right\|\leq  \left|\frac{m\pi_g\bar \pi_{g,0}}{m_g\hat\pi_{g,0}}-1\right|\frac{1}{\pi_g\bar \pi_{g,0}} \overset{\mathbb{P}}{\longrightarrow} 0 ,$$
and
\begin{equation}
\left\|  G^\infty_{W^\dagger}  -G^\infty_{W^*}    \right\|\leq\sum_g \pi_g \sup_u\left|D_g\left(\alpha u  W^\dagger_g(u)\right)  -D_g\left(\alpha u W^*_g(u)\right)      \right|\overset{\mathbb{P}}{\longrightarrow} 0. \label{eqvagueprime}
\end{equation}

With~\eqref{eq_cv_p_g} and~\eqref{eqvague}, we deduce that
\begin{align}
\left\|\widetilde G_{\widetilde W^\dagger} - \widetilde G_{ \widetilde W^*}\right\|&\leq \left\|\widetilde G_{\widetilde W^\dagger} - G^\infty_{\widetilde W^\dagger}\right\|
+ \left\|G^\infty_{\widetilde W^\dagger} - G^\infty_{ \widetilde W^*}\right\| \notag\\
&\quad+ \left\|G^\infty_{ \widetilde W^*} - \widetilde G_{ \widetilde W^*}\right\|   \notag\\
&\overset{\mathbb{P}}{\longrightarrow} 0,
\label{chapeauchapeauvague}
\end{align}
and likewise with~\eqref{eq_cv_p_g} and~\eqref{eqvagueprime} we have
\begin{align}
\left\|\widetilde G_{W^\dagger} - \widetilde G_{W^*}\right\|\overset{\mathbb{P}}{\longrightarrow} 0.
\label{chapeauvague}
\end{align}

Combining~\eqref{eq_cv_p_g},~\eqref{eqvague},~\eqref{chapeauchapeauvague},~\eqref{chapeauvague}, and the maximalities of $\widetilde G_{\widetilde W^*}(u)$ and $G^\infty_{ W^*}(u)$ will finish the proof. As a start, write
\begin{align*}
\left\|  G^\infty_{\widetilde W^*} -G^\infty_{W^*}    \right\| &\leq   \left\|  G^\infty_{ \widetilde W^*}  -G^\infty_{\widetilde W^\dagger}     \right\|  
+\left\|  G^\infty_{\widetilde W^\dagger}  -G^\infty_{W^*}     \right\|,
\end{align*}
with $ \left\|  G^\infty_{ \widetilde W^*}  -G^\infty_{\widetilde W^\dagger}     \right\|  \overset{\mathbb{P}}{\to}0$ by~\eqref{eqvague}, and, for all $u$,
$$\left|  G^\infty_{\widetilde W^\dagger}(u)  -G^\infty_{W^*}(u)     \right| =G^\infty_{W^*}(u)- G^\infty_{\widetilde W^\dagger}(u) ,$$
 by maximality of $G^\infty_{W^*}(u)$ over $K^\infty$. Then

\begin{align*}
\sup_u\left(G^\infty_{W^*}(u)- G^\infty_{\widetilde W^\dagger}(u)  \right) &\leq     \sup_u\left(G^\infty_{W^*}(u)-  \widetilde G_{W^*}(u)   \right) \\
&\quad+ \sup_u \left(\widetilde G_{W^*}(u)- \widetilde G_{\widetilde W^\dagger}(u)\right)\\
&\quad+\sup_u\left(  \widetilde G_{\widetilde W^\dagger}(u)  -    G^\infty_{\widetilde W^\dagger}(u)           \right),
\end{align*}
with $\sup_u\left(G^\infty_{W^*}(u)-  \widetilde G_{W^*}(u)   \right) \overset{\mathbb{P}}{\to}0$ and $\sup_u\left(  \widetilde G_{\widetilde W^\dagger}(u)  -    G^\infty_{\widetilde W^\dagger}(u)           \right) \overset{\mathbb{P}}{\to}0$ by~\eqref{eq_cv_p_g}.

Finally,
\begin{align*}
 \sup_u \left(\widetilde G_{W^*}(u)- \widetilde G_{\widetilde W^\dagger}(u)\right)&\leq  \sup_u \left(\widetilde G_{W^*}(u)-\widetilde G_{W^\dagger}(u)\right)   \\
&\quad+    \sup_u \left(\widetilde G_{W^\dagger}(u)- \widetilde G_{ \widetilde W^*}(u)\right)            \\
&\quad+              \sup_u \left( \widetilde G_{ \widetilde W^*}(u)-\widetilde G_{\widetilde W^\dagger}(u)\right),
\end{align*}
with $  \sup_u \left(\widetilde G_{W^*}(u)-\widetilde G_{W^\dagger}(u)\right)   \overset{\mathbb{P}}{\to}0$ \eqref{chapeauvague} and $   \sup_u \left( \widetilde G_{ \widetilde W^*}(u)-\widetilde G_{\widetilde W^\dagger}(u)\right)\overset{\mathbb{P}}{\to}0$ \eqref{chapeauchapeauvague}. As a consequence there exists a random variable $V_m\overset{\mathbb{P}}{\to}0$ such that
$$\left\| \widetilde G_{\widetilde W^*}  -G^\infty_{W^*}  \right\|  \leq   \sup_u \left(\widetilde G_{W^\dagger}(u)- \widetilde G_{ \widetilde W^*}(u)\right)   +  V_m   ,$$
but $\widetilde G_{W^\dagger}(u)- \widetilde G_{ \widetilde W^*}(u) \leq 0$ by maximality of $\widetilde G_{\widetilde W^*}(u) $ over $\Km$, so
$$\left\|  \widetilde G_{\widetilde W^*}  -G^\infty_{W^*}  \right\|  \leq V_m \overset{\mathbb{P}}{\to}0 .  \qedhere$$
\end{proof}

Next Lemma is a direct application of Lemma~\ref{seuil}. Recall that $u^*=u^\infty_{W^*}$ (see~\eqref{def_uwinfoptim}) and let
\begin{equation}
\tilde u=\tilde u_{\widetilde W^*}=\mathcal{I}\left(\widetilde G_{\widetilde W^*}  \right)
\label{def_uwtilde}
\end{equation}

\begin{lemma}
We have the following convergences in probability:
\begin{equation*}
\left\{
\begin{array}{rcl}
\tilde u &\overset{\mathbb{P}}{\longrightarrow}&  u^*\\
\widetilde G_{\widetilde W^*}(\tilde u) &\overset{\mathbb{P}}{\longrightarrow} &G^\infty_{W^*}(u^*).
\end{array}
\right.
\end{equation*}
\label{BOUM}
\end{lemma}
\begin{proof} $u\mapsto G^\infty_{W^*}(u)/u$ is nondecreasing and $G^\infty_{W^*}$ is continuous by Lemma~\ref{prop_max} so by Lemma~\ref{seuil} $\mathcal{I}(\cdot)$ is continuous in $G^\infty_{W^*}$: let $\gamma>0$ and $\eta_\gamma$ as in the proof of Lemma~\ref{seuil}, then
\begin{align*}
\Pro{|\tilde u-u^*|\leq\gamma}\geq\Pro{\left\|\widetilde G_{\widetilde W^*}-G^\infty_{W^*}\right\|\leq\eta_\gamma}\underset{\text{Lemma~\ref{supG}}}{\overset{}{\longrightarrow}}1.
\end{align*}

Second result follows immediately because $\widetilde G_{\widetilde W^*}(\tilde u)=\tilde u$ and $G^\infty_{W^*}(u^*)=u^*$ by Lemma~\ref{IG}.
\end{proof}

\begin{lemma}~\par
\emph{(i)} If $\alpha \leq\bar \pi_0$, $\widetilde W^*(\tilde u)\overset{\mathbb{P}}{\longrightarrow} W^*(u^*)$.

\emph{(ii)} If $\alpha \geq\bar \pi_0$, the inferior limit in probability of $\alpha \tilde u \widetilde W_g(\tilde u)$ is greater than or equal to 1, uniformly in $g$, which reads formally:
$$\forall \epsilon>0,\, \Pro{\forall g,\, \alpha \tilde u \widetilde W^*_g(\tilde u){>}1-\epsilon}\longrightarrow1.$$
\label{chapeauchapeau}
\end{lemma}
\begin{proof}
First, we use the same trick as in the proof of Lemma~\ref{supG}: let $\widetilde W^\dagger_g(u)=\frac{m_g\bar \pi_{g,0}}{m\pi_g\bar \pi_{g,0}} \widetilde W^*_g(u)$ such that $\widetilde W^\dagger(u)\in K^\infty$ and $\|\widetilde W^*_g-\widetilde W^\dagger_g  \| \overset{\mathbb{P}}{\longrightarrow} 0$. 

Let us show that $   \left|  G^\infty_{\widetilde W^\dagger(\tilde u)}(u^*) -G^\infty_{W^*}(u^*)  \right|  \overset{\mathbb{P}}{\longrightarrow} 0$  to apply then Lemma~\ref{lm_cv_w} (always possible because $u^*>0$). We have
\begin{align*}
 \left|  G^\infty_{\widetilde W^\dagger(\tilde u)}(u^*) -G^\infty_{W^*}(u^*)  \right|  &\leq   \left|   G^\infty_{\widetilde W^\dagger(\tilde u)}(u^*)  -G^\infty_{\widetilde W^*}(\tilde u)\right|   \\
  &\quad+  \left|  G^\infty_{\widetilde W^*}(\tilde u)-\widetilde G_{\widetilde W^*}(\tilde u)  \right| \\
  &\quad+ \left| \widetilde G_{\widetilde W^*}(\tilde u)  -G^\infty_{W^*}(u^*)   \right|.
\end{align*}
First term converges to 0 because for all $g$, $D_g$ is uniformly continuous and
\begin{align}
 \left|\alpha u^* \widetilde W^\dagger_g(\tilde u)-\alpha \tilde u  \widetilde W^*_g(\tilde u)    \right| &\leq \left|\alpha u^* \widetilde W^\dagger_g(\tilde u)-\alpha u^*  \widetilde W^*_g(\tilde u)    \right|+\left|\alpha u^* \widetilde W^*_g(\tilde u)-\alpha \tilde u  \widetilde W^*_g(\tilde u)    \right|   \notag\\
 &\leq  \|\widetilde W^\dagger_g-\widetilde W^*_g  \|  +\left|  u^* - \tilde u \right|\frac{m}{m_g\hat \pi_{g,0}}\overset{\mathbb{P}}{\longrightarrow}0.
\label{cv_inter}
\end{align}
Apply~\eqref{eq_cv_p_g} to the second term and Lemma~\ref{BOUM} to the third.

\emph{(i)} If $\alpha \leq\bar \pi_0$, then $\alpha u^* \leq\bar \pi_0$ and by {equation~\eqref{cas1},} $\widetilde W^\dagger(\tilde u)\overset{\mathbb{P}}{\longrightarrow} W^*(u^*)$. But for all $g$
\begin{align*}
\left|\widetilde W^*_g( \tilde u) - W^*_g(u^*) \right|\leq  \|\widetilde W^*_g-\widetilde W^\dagger_g  \|   + \left| \widetilde W^\dagger_g(\tilde u)- W^*_g(u^*) \right|,
\end{align*}
and then $\widetilde W^*(\tilde u)\overset{\mathbb{P}}{\longrightarrow} W^*(u^*)$.

\emph{(ii)} If $\alpha \geq\bar \pi_0$, $u^*=1$ by Lemma~\ref{prop_max} and by {equation~\eqref{cas2},}
$$\forall \epsilon>0,\, \Pro{\forall g,\, \alpha u^* \widetilde W^\dagger_g(\tilde u){>}1-\frac{\epsilon}{2}}\longrightarrow1.$$
By equation~\eqref{cv_inter} we also have
$$\forall \epsilon>0,\, \Pro{\forall g,\,  \left|\alpha u^* \widetilde W^\dagger_g(\tilde u)-\alpha \tilde u  \widetilde W^*_g(\tilde u)    \right|\leq\frac{\epsilon}{2}}\longrightarrow1,$$
and by combining the two we get the desired result.
\end{proof}

\begin{lemma}
We have the following convergences in probability:
\begin{equation*}
\widehat G_{\widetilde W^*}(\tilde u)\overset{\mathbb{P}}{\longrightarrow}G^\infty_{W^*}(u^*),
\end{equation*}
\begin{equation*}
\widehat H_{\widetilde W^*}(\tilde u)\overset{\mathbb{P}}{\longrightarrow}H^\infty_{W^*}(u^*).
\end{equation*}
\label{BOUMH}
\end{lemma}
\begin{proof}
We have 
\begin{align*}
\left| \widehat G_{\widetilde W^*}(\tilde u)-G^\infty_{W^*}(u^*)    \right|  &\leq \sup_{w\in\mathbb{R}_+^G} \left\|  \widehat G_{w}-G^\infty_{w}    \right\| +\left|G^\infty_{\widetilde W^*}(\tilde u) -G^\infty_{W^*}(u^*)     \right|.
\end{align*}
Hence, by Lemma~\ref{supas}, we only need to show that $\begin{version3}G\end{version3}^\infty_{\widetilde W^*}(\tilde u)\overset{\mathbb{P}}{\longrightarrow}\begin{version3}G\end{version3}^\infty_{W^*}(u^*)$.

\emph{(i)} If $\alpha \leq\bar \pi_0$, $\tilde u \overset{\mathbb{P}}{\longrightarrow}  u^*$ and $\widetilde W^*(\tilde u) \overset{\mathbb{P}}{\longrightarrow}  W^*(u^*)$ by Lemma~\ref{chapeauchapeau}. Then $\alpha \tilde u \widetilde W^*(\hat u) \overset{\mathbb{P}}{\longrightarrow} \alpha u^* W^*(u^*) $. We get the desired convergence by $D_g$'s continuity.

\emph{(ii)} If $\alpha \geq\bar \pi_0$, $u^*=1$ and $\alpha u^*W^*_g(u^*)\geq1$ for all $g$ so $G^\infty_{W^*}(u^*)=1$. Then by Lemma~\ref{chapeauchapeau} $D_g\left(\alpha\tilde u\widetilde W^*_g(\hat u)\right)\overset{\mathbb{P}}{\longrightarrow}1$ which means that $G^\infty_{\widetilde W^*}(\tilde u)\overset{\mathbb{P}}{\longrightarrow}\sum_g\pi_g1=1$.

The proof for $\hat H$ is similar, just replace $D_g$ by $\pi_{g,0}U$.
\end{proof}

{
The last lemma states that $\LCM(\widehat D_g)$ is a valid estimator of $D_g$ to use in GADDOW.
\begin{lemma}
Assume that $\widetilde D_g=\LCM(\widehat D_g)$. Then $\widetilde D_g$ is nondecreasing, $\widetilde D_g(0)=0$, $\widetilde D_g(1)=1$ and $\left\|\widetilde D_g-D_g  \right\|\overset{\mathbb{P}}{\longrightarrow}0$.
\label{lm_lcm}
\end{lemma}
\begin{proof}
$\widetilde D_g(0)=\widehat D_g(0)=0$ and $\widetilde D_g(1)=\widehat D_g(1)=1$ from the closed form given in Lemma~1 in \citet{carolan2002least}. Let $a,b\in[0,1]$, $a<b$, and let 
\begin{equation*}
C(t)=\left\{\begin{array}{rcl}
\widetilde D_g(t+b-a)  &  \text{if} & t+b-a\leq1  \\
 1 & \text{if}  &   t+b-a\geq1 .
\end{array}\right.
\end{equation*}
Then,
\begin{equation*}
C(t)\geq\left\{\begin{array}{rcccl}
\widehat D_g(t+b-a)&\geq &\widehat D_g(t)  &  \text{if} & t+b-a\leq1  \\
 1&\geq &\widehat D_g(t) & \text{if}  &   t+b-a\geq1 ,
\end{array}\right.
\end{equation*}
because $\widehat D_g$ is non decreasing. 

\begin{version3}Furthermore, C is concave. The inequality $C(\lambda x+(1-\lambda)y)\geq \lambda C(x)+(1-\lambda)C(y)$ for $x, y, \lambda\in[0,1]$, $x\leq y$, is trivial except in the case where $\lambda x+(1-\lambda)y\leq A\leq y$, where $A=1-(b-a)$. In this case, let $\lambda'$ such that $\lambda x+(1-\lambda)y=\lambda' x+(1-\lambda')A$. Then 
\begin{align*}
C(\lambda x+(1-\lambda)y)&=C(\lambda' x+(1-\lambda')A)\\
&\geq \lambda' C(x)+(1-\lambda') C(A)\\
&\quad=\lambda' C(x)+(1-\lambda')\\
&\quad\geq \lambda C(x)+(1-\lambda)= \lambda C(x)+(1-\lambda)C(y)
\end{align*}
where the last inequality is true because $\lambda'\leq \lambda$ and $C(x)\leq1$. Indeed, we have
\begin{align*}
\lambda'&=\frac{\lambda x+(1-\lambda)y-A}{x-A}\\
&=\frac{\lambda(y-x)-(y-A)}{A-x}\\
&=\frac{\lambda(A-x)+\lambda(y-A)-(y-A)}{A-x}\\
&=\lambda-(1-\lambda)\frac{y-A}{A-x}\leq\lambda.
\end{align*}
\end{version3} 

So by definition of the LCM, $C(t)\geq \widetilde D_g(t)$ for all $t\in[0,1]$. In particular,
\begin{equation*}
\widetilde D_g(b)=C(a)\geq \widetilde D_g(a),
\end{equation*}
and $\widetilde D_g$ is nondecreasing. Finally, the convergence comes from $\|\widetilde D_g-D_g  \|\leq\|\widehat D_g-D_g  \|$, see also \citet{carolan2002least}. 
\end{proof}
}

\section{Proof of Corollary~\ref{cor_zz} for $\ZZPro1$}
\label{subsection_proof_pro1}
First, $\hat w^{(1)}\overset{\mathbb{P}}{\longrightarrow}w^{(1)}$ where $w^{(1)}=\left(\frac{1}{\bar \pi_0}   , \dotsc,\frac{1}{\bar \pi_0}  \right)$ and $ \hat w^{(2)}\overset{\mathbb{P}}{\longrightarrow} w^{(2)}$ where, for all $g$, $w^{(2)}_g=\frac{\bar \pi_{g,1}}{\bar \pi_{g,0}(1-\bar \pi_0)}$. By using Lemma~\ref{supas} and the continuity of $D_g$, we get that $\|\widehat G_{\hat w^{(1)}}-G^\infty_{w^{(1)}}\|\overset{\mathbb{P}}{\longrightarrow}0$ and $\|\widehat G_{\hat w^{(2)}}-G^\infty_{w^{(2)}}\|\overset{\mathbb{P}}{\longrightarrow}0$ and then by Lemma~\ref{seuil} we get that $\hat u_{\hat w^{(1)}}\overset{\mathbb{P}}{\longrightarrow} u^\infty_{ w^{(1)}}$ and $\hat u_{\hat w^{(2)}}\overset{\mathbb{P}}{\longrightarrow} u^\infty_{ w^{(2)}}$ so $\hat u_M \overset{\mathbb{P}}{\longrightarrow} u_M$ where $u_M=\max( u^\infty_{ w^{(1)}},u^\infty_{ w^{(2)}} )$.

Define again $\widehat W^\dagger_g(u)=\frac{m_g\hat\pi_{g,0}}{m\pi_g\bar \pi_{g,0}}\widehat W^*_g(u)$ and note that the power of Pro1 is $\Esp{\widehat P_{\widehat W^*}(\hat u_M)}$. We have
\begin{align*}
\widehat P_{\widehat W^*}(\hat u_M)&\leq \sup_{w\in \mathbb{R}^G_+}\left\|\widehat P_w-P^\infty_w\right\|+\left\| P^\infty_{\widehat W^*}-P^\infty_{\widehat W^\dagger} \right\| +P^\infty_{\widehat W^\dagger}(\hat u_M)\\
&\leq  \sup_{w\in \mathbb{R}^G_+}\left\|\widehat P_w-P^\infty_w\right\|+\left\| P^\infty_{\widehat W^*}-P^\infty_{\widehat W^\dagger} \right\| +P^\infty_{W^*}(\hat u_M)   \\
&\quad\overset{\mathbb{P}}{\longrightarrow} P^\infty_{W^*}( u_M),
\end{align*}
because $P^\infty_{W^*}$ is continuous by Lemma~\ref{ASTUCE}.

Note that $u^*\geq u_M$ (because $G^\infty_{W^*}\geq G^\infty_{w^{(1)}}$ and $G^\infty_{W^*}\geq G^\infty_{w^{(2)}}$) to conclude.

\begin{version3}

\section{Proof of Theorem~\ref{thm_cross}}
\label{section_proof_cross}
First, note that, by the independence provided by Assumption~\ref{dep_folds}, we can work conditionally to $(\mathbb{F}_m)_m$ and consider this sequence fixed and deterministic. Second, note that $m^{-1}| \{1\leq i\leq m_g : \mathbb{F}_m(g,i) =f \} |   \underset{m\to\infty}{\longrightarrow}\frac{\pi_g}{F}$. 
Finally, note that the empirical function used for $\crossADDOW$ is defined by
\begin{equation}
 \widehat G^{\cross} : u \mapsto m^{-1}\sum_{g=1}^G\sum_{f=1}^F\sum_{\substack{1\leq i\leq m_g:\\\mathbb{F}_m(g,i)= f}}\mathds{1}_{\{ p_{g,i}\leq \alpha u w^*_{g,f} \}} .\label{eq:Gcross}
\end{equation}

By using Assumption~\ref{weak_dep_folds}, we can apply Lemmas~\ref{supas}, \ref{supG}, \ref{BOUM} and~\ref{chapeauchapeau} to $\ADDOW_{-f}$ for each $f$, so we get that $w^*_{g,f} \overset{\mathbb{P}}{\longrightarrow}  W^*_g( u^*)$. Using again Lemma~\ref{supas} for each fold, and then using the continuity of $G^\infty$, we get that
$$\left\| \widehat G^{\cross} -G^\infty_{W^*(u^*)}\right\|\overset{\mathbb{P}}{\longrightarrow}0 .$$
Then, by Remark~\ref{weight_vector},
$$\hat u^{cross}  \overset{\mathbb{P}}{\longrightarrow} u^\infty_{W^*(u^*)} ,$$
where $\hat u^{cross} =\mathcal{I}\left(\widehat G^{\cross}  \right) $ is the step-up threshold of $\crossADDOW$. Now, by using the asymptotic counterpart of Remark~\ref{MWBH<wbh}, we have that $u^\infty_{W^*(u^*)}=u^*$. Then we can proceed as in the proofs of Lemma~\ref{BOUMH} and Theorems~\ref{thm_fdr} and~\ref{thm_pow} to finally obtain that
$$ \lim_{m\to\infty} \FDR\left( \crossADDOW \right) =  \frac{H^\infty_{W^*}(u^*)}{u^*}$$
and
$$\lim_{m\to\infty} \Pow\left( \crossADDOW \right) = P^\infty_{W^*}(u^*)  ,$$
which concludes.


\begin{remark}
To keep notation light, we did not introduce the generalization of Section~\ref{subsection_notation} but we could use it to get a Theorem with a generalized $\crossADDOW$.
\end{remark}
\begin{remark}
It is easy to see that $\crossADDOW$ is a WBH procedure with $G\times F$ groups when using the expression given by Equation~\eqref{eq:Gcross}.
\end{remark}
\end{version3}

\end{document}